\documentclass[11pt]{article}
\usepackage{amsmath} \usepackage{amssymb} \usepackage{latexsym} \usepackage{enumitem} \usepackage{mathrsfs} \usepackage{comment}
\usepackage{color}
\usepackage{accents}
\usepackage{bbm}

\usepackage[colorlinks=true,urlcolor=blue,
citecolor=red,linkcolor=blue,linktocpage,pdfpagelabels,bookmarksnumbered,bookmarksopen]{hyperref}

\newtheorem{assumption}{Assumption}

\def\qed{ \ \vrule width.2cm height.2cm depth0cm\smallskip}
\newenvironment{proof}{\noindent {\bf Proof.\/}}{$\qed$\vskip 0.1in}

\newcommand{\ba}{\begin{array}}
\newcommand{\ea}{\end{array}}
\newcommand{\be}{\begin{equation}}
\newcommand{\ee}{\end{equation}}
\newcommand{\bea}{\begin{eqnarray}}
\newcommand{\eea}{\end{eqnarray}}
\newcommand{\beaa}{\begin{eqnarray*}}
\newcommand{\eeaa}{\end{eqnarray*}}

\def\dbE{\mathbb{E}}
\def\dbF{\mathbb{F}}

\def\dbL{\mathbb{L}}

\def\dbP{\mathbb{P}}
\def\dbR{\mathbb{R}}

\def\a{\alpha}
\def\b{\beta}
\def\g{\gamma}
\def\d{\delta}
\def\e{\varepsilon}

\def\k{\kappa}
\def\l{\lambda}
\def\m{\mu}

\def\si{\sigma}

\def\f{\varphi}
\def\th{\theta}

\def\G{\Gamma}
\def\D{\Delta}
\def\Th{\Theta}
\def\L{\Lambda}

\def\O{\Omega}
\def\cA{{\cal A}}

\def\cC{{\cal C}}

\def\cF{{\cal F}}
\def\cG{{\cal G}}
\def\cH{{\cal H}}

\def\cL{{\cal L}}
\def\cM{{\cal M}}

\def\cP{{\cal P}}

\def\cU{{\cal U}}

\def\ccD{\mathfrak{D}}

\def\no{\noindent}

\def\ms{\medskip}
\def\bs{\bigskip}
\def\q{\quad}
\def\qq{\qquad}

\def\pa{\partial}
\def\cd{\cdot}
\def\cds{\cdots}

\def\td{\nabla}

\def\bm{{\bf m}}

\def\tr{\hbox{\rm tr}}

\def\qed{ \hfill \vrule width.25cm height.25cm depth0cm\smallskip}

\newcommand{\basa}{\begin{assumption}}
\newcommand{\easa}{\end{assumption}}

\newcommand{\bas}{\begin{assum}}
\newcommand{\eas}{\end{assum}}

\def\pa{\partial}

 \def\cd{\cdot}
\def\cds{\cdots}

\def\supp{\hbox{\rm supp$\,$}}

\def\tr{\hbox{\rm tr$\,$}}

\def\dis{\displaystyle}

\def\wh{\widehat}

\def\1{{\bf 1}}

\def\:{\!:\!}
\def\reff{\eqref}
\def \proof{{\noindent \bf Proof.\quad}}

at 9pt

\definecolor{alp}{rgb}{0.0, 0.5, 0.0}

\newtheorem{thm}{Theorem}[section]
\newtheorem{lem}[thm]{Lemma}

\newtheorem{prop}[thm]{Proposition}
\newtheorem{rem}[thm]{Remark}
\newtheorem{eg}[thm]{Example}
\newtheorem{defn}[thm]{Definition}
\newtheorem{assum}[thm]{Assumption}

\begin{document}

\title{\bf  Second-order monotonicity conditions and mean field games with volatility control}
%Global Well-posedness of Master Equations for Extended Mean Field Games}
\author{Chenchen Mou\thanks{\noindent  Department of Mathematics,
City University of Hong Kong. E-mail: \href{mailto:chencmou@cityu.edu.hk}{chencmou@cityu.edu.hk}. This author is supported in part by Hong Kong RGC Grants GRF 11311422 and GRF 11303223.}, ~ Jianfeng Zhang\thanks{\noindent  Department of Mathematics,
University of Southern California. E-mail:
\href{mailto:jianfenz@usc.edu}{jianfenz@usc.edu}. This author is supported in part by NSF grant DMS-2205972.
} ~ and ~
Jianjun Zhou \thanks{\noindent  College of Science,
Northwest A\&F University. E-mail:
\href{mailto:jianfenz@usc.edu}{zhoujianjun@nwsuaf.edu.cn}.
}
}
\date{}
\maketitle

\begin{abstract}
In this manuscript we study the well-posedness of the master equations for mean field games with volatility control. This infinite dimensional PDE is nonlinear with respect to both the first and second-order derivatives of its solution. For standard mean field games with only drift control, it is well-known that certain monotonicity condition is essential for the uniqueness of mean field equilibria and for the global well-posedness of the master equations.  
To adapt to the current setting with volatility control, we propose a new notion called second-order monotonicity conditions. Surprisingly, the second-order Lasry-Lions monotonicity is equivalent to its standard (first-order) version, but such an equivalency fails for displacement monotonicity. When the Hamiltonian is separable and the data are Lasry-Lions monotone, we show that the Lasry-Lions monotonicity propagates and the master equation admits a unique classical solution. This is the first work for the  well-posedness, both local and global, of master equations when the volatility is controlled.
\end{abstract}

\no{\bf Keywords.}  Mean field games, volatility control, master equation, monotonicity conditions.

\ms
\no{\it 2020 AMS Mathematics subject classification:}  35Q89, 35R15, 49N80,  91A16, 60H30 %/93E20

\vfill\eject

%\tableofcontents

\section{Introduction}
\label{sect-Introduction}
\setcounter{equation}{0}

Mean Field Games (MFGs) represents a rapidly evolving field within mathematical game theory that focuses on the study of strategic decision-making among large populations of interacting agents. A general mathematical study of such games was initiated by Caines-Huang-Malham\'e \cite{HCM06} and Lasry-Lions \cite{LL07a} independently. For a comprehensive introduction to recent progress in this domain, we refer to the monographs by Carmona-Delarue \cite{CD1, CD2}, and the lecture notes by Cardaliaguet-Porretta \cite{CP}. MFG theory has proven to be extremely rich in applications, encompassing areas such as economics, engineering, finance, social science, and many others.

The MFG master equation was introduced by Lions in his lecture \cite{Lions} at Coll\`ege de France for the study of MFGs. It is an infinite dimensional nonlocal partial differential equation, which has the following form (for notational simplicity we assume the state space is $\dbR$):
\begin{equation}\label{master}
\begin{split}
\partial_t V+H(x,\partial_x V,\partial_{xx}V,\mu)+\mathcal{N}V=0,\,\, V(T,x,\mu)=G(x,\mu),\,\,\text{where}\qquad \\
\mathcal{N}V(t,x,\mu):=\int_{\mathbb R^d}\partial_{\tilde x\mu}V(t,x,\mu,\tilde x) \partial_\gamma H(\tilde x,\partial_xV(t,\tilde x,\mu),\partial_{xx}V(t,\tilde x,\mu),\mu)\\
+\partial_{\mu}V(t,x,\mu,\tilde x)\partial_z H(\tilde x,\partial_xV(t,\tilde x,\mu),\partial_{xx}V(t,\tilde x,\mu),\mu)\mu(d\tilde x),\qquad\,
\end{split}
\end{equation}
with the Hamiltonian $H:(x,z,\gamma,\mu)\in\mathbb R^3\times\mathcal{P}_2(\mathbb R)\mapsto \mathbb R$ and the terminal cost $G:(x,\mu)\in\mathbb R\times\mathcal{P}_2(\mathbb R)\mapsto \mathbb R$. Here $\partial_{\mu}$ stands for the Lions-derivative, also called Wasserstein derivative. Besides its independent interest for the construction of closed-loop Nash equilibria in MFGs, one of the main motivations for the introduction of master equation \eqref{master} is to provide quantitative convergence rates for closed-loop Nash equilibria of finite player games, when the number of involving players tends to infinity, see e.g. Cardaliaguet-Cirant-Porretta \cite{CCP0}, Cardaliaguet-Delarue-Lasry-Lions \cite{CDLL}, Carmona-Delarue \cite{CD2}, Delarue-Lacker-Ramanan \cite{DelarueLackerRamanan1,DelarueLackerRamanan2}. The master equation \eqref{master} corresponds to MFGs with both drift and volatility controls, so it is nonlinear in $\partial_xV$ and $\partial_{xx}V$. Recall that standard MFGs only involve drift control and thus their master equations are nonlinear in $\partial_{x}V$ but linear in $\partial_{xx}V$. To the best of our knowledge, the local and global well-posedness for the fully nonlinear master equation \eqref{master} remains completely open.

Mean field equilibria are typically characterized through the MFG system, which is a coupled forward backward system of PDEs or McKean-Vlasov SDEs. The solution $V$ of the master equation serves as the decoupling field for the MFG system, and on the opposite direction, the solution to the MFG system provides a representation for $V$ and hence leads to the well-posedness of the master equation. The idea is analogous to that of the following two-point boundary problem: for any $(t_0,x)\in [0,T]\times\dbR$,
\begin{equation}\label{ODE}
\dot x_t=b(x_t,y_t),\,\, x_{t_0}=x;\qq \dot y_t=f(x_t,y_t),\,\,y_T=g(x_T),
\end{equation}
where Cauchy-Lipschitz theory holds only in a small duration (depending on the Lipschitz constants of the data $b, f$ and $g$). Thus, in the small duration, we can find its decoupling field $v:[t_0,T]\times \dbR\to\dbR$ such that $y_t=v(t,x_t)$ for any $t\in [t_0,T]$, and $v(T,\cdot)=g(\cdot)$. To obtain a global solution of \eqref{ODE}, we need to glue the local ones to a global one, and the necessary ingredient for the gluing procedure is to control the Lipschitz constant of $v(t_0,\cdot)$ with respect to initial data $x$ for any $t_0\in[0,T]$. In our setting, note that both $\pa_xV$ and $\pa_{xx}V$ are nonlinearly involved in the Hamiltonian $H$ and hence in the corresponding MFG system. Therefore, similar to \eqref{ODE}, a priori Lipschitz estimates of $\pa_xV,\pa_{xx}V$ with respect to the measure variable $\mu$ are crucial to the global well-poesdness of the MFG system, which is the initial step towards to the global well-posedness of master equation \eqref{master}. For standard MFGs (with only drift controls), however, only the a priori Lipschitz estimate of $\pa_xV$ is required for its well-posedness.

Significant progress has been made on standard MFG master equations in recent years. There are various well-posedness results for master equations over short time duration, see e.g. Bensoussan-Yam \cite{BY}, Cardaliaguet-Cirant-Porretta \cite{CCP}, Carmona-Delarue \cite{CD2}, Gangbo-Swiech \cite{GS}.  As mentioned earlier, to obtain the global (in time) well-posedness results of standard MFG master equations, we need the a priori Lipschitz estimate of $\pa_xV$ with respect to $\mu$. It is now well-known that the key to ensuring such Lipschitz estimate is to show the propagation of certain monotonicity condition. For standard MFGs, these monotonicity conditions are typically classified into two main groups: the Lasry-Lions monotonicity and the displacement monotonicity. Initially, the Lasry-Lions monotonicity condition was widely used in the literature, 
see e.g. Bertucci \cite{B1}, Bertucci-Cecchin \cite{BC}, Bertucci-Lasry-Lions \cite{BLL1}, Cardaliaguet-Delarue-Lasry-Lions \cite{CDLL}, Cardaliaguet-Souganidis \cite{CarSou}, Carmona-Delarue \cite{CD2}, Chassagneux-Crisan-Delarue \cite{CCD}, Mou-Zhang \cite{MZ2}. Displacement monotonicity was later identified as an effective alternative condition,
see e.g. Bensoussan-Graber-Yam \cite{BGY2}, Gangbo-M\'esz\'aros \cite{GM}, Gangbo-M\'esz\'aros-Mou-Zhang \cite{GMMZ}. It is worth mentioning that other types of monotonicity conditions have recently been discovered to be useful for the global well-posedness of master equations, see e.g. Bansil-M\'esz\'aros \cite{BM}, Graber-M\'esz\'aros \cite{GraberM2}, Mou-Zhang \cite{MZ3}. Without any monotonicity conditions, master equations are globally ill-posed in general, see e.g. Mou-Zhang \cite{MZ1}, and their well-posedness can be shown before the blow-up of its characteristics, see Bertucci-Lasry-Lions \cite{BLL}. However, for potential MFGs, the global well-posedness of the corresponding master equations may still hold in certain weak sense, see Cecchin-Delarue \cite{CecchinDelarue1}.

Recall that the Lasry-Lions monotonicity condition can be formulated as follows:
%\begin{equation}\label{eq:LLmonint}
%\dbE\big[U(\xi_1,\cL_{\xi_1})+U(\xi_2,\cL_{\xi_2})-U(\xi_1,\cL_{\xi_2})-U(\xi_2,\cL_{\xi_1})\big]\geq 0,\,\,\text{for any random variables $\xi_1,\xi_2$},
%\end{equation}
%which is equivalent to the following differential form:
\begin{equation}\label{eq:1stLLmonint}
\dbE\big[\pa_{x\mu}U(\xi,\cL_{\xi},\tilde \xi)\eta\tilde\eta\big]\geq 0,
\end{equation}
for any square integrable random variables $(\xi,\eta)$ and their independent copy $(\tilde\xi,\tilde\eta)$. The displacement monotonicity condition can be expressed as:
%\begin{equation}\label{eq:dismonint}
%\dbE\big[(\pa_xU(\xi_1,\cL_{\xi_1})-\pa_xU(\xi_1,\cL_{\xi_1}))(\xi_1-\xi_2)\big]\geq 0,\,\,\text{for any random variables $\xi_1,\xi_2$},
%\end{equation}
%which is equivalent to the following differential form:
\begin{equation}\label{eq:1stdismonint}
\dbE\big[\pa_{x\mu}U(\xi,\cL_{\xi},\tilde \xi)\eta\tilde\eta +\pa_{xx}U(\xi,\cL_\xi)\eta^2\big]\geq 0. %,\,\,\text{for any random variables $\xi,\eta$}.
\end{equation}
Notice that the propagation of the above monotonicity conditions is sufficient to derive the a priori Lipschitz estimates for the first-order derivative $\pa_x V$ with respect to $\mu$ in standard MFGs. However, these conditions are in general not sufficient to derive the same estimates for $\pa_x V$ and the second-order derivative $\pa_{xx} V$ in MFGs with volatility control. 

One of the main contributions of this manuscript is the discovery of the second-order monotonicity  conditions, which would ensure the desired uniform Lipschitz estimates for $\pa_x V, \pa_{xx} V$ in the volatility control setting. For a function $U:\dbR\times\cP_2(\dbR)\to \dbR$, introduce\footnote{When extending to multidimensional setting with $(x, \mu)\in \dbR^d\times \cP_2(\dbR^d)$, $\ccD_x U$ and $\ccD_\mu U$ should be interpreted as bilinear forms on $(\dbR^d, \dbR^{d\times d})$.}: 
%Throughout the paper, for every
\begin{equation}\label{eq:Dxmu}
\left.\ba{c}
\dis\mathfrak{D}_x U(x,\mu):=\left(
  \begin{array}{cc}
    \partial_{xx}U(x,\mu) & \partial_{xxx}U(x,\mu) \\
   \partial_{xxx}U(x,\mu) &\partial_{xxxx}U(x,\mu) \\
  \end{array}
\right),\,\,\forall (x,\mu)\in\dbR\times\cP_2(\dbR),\ms\\
%\end{equation}
%\begin{equation}\label{eq:Amu}
\dis \mathfrak{D}_\mu U(x,\mu,\tilde x):=\left(
  \begin{array}{cc}
    \partial_{x\mu}U(x,\mu,\tilde x) & \partial_{\tilde{x}x\mu}U(x,\mu,\tilde x) \\
   \partial_{xx\mu}U(x,\mu,\tilde x) &\partial_{\tilde{x}xx\mu}U(x,\mu,\tilde x) \\
  \end{array}
\right),\,\,\forall (x,\mu,\tilde x)\in\dbR\times\cP_2(\dbR)\times\dbR.
\ea\right.
\end{equation}
We say $U$ satisfies the second-order Lasry-Lions monotonicity condition if, for any $\xi, \eta_1, \eta_2$,
\begin{equation}\label{eq:2ndLLmonint}
\dbE\Big[\big(\eta_1,\eta_2\big)\mathfrak{D}_\mu U(\xi,\cL_{\xi},\tilde\xi) \big(\tilde\eta_1,\tilde\eta_2\big)^\top\Big]\geq 0,
\end{equation}
and we say $U$ satisfies the second-order displacement monotonicity condition if
\begin{equation}\label{eq:2nddismonint}
\dbE\Big[\big(\eta_1,\eta_2\big)\mathfrak{D}_\mu U(\xi,\cL_{\xi},\tilde\xi) \big(\tilde\eta_1,\tilde \eta_2\big)^\top + \big(\eta_1,\eta_2\big)\mathfrak{D}_x U(\xi,\cL_{\xi}) \big(\eta_1,\eta_2\big)^\top\Big]\geq 0.
\end{equation}
It is evident that second-order monotonicity conditions \eqref{eq:2ndLLmonint} and \eqref{eq:2nddismonint} imply monotonicity conditions \eqref{eq:1stLLmonint} and \eqref{eq:1stdismonint}, respectively. It is somewhat surprising that, provided sufficient regularity of $U$, the second-order Lasry-Lions monotonicity  condition \reff{eq:2ndLLmonint} is actually equivalent to its first-order version \reff{eq:1stLLmonint}. However, such an equivalency fails for the displacement monotonicity between \reff{eq:1stdismonint} and \reff{eq:2nddismonint}. Due to this fact, this manuscript focuses on the study of global well-posedness of master equation \eqref{master} under the Lasry-Lions monotonicity condition. Similar to the drift control case, in order to use the Lasry-Lions monotonicity, we shall assume the Hamiltonian $H$ to be separable. We will study more general second-order monotonicity conditions, including the second-order displacement monotonicity,  and the global well-posedness of the master equation \reff{master} with general non-separable Hamiltonian in a subsequent work.
%general master other types of monotonicity conditions for future work. 

We show that, under Lasry-Lions monotonicity condition \eqref{eq:1stLLmonint} on our coefficients, the solution $V(t,\cd)$ of master equation \eqref{master} also satisfies \eqref{eq:1stLLmonint} for all $t\in [0,T]$, i.e. $V$ propagates Lasry-Lions monotonicity condition \eqref{eq:1stLLmonint}. We  follow the approach  in Gangbo-M\'esz\'aros-Mou-Zhang \cite{GMMZ}, demonstrating that if $V$ is a smooth solution of master equation \eqref{master}, then for any random variable $(\xi,\eta)$, there exist (characteristic) processes $(X,\delta X)$ starting at $(\xi,\eta)$ such that
\begin{equation*}
{d\over dt}\dbE\big[\pa_{x\mu}V(t,X_t,\cL_{X_t},\tilde X_t) \delta X_t\widetilde{\delta X_t}\big]\leq 0,
\end{equation*}
which implies the propagation of \reff{eq:1stLLmonint} immediately.
Then, with the equivalence between  \eqref{eq:1stLLmonint} and \eqref{eq:2ndLLmonint}, we also obtain the propagation of the second-order Lasry-Lions monotonicity condition \eqref{eq:2ndLLmonint} for $V$. We then show that  \eqref{eq:2ndLLmonint} for $V$ implies the a priori Lipschitz estimates of $\pa_xV,\pa_{xx}V$ with respect to $\mu$. Consequently, this results in the global well-posedness of master equation \eqref{master} by extending local solutions to a global one. As an important application of the well-posedness result, we also prove a stability result for master equation \reff{master}, which is quite natural but seems to be new even in the drift control case.

We should note that the local well-posedness of master equation \eqref{master} is also new, and is more involved than the drift control case. We establish it through the corresponding McKean-Vlasov FBSDE system. Unlike in standard MFGs, this FBSDE system consists of  one forward SDE for the characteristic process and three backward SDEs, corresponding to $V$, $\pa_x V$, and $\pa_{xx}V$, respectively. The key is a  pointwise representation formula for the Lions derivative $\pa_\mu V$, following the idea in Mou-Zhang \cite{MZ2} and Gangbo-M\'esz\'aros-Mou-Zhang \cite{GMMZ}.

Although our paper is the first one establishing well-posedness of master equations for MFGs with volatility control, there have been some works on MFGs with volatility control in the literature, see e.g. Barrasso-Touzi \cite{BT}, Djete \cite{Djete}, Jackson-Tangpi \cite{JT}, Lacker \cite{Lacker0}, Lacker-Webster \cite{LackerWebster},  mostly focusing on the existence of mean field equilibria or the convergence of the corresponding $N$-player games. In particular, interestingly  \cite{JT} obtained the uniqueness of mean field equilibrium, which is very closely related to the well-posedness of the master equation. The work \cite{JT} does not involve our second-order monotonicity conditions, instead, the uniqueness of mean field equilibrium is ensured by the (first-order) strict displacement monotonicity condition together with an additional structural condition involving certain convexity on the data with respect to the state and the control variables. The latter condition requires further understanding.
%Lagrangian.
% These conditions are both crucial for proving the uniqueness of mean field equilibria, and thus lead to the quantitative convergence result.

The rest of the paper is organized as follows. In Section \ref{sect-2ndmon} we introduce the second-order monotonicity conditions and prove the equivalency in the Lasry-Lions sense. In Section \ref{sect-MFG}  we present MFGs with volatility control, the MFG systems, and the connection with the master equation. In Section \ref{sect-main} we present the technical assumptions and the main results. The global wellposedness is proved in Section \ref{sect-global}, provided the local wellposedness, which includes the propagation of Lasry-Lions monotonicity and the uniform Lipschitz estimates. In Section \ref{sect-stability} we prove the stability result. Section \ref{sect-local} is devoted to the local well-posedness, including the pointwise representation formula for $\pa_\mu V$.  Finally, some technical calculations are reported in Appendix.

\section{Second-order monotonicity conditions}
\label{sect-2ndmon}
\setcounter{equation}{0}

In this section we introduce the basic setting and propose a new notion of second-order monotonicity condition, both in Lasry-Lions sense and in displacement sense, tailored to MFGs with volatility controls. To ease the presentation, throughout the paper we assume all the involved processes are one dimensional. We emphasize that all the results in this paper can be extended to the multiple dimensional setting straightforwardly, but with much heavier notation.

\subsection{The basic setting}
\label{sect-setting}
Fix a finite time horizon $T>0$, a filtered probability space $(\O,  \dbF, \dbP)$, and a Brownian motion $B$. We shall assume $\cF_0$ is rich enough to support any probability measure on $\dbR$ and $\dbF = \{\cF_t\}_{0\le t\le T}$ is generated by $\cF_0$ and $B$. As often in the MFG literature, let $(\tilde \O, \tilde\dbF, \tilde \dbP)$ and $(\bar\O, \bar \dbF, \bar \dbP)$ be copies of $(\O, \dbF, \dbP)$, and whenever needed, we shall consider their product space. For a random variable (or process) $\xi$ on $(\O, \dbF, \dbP)$, we may extend it to the product space in an obvious manner, still denoted as $\xi$, and we can easily construct independent  copies, denoted as $\tilde \xi$, $\bar \xi$. We abuse the notation $\dbE$ to denote all expectations $\dbE^\dbP$, $\dbE^{\dbP\times \tilde \dbP}$, $\dbE^{\dbP\times \tilde \dbP\times \bar\dbP}$.

For any $p\ge 1$ and $\si$-algebra $\cG\subset \cF_T$, let $\cP_p(\dbR)$ denote the set of probability measures $\mu$ on $\dbR$ with finite $p$-th moment, and $\dbL^p(\cG)$ the set of $\cG$-measurable random variables $\xi$ with finite $p$-th moment. Moreover, given $\mu\in \cP_p(\dbR)$, $\dbL^p(\cG;\mu)$ denote the subset of $\xi\in \dbL^p(\cG)$ such that $\cL_\xi = \mu$, where $\cL_\xi$ denotes the law of $\xi$ under $\dbP$.  The set $\dbL^p(\cG)$ is equipped with the $p$-th norm $\|\cd\|_p$ and $\cP_p(\dbR)$ is equipped with the $p$-Wasserstein distance $W_p$, and both are  complete metric spaces:
\beaa
\|\xi\|_p := \big(\dbE[|\xi|^p]\big)^{1\over p},\q    W_p(\mu_1,\mu_2):=\inf\big\{\big(\|\xi_1-\xi_2\|_p: \text{ for all  $\xi_i\in \dbL^p(\cF_T;\mu_i)$, $i=1,2$}\big\}.
\eeaa

For a $W_2$-continuous function $U: \mathcal{P}_2(\dbR) \to \mathbb R$, %its linear functional derivative ${\d\over \d\mu} U$ and 
its Lions-derivative, also called Wasserstein derivative, is a function $\pa_\mu U: (\mu, \tilde x)\in \mathcal{P}_2(\dbR)\times \mathbb R\to \mathbb R$ satisfying:
% but defined in Gateaux's sense and in Frechet's sense, respectively:
\begin{eqnarray}
\label{pamu}
\left.\ba{lll}
%\dis U(\mu + \e (\nu-\mu)) - U(\mu) = \e \int_\dbR {\d\over \d\mu} U(\mu, \tilde x) (\nu-\mu)(d\tilde x)  + o(\e),\q \forall \mu, \nu \in \cP_2(\dbR), \e>0;\\
\dis U(\mathcal{L}_{\xi +  \eta}) - U(\cL_\xi) = \mathbb E\big[ \partial_\mu U(\mu, \xi) \eta  \big] + o(\|\eta\|_2), \ \forall\xi, \eta\in\mathbb L^2(\mathcal{F}_T).
\ea\right.
\end{eqnarray}
%Under appropriate conditions, it is known that (cf. \cite{}): $\pa_\mu U(\mu, \tilde x) = \pa_{\tilde x} {\d\over \d\mu} U(\mu, \tilde x)$.
The value function for the mean field game will be defined on the state space
\beaa
\Theta:= [0, T]\times \mathbb R \times \mathcal{P}_2(\dbR).
\eeaa
For $U: (t,x,\mu)\in \Th\to \dbR$, we define $\pa_t U, \pa_x U: \Th\to\dbR$ in the standard sense and $\pa_\mu U: (t,x,\mu, \tilde x)\in\Th \times \dbR \to \dbR$ as in \reff{pamu}. We may define higher order derivatives, including $\pa_{\tilde x}\pa_\mu U$, in the obvious way. Tailored to our study of master equations, we introduce the following spaces of smooth functions. We first note that, since the Hamiltonian $H$ of the MFG has variables $(x, z, \g, \mu)$, we extend the state space to $\Th^d:= [0, T]\times \dbR^d \times \cP_2(\dbR)$. Consider functions $U: (t,y,\mu) \in \Th^d \to \dbR$,
%to general $\dbR^d$ with its variable denoted as $y\in \dbR^d$,
 and let $\pa^{(k)}_y U$ denote the $k$-th order derivative of $U$ with respect to $y$. 

\begin{defn}
\label{defn-space}
Let $U: (t,y, \mu)\in \Th^d \to \dbR$ be a continuous function.

\ms
\no(i) We say $U\in \cC^{1,2}(\Th^d)$ if all the derivatives with total parabolic order up to $2$ exist and are continuous:  $\pa_t U, \pa_y U, \pa_{yy} U, \pa_\mu U, \pa_{y\mu} U, \pa_{\tilde x\mu} U$.

\ms
\no(ii) Let $\cC^{1,4}(\Th^d)$ denote the subspace of $U\in \cC^{1,2}(\Th^d)$ such that $\pa_y U, \pa_{yy} U \in \cC^{1,2}(\Th^d)$. That is, the following derivatives exist and are continuous: 
\beaa
\pa^{(k)}_y U,~ 0\le k\le 4;\q \pa^{(k)}_y \pa_\mu U,~0\le k\le 3; \q \pa^{(k)}_y \pa_t U, ~\pa^{(k)}_y \pa_{\tilde x\mu} U, ~0\le k\le 2.
\eeaa

\no(iii) Let $\cC^{2,6}(\Th^d)$ denote the subspace of $U\in \cC^{1,4}(\Th^d)$ such that 
%$\pa_y U, \pa_{yy} U, \pa_\mu U, \pa_{\tilde x\mu} U$ are in $\cC^{1,4}$. That is, 
$\pa^{(k)}_{y} U \in \cC^{1,2}(\Th^d)$, $k=0,\cds, 4$; and $\pa^{(k)}_{y}\pa_\mu U,~ \pa^{(k)}_{y}\pa_{\tilde x\mu} U \in \cC^{1,2}(\Th^{d+1})$, $k=0, 1, 2$.

\ms
\no(iv) For all the above spaces $\cC^{k,n}(\Th^d)$, let $\cC^{k,n}_b(\Th^d)$ denote the subspace such that all the involved derivatives are uniformly bounded. We emphasize that we do not require $U$ itself to be bounded.
\end{defn}
When $U$ is independent of $t$, we denote $\cC^{k,n}(\dbR^d\times \cP_2(\dbR))$ and $\cC^{k,n}_b(\dbR^d\times \cP_2(\dbR))$  in the obvious sense.
The space $\cC^{1,2}$ is for the standard classical solution. Our second-order monotonicity condition requires the regularity in $\cC^{1,4}$. We raise the regularity further to $\cC^{2,6}$ for technical convenience, which can possibly be weakened, see Remark \ref{rem-stability} below.

One important property of functions $U\in \cC^{1,2}(\Th)$  is the following It\^o's formula, see \cite{BLPR, CD2, CCD}. For $i=1,2$, let $d X^i_t := b^i_t dt + \si^i_t dB_t,$ where $b^i, \si^i:[0,T]\times\Omega\to\mathbb R$  are $\dbF$-progressively measurable and bounded (for simplicity), and $m_t:= \cL_{X^2_t}$, then we have
\bea
\label{Ito}
\left.\ba{c}
\dis d U(t, X^1_t, m_t) =  \Big[\pa_t U + \pa_x U b^1_t + \frac{1}{2} \pa_{xx} U |\si_t^1|^2\Big](t, X^1_t, m_t) dt +\pa_xU(t,X^1_t, m_t)\si_t^1dB_t\ms\\
\dis  +\dbE_{\cF_t}\Big[\pa_\mu U(t,X^1_t,m_t,\tilde X^2_t) \tilde b^{2}_t+ \frac{1}{2} \pa_{\tilde x\mu} U(t, X^1_t, m_t, \tilde X^2_t)|\tilde \si_t^2|^2\Big]dt,
\ea\right.
\eea
where $(\tilde X^2, \tilde b^2, \tilde \si^2)$ is an independent copy of $(X^2, b^2, \si^2)$, and $\dbE_{\cF_t}$ is the conditional expectation.

\subsection{Second-order Lasry-Lions monotonicity condition}
\label{subsect-2ndLL}

We first recall the standard  Lasry-Lions monotonicity condition for a function $U:\dbR\times\cP_2(\dbR)\to\dbR$:
\begin{equation}\label{eq:LLmon}
\mathbb E\big[U(\xi_1,\mathcal{L}_{\xi_1})+U(\xi_2,\mathcal{L}_{\xi_2})-U(\xi_1,\mathcal{L}_{\xi_2})-U(\xi_2,\mathcal{L}_{\xi_1})\big]\geq 0\quad\forall \xi_1,\xi_2\in\dbL^2(\cF_T).
\end{equation}
When $U\in\cC_b^{1,2}(\dbR\times\cP_2(\dbR))$, this is equivalent to its differential form (cf. \cite[Remark 2.4]{GMMZ}):
\begin{equation}
\label{LLmon1}
\cM_{LL}^1 U(\xi, \eta):= \dbE\big[\pa_{x\mu}U(\xi,\cL_{\xi},\tilde\xi)\eta\tilde\eta\big]\geq 0,\q\forall \xi,\eta\in\dbL^2(\cF_T).
\end{equation}

To motivate the second-order condition, we consider the potential game case. That is, $U(x, \mu) = {\d\over \d\mu} \cU(\mu, x)$ for some $\cU: \cP_2(\dbR)\to \dbR$, where $ {\d\over \d\mu} \cU: \cP_2(\dbR) \times \dbR\to\dbR$ is the linear functional derivative:
% defined in Gateaux's sense: for $\e\downarrow 0$, 
\begin{equation}\label{eq:deltamu}
\cU(\nu)-\cU(\mu)= \int_0^1 \int_{\dbR}\frac{\delta }{\delta\mu}\cU\big(\mu + \th (\nu-\mu), x\big)(\nu-\mu)(dx)d\th,\q \forall\nu, \mu\in\cP_2(\dbR). 
\end{equation}
Then the above monotonicity condition is equivalent to the convexity of the function $\cU$, see e.g. \cite[Remark 5.75]{CD1}. Indeed, noting that ${\d\over \d\mu} {\d\over \d\mu} \cU(\mu, x, \tilde x) = {\d\over \d\mu} U(x, \mu, \tilde x)$,  the convexity implies:
\bea
\label{convex}
0 \le {d^2 \over d \e^2} \cU\big(\mu + \e(\nu-\mu)\big)\big|_{\e=0} = \int_{\dbR^2}  {\d\over \d\mu} U(x, \mu, \tilde x) (\nu-\mu)(dx)(\nu-\mu)(d\tilde x).
\eea
Given $\xi \in \dbL^2(\cF_0)$ and a Lipschitz continuous function $\phi: \dbR\to \dbR$, consider the random ODE:
\bea
\label{Xdrift}
X_t = \xi + \int_0^t \phi(X_s) ds.
\eea
Then, denoting $\mu:= \cL_\xi$, $\nu = \mu_t := \cL_{X_t}$, $\eta := \phi(\xi)$, we can show that (see Appendix):
\bea
\label{LL1}
 \lim_{t\downarrow 0} {1\over t^2} \int_{\dbR^2}  {\d\over \d\mu} U(x, \mu, \tilde x) (\mu_t-\mu)(dx)(\mu_t-\mu)(d\tilde x) = \cM_{LL}^1 U(\xi, \eta).%\dbE\big[\pa_{x\mu}U(\xi,\cL_{\xi},\tilde\xi)\eta\tilde\eta\big].
\eea
%Here we used the fact that  $\pa_x {\d\over \d\mu}U(x, \mu) = \pa_\mu U(\mu, x)$. 
Then \reff{LLmon1} follows from \reff{convex} for this $(\xi, \eta)$. 

To adapt to the volatility control, we change \reff{Xdrift} to an SDE:
\bea
\label{Xvolatility}
X_t = \xi + \int_0^t \phi_1(X_s) ds + \int_0^t \phi_2(X_s) dB_s,
\eea
where $\phi_1, \phi_2$ are Lipschitz continuous functions.
Denote $\mu:= \cL_\xi$, $\mu_t := \cL_{X_t}$, $\eta_1 := \phi_1(\xi)$, $\eta_2 :={1\over 2} |\phi_2(\xi)|^2$. Similarly to \reff{LL1}, we will show in Appendix that, when $U\in \cC^{1,4}_b(\dbR\times \cP_2(\dbR))$,
\bea
\label{LL2}
 \lim_{t\downarrow 0} {1\over t^2} \int_{\dbR^2}  {\d\over \d\mu} U(x, \mu, \tilde x) (\mu_t-\mu)(dx)(\mu_t-\mu)(d\tilde x) = \dbE\Big[ (\eta_1, \eta_2) \ccD_{\mu}U(\xi,\cL_{\xi},\tilde\xi)(\tilde\eta_1, \tilde \eta_2)^\top\Big].
\eea
This motivates us to introduce the following notion of second-order Lasry-Lions monotonicity, which will be crucial for a priori estimates for the value function of MFG with volatility control. 
\begin{defn}
\label{defn-2LL}
We say $U\in \cC_b^{1,4}(\dbR\times \cP_2(\dbR))$ satisfies the second-order Lasry-Lions monotonicity condition if, recalling \reff{eq:Dxmu}, 
\begin{equation}
\label{LLmon2}
\cM_{LL}^2 U(\xi, \eta_1, \eta_2) := \dbE\Big[\big(\eta_1,\eta_2\big)\mathfrak{D}_\mu U(\xi,\cL_{\xi},\tilde \xi)\big(\tilde \eta_1,\tilde\eta_2\big)^\top\Big]\geq 0,\q\forall\xi,\eta_1,\eta_2\in\dbL^2(\cF_T).
\end{equation}
\end{defn}
We remark that this is not restricted to the potential game case, and we do not require $\eta_2\ge 0$. 

%Recall \reff{eq:Dxmu}, 
By setting $\eta_2 \equiv 0$, clearly \reff{LLmon2} implies \reff{LLmon1}. It is somewhat surprising that the two conditions are actually equivalent.

\begin{prop}\label{prop:2nd}
Let $U\in \mathcal{C}_b^{1,4}(\mathbb R\times\mathcal{P}_2(\dbR))$. Then $U$ is Lasry-Lions monotone if and only if it satisfies the second-order Lasry-Lions monotonicity condition. %\eqref{LLmon2} holds for any random variables $\xi,\eta_1,\eta_2\in \dbL^{2}(\cF_{0}^1)$.
\end{prop}
\proof As said we shall only prove the only if direction. Assume \reff{LLmon1} holds true. 

{\it Step 1.} We first assume that $\xi\in\dbL^2(\cF_T)$ has positive density $\rho\in \cC^1(\dbR)$, $\eta_1 = \phi_1(\xi)$ for some Borel measurable function $\phi_1: \dbR\to \dbR$, and  $\eta_2 = \phi_2(\xi)$ for some $\phi_2\in \cC^1(\dbR)$ with compact support.  Then, by applying the integration parts formula repeatedly, we have
\begin{eqnarray*}
&&\dis \cM_{LL}^2 U(\xi, \eta_1, \eta_2)%\dbE\Big[\big(\eta_1,\eta_2\big)\mathfrak{D}_\mu U(\xi,\cL_{\xi},\tilde \xi)\big(\tilde\eta_1,\tilde\eta_2\big)^\top\Big]\\
=\int_{\dbR^2} \Big[\pa_{x\mu}U(x,\mu,\tilde x)\phi_1(x)\phi_1(\tilde x) + \pa_{x x\m}U(x,\mu,\tilde x)\phi_2(x)\phi_1(\tilde x)\\
&&\qq\qq +\pa_{ \tilde xx\mu}U(x,\mu,\tilde x)\phi_2(\tilde x)\phi_1(x)+\pa_{\tilde xxx \mu}U(x,\mu,\tilde x)\phi_2(\tilde x)\phi_2(x)\Big]\rho(x)\rho(\tilde x)dxd\tilde x\\
&&=\int_{\dbR^2}\pa_{x\mu}U(x,\mu,\tilde x)\Big[\phi_1(x)-\frac{\big(\phi_2(x)\rho(x)\big)'}{\rho(x)}\Big]\Big[\phi_1(\tilde x)-\frac{\big(\phi_2(\tilde x)\rho(\tilde x)\big)'}{\rho(\tilde x)}\Big]\rho(x)\rho(\tilde x)dxd\tilde x\\
&&=\dbE\big[\pa_{x\mu}U(\xi,\cL_{\xi},\tilde\xi)\eta\tilde \eta\big] = \cM_{LL}^1 U(\xi, \eta)\geq 0,
\end{eqnarray*}
where $\eta:=\phi_1(\xi)-\frac{(\phi_2\rho)'}{\rho}(\xi)$. In particular, since $\rho>0$ and $\phi_2$ has compact support, we see that $\frac{(\phi_2\rho)'}{\rho}$ is bounded and hence $\eta\in\dbL^2(\cF_T)$. 

{\it Step 2.} We next consider general $\eta_i$, but $\xi$ still has positive density $\rho\in \cC^1(\dbR)$. Note that 
\bea
\label{etacond}
\dbE\Big[\big(\eta_1,\eta_2\big)\mathfrak{D}_\mu U(\xi,\cL_{\xi},\tilde \xi)\big(\tilde \eta_1,\tilde\eta_2\big)^\top\Big] = \dbE\Big[\big(\zeta_1,\zeta_2\big)\mathfrak{D}_\mu U(\xi,\cL_{\xi},\tilde \xi)\big(\tilde \zeta_1,\tilde\zeta_2\big)^\top\Big],
\eea
where $\zeta_i := \dbE\big[\eta_i\big|\si(\xi)\big]$ is $\si(\xi)$-measurable, $i=1,2$. Then $\zeta_i = \phi_i(\xi)$ for some deterministic function $\phi_i$. Since $\xi$ has density, we can easily construct $\phi_2^{n} \in \cC^1(\dbR)$ with compact support such that $\lim_{n\to\infty}\dbE[|\zeta^{n}_2 -\zeta_2|^2] = 0$, where $\zeta^{n}_2:= \phi^{n}_2(\xi)$. By Step 1, $\cM_{LL}^2 U(\xi, \zeta_1, \zeta^{n}_2)\ge 0$. Send $n\to\infty$, we see that  $\cM_{LL}^2 U(\xi, \zeta_1, \zeta_2)\ge 0$. Then by \reff{etacond} we obtain \reff{LLmon2} for general $(\eta_1, \eta_2)$. 

{\it Step 3.} Finally, for general $\xi$ and $(\eta_1, \eta_2)$, one can easily construct $\xi^{n}$ satisfying the requirements in Step 2 and $\lim_{n\to\infty} \dbE[|\xi^n-\xi|^2]=0$. By Step 2, $\cM_{LL}^2 U(\xi^n, \eta_1, \eta_2)\ge 0$. Send $n\to \infty$, we obtain the desired monotonicity.
\qed

To help understand the above seemingly surprising equivalence, we remark that the condition \reff{LLmon1} requires certain underlying structure and that same structure implies the second-order  condition \reff{LLmon2}. To illustrate the idea, we consider the following simple example, and we refer to \cite[Section 3.4.2]{CD1} for more examples.

\begin{eg}
Let $U(x,\mu):=\varphi(x)\int_{\mathbb R}\psi(\tilde x)\mu(d\tilde x)$ for some functions $\varphi,\psi\in \cC_b^2(\mathbb R)$. Then 
\beaa
\pa_{x\mu}U(x,\mu,\tilde x)=\varphi'(x)\psi'(\tilde x),\q\mbox{and thus}\q \cM_{LL}^1 U(\xi, \eta) = \dbE[\f'(\xi)\eta] ~\dbE[\psi'(\xi)\eta].
\eeaa
Note that \eqref{LLmon1} holds true if and only if $\varphi'\equiv \psi'$. % (which is a structural assumption). 
Then, for any $\xi,\eta_1,\eta_2\in\dbL^2(\cF_T)$,
\begin{eqnarray*}
&&\dis \cM_{LL}^2 U(\xi, \eta_1, \eta_2) % \dbE\big[(\eta_1,\eta_2)\ccD_{\mu}U(\xi,\cL_{\xi},\tilde \xi)(\tilde\eta_1,\tilde\eta_2)^\top\big]\\
= \dbE\bigg[(\eta_1,\eta_2)\left(
  \begin{array}{cc}
    \varphi'(\xi) \varphi'(\tilde \xi) &     \varphi'(\xi) \varphi''(\tilde \xi)   \\
        \varphi''(\xi) \varphi'(\tilde \xi)  &     \varphi''(\xi) \varphi''(\tilde \xi)  \\
  \end{array}
\right)\left(
  \begin{array}{c}
    \tilde\eta_1  \\
    \tilde\eta_2 \\
  \end{array}
\right)\bigg]\\
&&\dis =\Big(\dbE\big[\varphi'(\xi)\eta_1\big]\Big)^2+2\dbE\big[\varphi'(\xi)\eta_1\big]\dbE\big[\varphi''(\xi)\eta_2\big]+\Big(\dbE\big[\varphi''(\xi)\eta_2\big]\Big)^2\\
&&\dis =\Big(\dbE\big[\varphi'(\xi)\eta_1+\varphi''(\xi)\eta_2\big]\Big)^2\geq 0.
\end{eqnarray*}
This shows that the structural assumption $\varphi'\equiv \psi'$, which implies \eqref{LLmon1}, also leads to \eqref{LLmon2}.
\end{eg}

\subsection{Second-order displacement monotonicity condition}
\label{sect-2nddis}
%\subsubsection{Displacement monotonicity condition}
We now turn to the displacement monotonicity of $U$. First recall its standard definition:
\begin{equation}\label{eq:dismon}
\mathbb E\Big[(\pa_xU(\xi_1,\mathcal{L}_{\xi_1})-\pa_xU(\xi_2,\mathcal{L}_{\xi_2})(\xi_1-\xi_2)\Big]\geq 0,\quad\forall \xi_1,\xi_2\in\dbL^2(\cF_T).
\end{equation}
When $U\in \cC_b^{1,2}(\dbR\times \cP_2(\dbR))$, this is equivalent to its differential form (cf. \cite[Remark 2.4]{GMMZ})):
\begin{eqnarray}\label{DPmon1}
\cM^1_{DP} U(\xi, \eta) := \dbE\Big[\pa_{x\mu}U(\xi,\cL_{\xi},\tilde \xi)\eta\tilde\eta + \pa_{xx}U(\xi,\cL_{\xi})\eta^2\Big] \ge 0, \quad\forall \xi,\eta\in\dbL^2(\cF_T).
\end{eqnarray}

To motivate the second-order displacement monotonicity, we again consider the potential game case: $U(x, \mu) := {\d\over \d\mu} \cU(\mu, x)$. It is known that, in this case $U$ is displacement monotone if and only if $\cU$ is displacement convex, that is, the mapping $\xi\in \dbL^2(\cF_T) \to \cU(\cL_\xi)$ is convex, see e.g. \cite[Lemma 5.72]{CD1}. Indeed, the displacement convexity implies that, given $\xi, \eta\in \dbL^2(\cF_0)$, the mapping $t\to \cU\big(\cL_{\xi + t \eta}\big)$ is convex and thus ${d^2 \over d t^2} \cU\big(\cL_{X_t}\big)\big|_{t=0}\ge 0$. We can show that (see Appendix):
\bea
\label{DP1}
{d^2 \over d t^2} \cU\big(\cL_{X_t}\big)\big|_{t=0}  = \cM^1_{DP} U(\xi, \eta),\q \mbox{where}\q X_t := \xi + t \eta.
\eea
Then the desired convexity implies \reff{DPmon1}. We remark that for displacement monotonicity one needs to consider the above $X$ in order to use the convexity, instead of the random ODE in \reff{Xdrift}. 

To adapt to the volatility control, we again change the $X$ in \reff{DP1} to a diffusion: 
\bea
\label{Xvolatility2}
X_t = \xi + \eta_1 t + \eta_2 B_t,\q \xi, \eta_1, \eta_2\in \dbL^2(\cF_0).
\eea
We shall show in Appendix that, for $U\in \cC^{1,4}_b(\dbR\times \cP_2(\dbR))$, and denoting  $\eta_2' := {1\over 2} |\eta_2|^2$,
\bea
\label{DP2}
{d^2 \over d t^2} \cU\big(\cL_{X_t}\big)\big|_{t=0}  = \dbE\Big[ (\eta_1, \eta_2') \ccD_{\mu}U(\xi,\cL_{\xi},\tilde\xi)(\tilde\eta_1, \tilde \eta_2')^\top + (\eta_1, \eta_2') \ccD_{x}U(\xi,\cL_{\xi})(\eta_1, \eta_2')^\top\Big].
\eea
This motivates us to introduce the following notion of second-order displacement monotonicity. 
\begin{defn}
\label{defn-2DP}
We say $U\in \cC_b^{1,4}(\dbR\times \cP_2(\dbR))$ satisfies the second-order displacement monotonicity condition if, for any $\xi, \eta_1, \eta_2\in \dbL^2(\cF_T)$, 
\begin{equation}
\label{DPmon2}
\!\!\!\!\!\cM_{DP}^2 U(\xi, \eta_1, \eta_2) := \dbE\Big[\big(\eta_1,\eta_2\big)\mathfrak{D}_\mu U(\xi,\cL_{\xi},\tilde \xi)\big(\tilde \eta_1,\tilde\eta_2\big)^\top+\big(\eta_1,\eta_2\big)\mathfrak{D}_x U(\xi,\cL_{\xi})\big(\eta_1, \eta_2\big)^\top\Big]\geq 0.
\end{equation}
\end{defn}
Again this is not restricted to the potential game case, and we do not require $\eta_2\ge 0$. 

By simply setting $\eta_2 = 0$, it is clear that \reff{DPmon2} implies \reff{DPmon1}. However, unlike the Lasry-Lions monotonicity, in general the displacement monotonicity  \reff{DPmon1} does not imply the second-order condition  \reff{DPmon2}, as we see in the following example, whose proof is postponed to Appendix. 

\begin{eg}
\label{eg-DP}
Let $U(x, \mu) = \sin(x) \int_\dbR \cos(\tilde x) \mu(d\tilde x) + {\k\over 2} x^2$, where $\k\ge 0$ is a constant. Then $U$ satisfies \reff{DPmon1} when $\k\ge 2$, but it does not satisfy \reff{DPmon2}. 
\end{eg}

We remark that the second-order monotonicity condition, both in Lasry-Lions sense and in displacement sense, will ensure some crucial a priori estimate for the value function of MFG with volatility control and hence lead to the global well-posedness of the corresponding master equation. Thus the key is to find appropriate conditions on the coefficients for the propagation of second-order monotonicity conditions. Due to the equivalency in Proposition \ref{prop:2nd}, in the Lasry-Lions case it is sufficient to propagate the standard monotonicity condition \eqref{LLmon1}, while in the displacement case we have to propagate the second-order condition \reff{DPmon2}, which involves the fourth-order derivatives and is much more involved. In this paper we shall focus on the main ideas and thus restrict to the Lasry-Lions case. We will investigate the global well-posedness of general MFGs with volatility control under more general monotonicity conditions, including second-order displacement monotonicity, in a subsequent work.

\section{Mean field games with  volatility control}
\label{sect-MFG}
\setcounter{equation}{0}
In this section we introduce MFGs with volatility control, the MFG systems, and their connections with the master equation \reff{master} in a heuristic way. 

Consider the setting in Section \ref{sect-setting}. Our MFG involves the following data:
$$
b, \sigma, F: (x, \mu, a)\in \mathbb{R}\times \mathcal{P}_2(\mathbb R)\times  A \rightarrow \mathbb{R}, \qq G:\mathbb{R}\times \mathcal{P}_2(\mathbb R)\rightarrow \mathbb{R},
$$
where $A\subset \dbR$ is the control set; $b, \si$  are bounded and uniformly Lipschitz continuous in $(x, \mu, a)$; $F$ and $G$ are continuous and bounded; and $\si\ge 0$. In particular, the volatility $\si$ may depend on the control $a$.
Let $\cA$ denote the set of closed loop admissible controls: $\a:  (t, x)\in [0, T]\times \mathbb{R}\to  A$ which is continuous in $t$ and uniformly Lipschitz continuous in $x$.\footnote{We may consider admissible controls in the more general form $\a(t,x,\mu)$ as well. However, since $\mu_t = \cL_{X_t}$ will be deterministic, all the results in the paper will remain the same.}  As usual, the Hamiltonian $H$ is defined as follows, which in particular is increasing in $\g$ and concave in $(z,\g)$:
\begin{eqnarray}\label{Hamiltonian}
H(x, z, \gamma,\mu):=\inf_{a\in A}\left[\frac{1}{2}\sigma^2(x,\mu, a)\gamma+b(x,\mu, a)z+F(x,\mu,a)\right].
\end{eqnarray}
%The main feature here is that the volatility $\si$ may depend on $a$, and thus $H$ is nonlinear in $\g$. 

We now consider the MFG on $[0, T]$, with initial time $0$. Fix $\mu\in \cP_2(\dbR)$ and $\xi \in \dbL^2(\cF_0, \mu)$. Given a flow of measures $\bm = \{m_t\}_{0\le t\le T} \subset \cP_2(\dbR)$, consider the following controlled SDE: 
%on $[0, T]$ has a unique solution:
\begin{eqnarray}\label{Xxia}
 X^{\bm; \xi,\alpha}_t=\xi+\int^t_{0}b\big(X^{\bm;\xi,\alpha}_s, m_s, \alpha(s,X^{\bm;\xi,\alpha}_s)\big)ds
                            +\int^t_{0}\sigma\big(X^{\bm;\xi,\alpha}_s, m_s, \alpha(s,X^{\bm;\xi,\alpha}_s)\big)dB_s,
\end{eqnarray}
and the following standard stochastic control problem :
\begin{eqnarray}\label{Xxa'}
\left.\ba{c}
\dis u({\bf m}; 0, x) := \inf_{\a\in \cA} J({\bf m}; x,\alpha),\q\mbox{where}\\
\dis J({\bf m}; x,\alpha):= \mathbb{E}\left[G(X_T^{{\bf m}; x,\alpha},m_T)
       +\int^T_{0}F\big(X_t^{{\bf m}; x,\alpha},m_t,\alpha(t, X_t^{{\bf m}; x,\alpha})\big)dt\right].
       \ea\right.                   
\end{eqnarray}

\begin{defn}
\label{MFE}
          We say $\bm$ is an equilibrium measure of the MFG \reff{Xxia}-\reff{Xxa'} at $(0,\mu)$ if there exists $\alpha^*\in \mathcal{A}$, called a mean field equilibrium (MFE, for short) of the MFG, such that
          $$
   \bm =    \cL_{X^{\bm; \xi,\a^*}},\q\mbox{and}\q    u( \bm; 0, x) = J(\bm; x,\alpha^*),~ \mbox{for}\ \mu\mbox{-a.e.}  \ x\in \mathbb{R}. %,\q\mbox{where}\q \bm^*_t := \cL_{X^{\xi,\a^*}_t}.
          $$    
\end{defn}

%\subsection{MFG systems, master equations, and the stability}

Equilibrium measures $\bm$ are typically characterized through MFG systems. 
Clearly $u(t,x):= u(\bm; t, x)$ satisfies the following fully nonlinear HJB equation with parameter $\bm$:
\bea
\label{HJB}
\pa_t u(t, x)  + H\big(x,  \pa_x u(t, x), \pa_{xx} u(t, x), m_t \big)=0,\q  u(T, x) = G(x, m_T). 
\eea
On the other hand, assuming $u$ is a classical solution of \reff{HJB}, by \reff{Hamiltonian} we have
\bea
\label{Hoptimal}
\left.\ba{c}
\dis b(x, m_t, \a^*(t,x)) =\pa_z\cH_u(t,x,m_t),\q  {1\over 2}\si^2(x, m_t, \a^*(t,x)) = \pa_\g\cH_u(t,x,m_t),\\
\dis \mbox{where}\q \pa_z\cH_u(t,x, m_t) :=\pa_zH\big(x,  \pa_x u(t, x), \pa_{xx} u(t, x), m_t \big), ~\mbox{and similarly for $\pa_\g \cH_u$}.
%\big(x,  \pa_x u( t, x), \pa_{xx} u(t, x), m_t \big);\\
%\dis {1\over 2}\si^2(x, m_t, \a^*(t,x)) = \pa_\g H\big(x,  \pa_x u(t, x), \pa_{xx} u(t, x), m_t \big).
\ea\right.
\eea
Then,  by \reff{Xxia} we see that $\bm$ is a weak solution of the following Fokker-Planck equation:
\bea
\label{Fokker-Planck}
\left.\ba{c}
\dis \pa_t m(t,x) = \pa_{xx}\big(\pa_\g\cH_u(t, x, m_t)  m(t,x) \big) - \pa_x \Big(\pa_z\cH_u(t,x, m_t) m(t,x) \Big) =0,\q m_0 = \mu.
%\dis \mbox{where}\q \wh\f_u(t,x, m_t) :=\f \big(x,  \pa_x u(t, x), \pa_{xx} u(t, x), m_t \big),\q \f=\pa_z H, \pa_\g H.
\ea\right.
\eea
We emphasize that, although the $\cH_u$ above depends on $\bm$, such dependence is globally on the whole measure, not locally on the density $m(t,x)$. The coupled forward-backward PDE system \reff{HJB}-\reff{Fokker-Planck} is known as the MFG system. 

 We next derive the MFG system in the probabilistic form, namely a system of forward-backward SDEs. This issue is actually more subtle than the drift control case. Given $\xi\in \dbL^2(\cF_0,\mu)$ and assuming $\bm$ is an equilibrium measure, by \reff{Hoptimal} and \reff{Fokker-Planck} we introduce:
 \bea
 \label{XxiwhH}
 X^\xi_t = \xi + \int_0^t \pa_z\cH_u(s, X^\xi_s, m_s)ds + \int_0^t \sqrt{2 \pa_\g\cH_u(s, X^\xi_s, m_s)}dB_s,
 \eea
so that $\bm = \cL_{X^\xi}$. %Here we assume $\si\ge 0$, and note that by \reff{Hamiltonian} we have $\pa_\g H\ge 0$. 
 Inspired by the drift control case,  one would set 
 \bea
 \label{Yxi0}
 Y^{\xi,0}_t := u(t, X^{\xi}_t), ~\mbox{and thus, by It\^{o}'s formula},~  Z^{\xi,0}_t := \pa_x u(t, X^{\xi}_t) \sqrt{2\pa_\g\cH_u(t, X^\xi_t, m_t)}.
 \eea
  However, in the volatility control case, since $\pa_z\cH_u$ and $\pa_\g\cH_u$ involve $\pa_{xx} u$, which cannot be expressed in terms of $(Y^{\xi,0}, Z^{\xi,0})$, and thus \reff{Yxi0} does not lead to a self-contained FBSDE. In light of this, we set
 \bea
 \label{Yxi1}
 Y^{\xi,1}_t := \pa_x u(t, X^{\xi}_t),\qq Z^{\xi,1}_t := \pa_{xx} u(t, X^{\xi}_t) \sqrt{2\pa_\g\cH_u(t, X^\xi_t, m_t)}.
\eea
We note again that at above $ \pa_\g\cH_u(t, X^\xi_t, m_t) = \pa_\g H\big(X^\xi_t,  \pa_x u(t, X^\xi_t), \pa_{xx} u(t, X^\xi_t), m_t \big)$. Then, under appropriate conditions on $\pa_\g H$, one may express $ \pa_{xx} u(t, X^\xi_t)$ in terms of $(Y^{\xi,1}_t, Z^{\xi,1}_t)$ and thus obtain a self-contained FBSDE. However, 
%besides the extra condition on $\pa_\g H$ in order to get the above inverse functions, 
in this way the diffusion term in \reff{XxiwhH} would involve $Z^{\xi,1}_t$, and then the coupled FBSDE is well known to be hard to solve, see e.g. \cite[Chapter 8]{Zhang}. To get around of this difficulty, we introduce further that
 \bea
 \label{Yxi2}
 Y^{\xi,2}_t = \pa_{xx} u(t, X^{\xi}_t),\qq Z^{\xi,2}_t = \pa_{xxx} u(t, X^{\xi}_t)\sqrt{2\pa_\g\cH_u(t, X^\xi_t, m_t)}.
\eea
Then $(\pa_x u, \pa_{xx} u) = (Y^{\xi,1}_t, Y^{\xi,2}_t)$. By applying the standard It\^{o}'s formula, one may derive directly from \reff{HJB} and \reff{XxiwhH} the following self-contained coupled McKean-Valsov FBSDE, whose forward diffusion term does not depend on $Z$ anymore\footnote{This trick does not help for standard fully nonlinear PDEs. The introduction of $Y^{\xi,2}$ requires higher order regularity of $u$ in $x$, which is typically the main target in the PDE theory. In the mean field case, however, our main focus is the regularity of the value function in $\mu$. We thus consider the regularity in $x$ as given, see Assumption \ref{assum-regu} below. Then this FBSDE system is very helpful for our study of the well-posedness of the master equation. }:  denoting $\Xi^\xi_t := (X^\xi_t, Y^{\xi,1}_t, Y^{\xi,2}_t, \cL_{X^\xi_t})$, 
\begin{eqnarray}\label{FBSDE}
\begin{cases}
                         \dis  X^{\xi}_t= \xi + \int_0^t H_1(\Xi^\xi_s)ds +\int_0^tH_2(\Xi^\xi_s)dB_s,\\
                         \dis  {Y}^{\xi,0}_t=G(X^{\xi}_T, \cL_{X^\xi_T}) +\int_t^T H_3(\Xi^\xi_s)ds - \int_t^T {Z}^{\xi,0}_sdB_s,\\
                        \dis {Y}^{\xi,1}_t= \partial_{x}G(X^{\xi}_T, \cL_{X^\xi_T}) +\int_t^T \partial_xH( \Xi^\xi_s)ds -\int_t^T {Z}^{\xi,1}_sdB_s,\\
                          \dis  {Y}^{\xi,2}_t
                          =\partial_{xx}G(X^{\xi}_T, \cL_{X^\xi_T}) +\int_t^T  H_4(\Xi^\xi_s, {Z}^{\xi,2}_s) ds -\int_t^T {Z}^{\xi,2}_sdB_s,                          \end{cases}
\end{eqnarray}
where, for $\Xi := (x, z, \g, \mu) = (x, y_1, y_2,\mu)$ (we shall use the notation $\pa_zH, \pa_\g H$ to denote $\pa_{y_1}H, \pa_{y_2} H$), 
\bea
\label{cH}
\left.\ba{lll}
\dis H_1(\Xi) := \pa_z H(\Xi),\q H_2(\Xi) := \sqrt{2\pa_\g H(\Xi)};\\
\dis H_3(\Xi) := H( \Xi) - \pa_z H(\Xi) y_1 - \pa_\g H(\Xi) y_2;\\
\dis H_4(\Xi, z_2) :=\partial_{xx}H(\Xi) +2\partial_{xz}H(\Xi) y_2 +\partial_{zz}H(\Xi)|y_2|^2    \\
\dis\q +2\big[\partial_{x\gamma}H(\Xi) + \partial_{z\gamma}H(\Xi)y_2\big] z_2\slash \sqrt{2\partial_{\gamma}H(\Xi)}+\partial_{\gamma\gamma}H(\Xi)|z_2|^2\slash \big(2\partial_{\gamma}H(\Xi)\big).
\ea\right.
\eea
This coupled FBSDE system is also called the MFG system. We remark that the equation for $Y^{\xi,0}$ is decoupled from the rest and thus the coupled system is mainly for $(X^\xi, Y^{\xi,1}, Y^{\xi,2})$.
 
When there is a unique equilibrium measure $\bm$ at $(0,\mu)$, or equivalently when the MFG system \reff{HJB}-\reff{Fokker-Planck} or \reff{FBSDE} has a unique solution, we may define the value of the MFG:
\begin{eqnarray}\label{MFGV}
   V(0,x,\mu):=u(\bm; 0, x).
\end{eqnarray}
Similarly we may define the dynamic value function $V(t,x,\mu)$. Then $V$ is associated to the master equation \reff{master}, and it serves as the decoupling field of the  MFG systems:
\bea
\label{decoupling}
u(t, x) = V(t, x, m_t), \qq  Y^{\xi,k}_t =  \pa^{(k)}_x V(t, X^\xi_t, \cL_{X^\xi_t}),~k=0,1,2. 
\eea

\section{The main results}
\label{sect-main}
\setcounter{equation}{0}
In this section we  shall present the required assumptions and state the main  well-posedness results.  All proofs are postponed to the remaining sections. We note that the order of the proofs may not follow the order of the theorems, however, there will be no danger of cycle proofs.

\subsection{The stability and the local well-posedness}
\label{subsect-local}
We first specify some technical conditions. 
For any $R>0$, denote
\bea
\label{DR}
\left.\ba{c}
\dis D_R:= \big\{(z,\g)\in \dbR^2: |z|\le R, |\g|\le R\big\}.
\ea\right.
\eea
For a function $U$ on an appropriate domain, let $\|U\|$ denote its uniform norm, and $\|U\|_R$  the supremum of $U$ when restricting the variables $(z, \g)$ to $D_R$.

\begin{assum}\label{assum-regHG} (i) $G\in \mathcal{C}^{2,6}_b(\dbR\times\mathcal{P}_2(\dbR))$ and there exists a constant $L_G>0$ such that
\bea
\label{LG}
%\|\mathfrak{D}_xG\|\leq L^{G}_0,\q \|\mathfrak{D}_\mu G\|\leq L^{G}_1.
\|\pa^{(k)}_x G\|, \|\pa^{(l)}_x \pa_\mu G\|,  \|\pa^{(l)}_x \pa_{\tilde x\mu} G\|\le L_G,\q ~ k=1,2,3,4;~ l=1,2.
\eea

%\ms
 \no(ii) $H\in \mathcal{C}^{2,6}(\dbR^3\times\mathcal{P}_2(\dbR))$ and, for $\forall R>0$, all the involved derivatives are bounded when restricting $(z,\g)$ to $D_R$, and there exists a constant $L_H(R)>0$ such that 
\beaa
\|\pa^{(k)}_{(x,z,\g)} H\|_R,  \|\pa^{(l)}_{(x,z,\g)} \pa_\mu H\|_R, \|\pa^{(l)}_{(x,z,\g)} \pa_{\tilde x\mu} H\|_R \le L_H(R),\q k=1,2,3,4; ~l=1,2.
\eeaa 

%\ms
\no (iii)  $\si \ge \sqrt{2c_0}$ for some $c_0>0$, and consequently $\pa_\g H\ge c_0>0$.
\end{assum}

Next, given $\xi\in \dbL^2(\cF_0)$ and $\a\in \cA$, consider the McKean-Vlasov SDE:
\bea
\label{Xxia2}
X^{\xi,\a}_t = \xi + \int_0^t b(X^{\xi,\a}_s, \cL_{X^{\xi,\a}_s}, \a(s, X^{\xi,\a}_s)) ds +   \int_0^t \si(X^{\xi,\a}_s, \cL_{X^{\xi,\a}_s}, \a(s, X^{\xi,\a}_s)) dB_s.
\eea

\begin{assum}\label{assum-regu}
For any $\xi\in \dbL^2(\cF_0)$ and $\a\in \cA$, the HJB equation \reff{HJB} with parameter $\bm := \cL_{X^{\xi,\a}}$ has a unique classical solution $u \in C^{1,6}_b([0, T]\times \dbR)$, that is, the derivatives $\pa_t u$ and $\pa^{(k)}_x u$, $1\le k\le 6$, exist and are bounded and continuous. Moreover, there exists a constant $L^*_T>0$, possibly depending on $T$ but independent of $\xi, \a$, such that
\bea
\label{eq:L0u}
\|\pa^{(k)}_xu\|\leq L^*_T,\q k=1,2,3,4.
\eea
\end{assum}

\begin{rem}
\label{rem-PDE1}
(i) Under Assumption \ref{assum-regu}, for each  equilibrium measure $\bm$ of the MFG at $(0,\mu)$, there exists unique $u \in C^{1,6}_b([0, T]\times \dbR)$ such that $(\bm, u)$ solves the MFG system \reff{HJB}-\reff{Fokker-Planck}. Moreover, in this case the solutions to \reff{HJB}-\reff{Fokker-Planck} have one to one correspondence to the solutions of the MFG system \reff{FBSDE}, with $\bm = \cL_{X^\xi}$ and \reff{Yxi1}-\reff{Yxi2} holding true. 

\ms
\no (ii) To avoid solving the McKean-Vlasov SDE \reff{Xxia2}, one may strengthen Assumption \ref{assum-regu} slightly by assuming the required properties for all $\bm = \cL_X$, where $X_t := \xi + \int_0^t b_s ds + \int_0^t \si_s dB_s$, and the processes $b, \si$ here are bounded by the same bound of the coefficients $b, \si$ in \reff{Xxia2}.
\end{rem}

\begin{rem}
\label{rem-PDE2}
(i) Since the HJB equation \reff{HJB} is fully nonlinear, its classical solution theory is highly nontrivial. Nonetheless, given the uniform non-degeneracy condition in Assumption \ref{assum-regHG}-(iii), this is possible due to the well known Evans-Krylov theory. In this paper we try to focus on the regularity of $V$ in $\mu$, and thus do not discuss sufficient conditions for the regularity in $x$. We refer to the book \cite{Lieberman} and the references therein for the theory on fully nonlinear HJB equations. 

\ms
\no (ii) The estimate \eqref{eq:L0u} will serve as the a priori regularity estimate for the solution $V$ to the master equation \eqref{master} with respect to $x$. We permit the constant $L^*_T$  in \eqref{eq:L0u} to depend on $T$. This dependence has undesirable consequence in \cite{MZ3},  because the anti-monotonicity condition there relies on $L^*_T$ and then its dependence on $T$ will imply that the well-posedness holds true only for small $T$. That is why \cite{MZ3} made serious efforts to obtain the estimate \reff{eq:L0u} independent of $T$, under additional conditions, so as to obtain the global well-posedness of the master equation (in the drift control case). However, we emphasize that there is no such issue here, as there are no additional assumptions based on $L^*_T$ when we establish the global well-posedness.

\ms
\no(iii) Under Assumption \ref{assum-regHG}, in Assumption \ref{assum-regu} actually it suffices to assume $u\in C^{1,2}_b([0, T]\times \dbR)$. Then the higher regularity of $u$ will follow from the bootstrap arguments. 
\end{rem}

Our first result concerns the uniqueness and stability of the master equation.  
\begin{thm}
\label{thm-stability}
Under Assumptions  \ref{assum-regHG}  and \ref{assum-regu}, the master equation \reff{master} has at most one classical solution $V\in \cC^{1,4}_b(\Th)$. Moreover, when such an $V$ exists, the following hold.

\no (i)  For any $\mu\in \cP_2(\dbR)$, the MFG has a unique equilibrium measure $\bm$ at $(0,\mu)$.\footnote{The MFE $\a^*$ may not be unique in this case. For example, when the whole system is not controlled, then all controls $\a\in \cA$ can be viewed as MFEs.} Consequently, the MFG systems \reff{HJB}-\reff{Fokker-Planck} and \reff{FBSDE} have a unique solution,
% with $u\in \cC^3_b([0, T]\times \dbR)$ and $Y^{\xi,1}, Y^{\xi,2}, Z^{\xi,1}, Z^{\xi,2}$ bounded, 
and \reff{decoupling} holds true;

\no(ii) Assume there is another system $(b', \si', F', G')$ such that 
%$(b', \si')$ satisfy the conditions of $(b, \si)$ above and 
the corresponding Hamiltonian $H'$ also satisfies Assumptions  \ref{assum-regHG}  and \ref{assum-regu}. Then, for any equilibrium measure $\bm'$ (possibly non-unique) to this MFG at $(0,\mu)$, the corresponding $u'$ satisfies
% and the MFG system \reff{HJB}-\reff{Fokker-Planck} has a solution $(\bm', u')$ (not necessarily unique), then we have
\bea
 \label{stability1}
\left.\ba{c}
%\dis |\pa_x\D u(0,x)| + |\pa_{xx} \D u(0,x)|  \le C\big[\|\D H\|_{1, R} + \|\D H\|_{2, R}+\|\pa_x\D G\| + \|\pa_{xx}\D G\| \big],\\
\dis \sum_{k=0}^2 |\pa^{(k)}_x\D u(0,x)| \le C\sum_{k=0}^2 \big[\|\pa^{(k)}_{(x,z,\g)}\D H\|_{R} +\|\pa^{(k)}_x\D G\|  \big],\q\mbox{for $\mu$-a.e. $x$},\\
\dis  \mbox{where}\q \D u := u-u',~ \D H := H-H', ~\D G:= G-G', %~R := \max\big(\|\pa_x V\|, \|\pa_{xx} V\|, \|\pa_x u'\|, \|\pa_{xx} u'\|\big],
 \ea\right.
 \eea
 and the constant $C$ depends on $T$, $c_0$, $L^*_T$, and $L_H(R)$ with $R = L^*_T$.
 % = \max(\|\pa_x V\|, \|\pa_{xx} V\|, \|\pa_x u'\|, \|\pa_{xx} u'\|)$, and $\sum_{k=1}^4 \big[\|\pa^{(k)}_xV\| + \|\pa^{(k)}_xu'\|\big]$.  
 In particular, if the master equation corresponding to $(H', G')$ also has a (unique) classical solution $V'\in \cC^{1,4}_b(\Th)$, then we have, denoting $\D V := V-V'$,
 \bea
 \label{stability2}
\sum_{k=0}^2 \|\pa^{(k)}_x\D V\| \le C\sum_{k=0}^2\big[ \|\pa^{(k)}_{(x,z,\g)}\D H\|_{R} +\|\pa^{(k)}_x\D G\| \big]. %,\q\mbox{for all}~ (t,x,\mu)\in\Th.
 \eea
\end{thm}
We note that the uniqueness of $V$ and its decoupling property in (i) is standard in the drift control case. The stability result in (ii) is not surprising, but to the best of our knowledge, it is new even in the drift control case. This theorem will be proved in Section \ref{sect-stability} below.

\begin{rem}
\label{rem-stability}
In the above theorem we actually require only that $G\in \mathcal{C}^{1,4}_b(\dbR\times\mathcal{P}_2(\dbR))$, $H\in \mathcal{C}^{1,4}(\dbR^3\times\mathcal{P}_2(\dbR))$ with bounded derivatives when restricting $(z,\g)$ to $D_R$, and $u\in \cC^{1,4}_b([0, T]\times \dbR)$, and thus correspondingly $V\in C_b^{1,4}(\Th)$. For the existence result below, we raise all the regularity requirements to $\cC^{2,6}_b$, mainly because we want to apply  It\^{o}'s formula \reff{Ito} on %$\mathfrak{D}_x V(t, X_t, \cL_{X_t})$ and
 $\mathfrak{D}_\mu V(t, X_t, \cL_{X_t}, \tilde X_t)$ for some appropriate process $X$. This regularity requirement actually can be weakened by approximation arguments. We do not pursue this generality  in the paper for the ease of presentation.
\end{rem}

We emphasize that the uniqueness of classical solutions of the master equation \reff{master} is guaranteed by the above theorem. The uniqueness of the equilibrium measure $\bm$ is essentially equivalent to the existence of classical solution $V$ to the master equation. Our main goal of this paper is thus to study the existence of classical solution $V$.  We first establish the local well-posedness, whose proof will be reported in Section \ref{sect-local}.

\begin{thm}
\label{thm-local}
Let Assumptions  \ref{assum-regHG} and \ref{assum-regu} hold true. Then there exists a constant $\delta_0>0$, depending only  on $c_0$, $L_G$, $L^*_1$ (with $T=1$ in Assumption \ref{assum-regu}), and $L_H(R)$ with $R=L^*_1$, such that whenever $T\leq \delta_0$, the master equation \reff{master} has a unique classical solution $V\in \cC_b^{2,6}(\Theta)$. 
\end{thm}

\subsection{The global well-posedness} %Assumptions}
%\label{subsect-assum}
%\set counter{equation}{0}

We now turn to the main goal of the paper: the global well-posedness. As standard in the literature, this typically requires certain monotonicity condition. Since in this paper we focus only on the Lasry-Lions monotonicity, as in the drift control case we shall require  the Hamiltonian $H$ to be separable in the following sense. 
%We  will consider non-separable $H$ under more general second-order monotonicity conditions in a subsequent work. 

\begin{assum}
\label{assum-separable}
The Hamiltonian $H$ is separable, that is, there exist functions $H_0:\mathbb{R}^3\rightarrow \mathbb{R} $ and $F_0:\mathbb{R}\times \mathcal{P}_2(\dbR)\rightarrow \mathbb{R} $ such that
\bea
\label{separable}
H(x,z,\gamma,\mu)= H_0(x,z,\gamma)+ F_0(x,\mu),  \q  (x, z,\gamma, \mu)\in \dbR^3 \times \mathcal{P}_2(\dbR).
\eea
\end{assum}
In the case, the master equation \reff{master} becomes:
\bea
\label{master0}
\left.\ba{lll}
 \dis \partial_tV+H_0(x,\partial_{x}V,\partial_{xx}V)+F_0(x,\mu) +\mathcal{N}_0V=0,\q  V(T,x,\mu)=G(x,\mu),\q\mbox{where}\\ 
\dis \mathcal{N}_0V(t,x,\mu):=\int_{\mathbb R}\Big[\partial_{\mu}V(t,x,\mu,\tilde x)\partial_{z}
                          H_0(\tilde x,\partial_{x}V(t,\tilde x,\mu),\partial_{xx}V(t,\tilde x,\mu))\\
                          ~~~~~~~~~~~~~~~~~~+\partial_{\tilde{x}\mu}V(t,x,\mu,\tilde x)
                          \partial_{\gamma}H_0(\tilde x,\partial_{x}V(t,\tilde x,\mu),\partial_{xx}V(t,\tilde x,\mu))\Big]\mu(d\tilde x),
\ea\right.
\eea
and the HJB  equation  \eqref{HJB} becomes
\bea
\label{HJB0}
\pa_t u+H_0(x,\pa_xu,\pa_{xx}u)+F_0(x,m_t)=0,\quad u(T,x) = G(x, m_T).
\eea
Moreover, the Fokker-Planck equation \reff{Fokker-Planck} and the MFG system \reff{FBSDE} can also be simplified by the following facts: for any $k, l\ge 1$,
\bea
\label{HH0partial}
\left.\ba{c}
\dis \pa^{(k)}_{(z,\g)} H(x,z,\g, \mu) =  \pa^{(k)}_{(z,\g)} H_0(x,z,\g),\q  \pa^{(l)}_x\pa^{(k)}_{(z,\g)} H(x,z,\g, \mu) =  \pa^{(l)}_x\pa^{(k)}_{(z,\g)} H_0(x,z,\g),\\
\dis \pa^{(l)}_xH(x,z,\g, \mu) =  \pa^{(l)}_x H_0(x,z,\g) + \pa^{(l)}_xF_0(x, \mu).
\ea\right.
 \eea
In particular, it becomes slightly easier to verify Assumption \ref{assum-regu} for \reff{HJB0}.

We now specify the technical conditions on $H_0$ and $F_0$. For this purpose, we denote: 
\bea
\label{pa2H0}
\pa_{(z,\gamma)}^2H_0:=\left(
  \begin{array}{cc}
    \partial_{zz}H_0 & \partial_{z\gamma}H_0 \\
    \partial_{z\gamma}H_0 & \partial_{\gamma\gamma}H_0 \\
  \end{array}
\right).
\eea
We shall not confuse this notation with $\pa^{(2)}_{(z,\g)} H_0$ though, which refers generally  to the second-order derivatives $\pa_{zz}H_0$, $\pa_{z\g}H_0$, $\pa_{\g\g}H_0$, but not the specific matrix.
%Recall \reff{pa2H0}.

\begin{assum}\label{assum-regH0}
(i) $F_0\in \mathcal{C}_b^{2,6}(\dbR\times\mathcal{P}_2(\dbR))$ and there exists a constant $L_{F_0}>0$ such that
\begin{eqnarray*}
%\|\mathfrak{D}_x F_0\|\leq L^{F_0}_0, \ \
%  \|\mathfrak{D}_\mu F_0\|\leq L^{F_0}_1,\q \forall (x,\mu,\tilde x)\in \dbR\times\cP_2(\dbR)\times\dbR.
  \|\pa^{(k)}_x F_0\|, ~ \|\pa^{(l)}_x\pa_\mu F_0\|,   ~ \|\pa^{(l)}_x\pa_{\tilde \mu} F_0\|\leq L_{F_0},\q k=1,2,3,4;~ l=1,2.
\end{eqnarray*}

\no(ii)  $H_0\in C^6(\dbR^3)$  and, for any $R>0$, all the involved derivatives are bounded when restricting $(z, \g)$ to $D_R$, and  there exists a constant $L_{H_0}(R)>0$ depending on $R$ such that
\beaa
\|\pa^{(k)}_{(x,z,\g)} H_0\|_R \le L_{H_0}(R), \q k=1,2,3,4.
\eeaa
 
 \ms
\no (iii) $\pa_\g H_0\ge c_0$ for some $c_0>0$.

\ms
\no(iv) For any $R>0$, there exists $c_R>0$ such that, denoting by $I_{2\times 2}$ the $2\times 2$ identity matrix,
\begin{eqnarray}\label{pa2H0-negative}
\pa_{(z,\gamma)}^2H_0(x,z,\gamma)\leq-c_R I_{2\times 2} ,\,\, \forall(x,z,\gamma)\in \dbR\times D_R.
\end{eqnarray}
\end{assum}

Finally, we require the crucial monotonicity condition. 
%We now assume our data $G,f$ satisfy the Lasry-Lions monotonicity condition.
\begin{assum}
\label{assum-LLmon}
 $G$ and $F_0$ satisfy the Lasry-Lions monotonicity condition \eqref{LLmon1}.
\end{assum}

Our main result of the paper is  the following global well-posedness of master equation \reff{master0}.
% for  MFGs with volatility control.

\begin{thm}\label{thm-global}
Let Assumptions \ref{assum-regHG}-(i), \ref{assum-regu}, \ref{assum-separable}, \ref{assum-regH0},  and \ref{assum-LLmon} hold. Then for any $T$, the master equation \reff{master0}  admits a unique classical solution $V\in \cC_b^{2,6}(\Theta)$. 
\end{thm}

We shall prove this theorem in the next section.

\section{The global well-posedness: Proof of Theorem \ref{thm-global}}
\label{sect-global}
\setcounter{equation}{0}

%The main ideas to prove the theorem are as follows. 
%We conclude this section with the road map towards the global well-posedness of master equation \reff{master}.

We follow the approach in \cite{GMMZ} and proceed in three steps, assuming the local well-posdness. 

\emph{Step 1.} Prove that the  Lasry-Lions monotonicity condition \reff{LLmon1} 
%(which is equivalent to the second-order Lasry-Lions monotonicity condition) 
propagates along the classical solution $V$ to the master equation \reff{master0}. Consequently, by the equivalence in Proposition \ref{prop:2nd}, the second-order  Lasry-Lions monotonicity condition \reff{LLmon2} also propagates.

\emph{Step 2.} Show that the second-order  Lasry-Lions monotonicity condition \reff{LLmon2} implies the (a priori) uniform Lipschitz continuity of $\partial_xV$ and $\partial_{xx}V$ in the measure variable $\mu$;

\emph{Step 3.} Combine the local well-posedness result in Theorem \ref{thm-local} and the above uniform Lipschitz estimate  to obtain the global well-posedness of the master equation \reff{master0}.

The arguments for Step 3 are more or less standard. We remark that Step 2 above remains valid under the second-order displacement monotonicity condition \reff{DPmon2}.  However, since \reff{DPmon2} is not equivalent to the (first-order) displacement monotonicity condition \reff{DPmon1}, then in Step 1 we have to propagate \reff{DPmon2} directly, which is much more involved. We thus focus on the main ideas in this paper, by restricting to separable $H$ under Lasry-Lions monotonicity. We will study non-separable $H$ under more general second-order monotonicity conditions in a subsequent work.

\subsection{Propagation of monotonicity}
\label{sect-propagation}

In this subsection we show that, under our standing assumptions, the second-order Lasry-Lions monotonicity
condition propagates along the classical solution of master equation \reff{master0}. 

\begin{thm}\label{thm-propagation}
Let all the assumptions in Theorem \ref{thm-global} hold true, and $V\in \cC_b^{1,4}(\Theta)$ be a classical solution to the  master equation \reff{master0}. Then for any $t$,  $V(t,\cdot)$ satisfies \reff{LLmon1} for any $\xi, \eta\in\mathbb{L}^2(\mathcal{F}_t)$, and hence the second-order Lasry-Lions monotonicity
condition \reff{LLmon2} for any $\xi, \eta_1,\eta_2\in\mathbb{L}^2(\mathcal{F}_t)$.
\end{thm}
\proof Since $V(t,\cdot)\in\cC_b^{1,4}(\dbR\times\cP_2(\dbR))$, by  Proposition \ref{prop:2nd} we see that \reff{LLmon2} follows from \reff{LLmon1}. Thus we shall only verify \reff{LLmon1} for $V(t,\cd,\cd)$, and without loss of generality considering only $t=0$:
\begin{equation}\label{eq:V:1stLL0}
\cM^1_{LL} V(0, \xi, \eta) := \dbE\big[\pa_{x\mu}V(0,\xi,\cL_{\xi},\tilde\xi)\eta \tilde\eta\big]\geq 0,\,\,\forall \xi,\eta\in\dbL^2(\cF_0).
\end{equation}

For this purpose, we fix $\xi,\eta\in\dbL^2(\cF_0)$, and by Theorem \ref{thm-stability}, let $(\bm, u)$ be the unique solution to the MFG system \reff{Fokker-Planck}-\reff{HJB0}. Let $X:= X^\xi$ be the solution to SDE \reff{XxiwhH}, and by formally differentiating \reff{XxiwhH} with respect to $\xi$, we introduce:
%us consider the following McKean-Vlasov SDEs:
\bea\label{propagation-dX}
\left.\ba{lll}
\dis \delta X_t=\eta+\int^t_0 \Big[\th_1(X_s) \d X_s + \dbE_{\cF_s}[\th_2(X_s, \tilde X_s) \widetilde {\d X_s}]\Big]ds \\
\dis \qq\q +\int^{t}_0\Big[\b_1(X_s) \d X_s + \dbE_{\cF_s}[\b_2(X_s, \tilde X_s) \widetilde {\d X_s}]\Big]\slash \sqrt{2\partial_{\gamma}\cH_V (X_s)}dB_s,
\ea\right.
\eea
where,  in light of \reff{Hoptimal} and \reff{HH0partial}, 
%by denoting $\pa_{xz} \cH_V := \wh{\pa_{xz} H}_u$ and similarly for the other derivatives, 
\bea
\label{propagation-dX2}
\left.\ba{lll}
\dis \pa_{z} \cH_V (x) := \pa_{z} H_0(x, \pa_x V, \pa_{xx} V),~\mbox{and similarly for the other derivatives};\\
\dis \th_1(x) :=\pa_x[\pa_z \cH_V ](x) = \partial_{xz}\cH_V (x)+\partial_{zz}\cH_V (x)\partial_{xx}V(x)+\partial_{ z\gamma}\cH_V (x)\partial_{xxx}V(x);\\
\dis \th_2(x, \tilde x) := \pa_\mu[\pa_z \cH_V ](x, \tilde x)= \pa_{zz}\cH_V (x)\pa_{x\mu}V(x,\tilde x)+\pa_{z\gamma}\cH_V (x)\pa_{xx\mu}V(x,\tilde x);\\
\dis \b_1(x) := \pa_x[\pa_\g \cH_V ](x)= \pa_{x\gamma}\cH_V (x)+\pa_{z\gamma}\cH_V (x)\pa_{xx}V(x)+\pa_{\gamma\gamma}\cH_V (x)\pa_{xxx}V(x);\\
\dis \b_2(x,\tilde x) :=  \pa_\mu[\pa_\g \cH_V ](x, \tilde x)=\pa_{z\g}\cH_V (x)\pa_{x\mu}V(x,\tilde x)+\pa_{\g\gamma}\cH_V (x)\pa_{xx\mu}V(x,\tilde x). 
\ea\right.
\eea
Here and in the sequel, when the contexts are clear, for notational simplicity we omit the variables $(t, m_t)$ inside the derivatives of $\cH_V $ and $V$. 
Under our regularity conditions, clearly the linear McKean-Vlasov SDE \reff{propagation-dX} has a unique solution $\d X$. 

By straightforward calculation, with the details postponed to Appendix, we can show that
\bea
\label{propagation-claim}
&\dis {d\over dt}\cM^1_{LL} V(t, X_t, \d X_t) =  -\cM^1_{LL} F_0(X_t, \d X_t) +\dbE\Big[\G_t ~ \pa_{(z,\gamma)}^2\cH_V (X_t) ~\G_t^\top\Big],\\
&\dis \mbox{where}\q \G_t := \left( \dbE_{\cF_t}\big[\pa_{x\mu}V(X_t,\tilde X_t)\widetilde{\delta X_t}\big],~    \dbE_{\cF_t}\big[\pa_{xx\mu}V(X_t,\tilde X_t)\widetilde{\delta X_t}\big]\right).\nonumber
\eea
Then, by the monotonicity of $F_0$ in Assumption \ref{assum-LLmon}  and the concavity of $H_0$ in \reff{pa2H0-negative},  we have
\begin{equation*}
{d\over dt} \cM^1_{LL} V(t, X_t, \d X_t) \leq 0.
\end{equation*}
Moreover, by the monotonicity of $G$ in Assumption \ref{assum-LLmon}, we have 
$$\cM^1_{LL} V(T, X_T, \d X_T) = \cM^1_{LL} G(X_T, \d X_T) \ge 0.
$$ Then the above implies $\cM^1_{LL} V(0, \xi, \eta) = \cM^1_{LL} V(0, X_0, \d X_0)\ge 0$. 
\qed

\subsection{A priori $W_1$-Lipschitz continuity}
\label{sect-Lip}

The main result in this subsection is an a priori $W_1$-Lipschitz continuity of $\pa_{x}V,\pa_{xx}V$ with respect to $\mu$ under the second-order Lasry-Lions monotonicity condition \eqref{LLmon2}.
\begin{thm}\label{thm-Lip}
Let all the assumptions of Theorem \ref{thm-propagation} hold, with a possible exception of Assumption \ref{assum-LLmon} for $F_0$, and $V\in \cC_b^{2,6}(\Th)$ be a classical solution of master equation \reff{master0}. Assume further that $V(t,\cdot)$ satisfies \eqref{LLmon2} for any $t\in [0,T]$. Then $ \mathfrak{D}_\mu V(t,\cd)$ is uniformly bounded by some constant $L^V_T\ge L^*_T$,
 which depends only on $T$, $c_0$, $L_G$, $L_{F_0}$,  $L^*_T$, $L_{H_0}(R)$, and  $c_R$ in \reff{pa2H0-negative}, with $R:=L^*_T$.
\end{thm}
To prove the theorem, we first need a technical lemma. Fix $\xi\in \dbL^2(\cF_0)$ and $x\in \dbR$. As in the previous section,  let $(\bm, u)$ be the solution to the MFG system \reff{Fokker-Planck}-\reff{HJB0}, and $X:= X^\xi$ be the solution to SDE \reff{XxiwhH}. %Let $X^x$ denote the solution to SDE \reff{XxiwhH} with $\xi\equiv x$. 
Again we drop the variables $(t, m_t)$ when the contexts are clear. 
Moreover, to simplify the notations we introduce
\begin{equation*}
 \mathfrak{D}_\mu V(t,x,\mu,\tilde x)=: \left(
  \begin{array}{cc}
   U_1 & U_3 \\
 U_2 & U_4\\
  \end{array}
\right)(t,x,\mu,\tilde x).
\end{equation*}

\begin{lem}\label{lem-Lip}
Let all the assumptions of Theorem \ref{thm-Lip} hold true. Then 
\bea
\label{Lip-claim1}
&&\dis\!\!\!\!\!\!\!\!\!\!\!\!\!\!\!\!\!\! d\pa_{(z,\gamma)}^2\cH_V (X_t)=\Gamma_0^V(X_t)dt+\Gamma_1^V(X_t)dB_t;\\
\label{Lip-claim2}
&&\dis \!\!\!\!\!\!\!\!\!\!\!\!\!\!\!\!\!\! d \mathfrak{D}_\mu V(X_t,\tilde X_t) = \big[\Gamma_2^V+\Gamma_3^V + \G^V_4\big](X_t,\tilde X_t)dt +\Gamma_5^V(X_t,\tilde X_t)dB_t+\Gamma_6^V(X_t,\tilde X_t)d\tilde B_t,
\eea
where $\G^V_0, \G^V_1: \Th\to\dbR^{2\times 2}$ and $\Gamma_k^V:\Th\times \dbR\to\dbR^{2\times 2}$, $k=2,\cds,6$, satisfy, recalling \reff{propagation-dX2},
\bea
\label{Lip-Gamma}
\left.\ba{lll}
\dis |\Gamma_0^V(x)|, |\Gamma_1^V(x)|\leq C\q\text{and}\q |\Gamma_2^V(x,\tilde x)|\leq C\big[| \mathfrak{D}_\mu V(t,x,\tilde x)|+1\big],\ms\\
\dis \Gamma_3^V(x,\tilde x):= 2\left(
  \begin{array}{cc}
   0 &\b_1(\tilde x) \pa_{\tilde x}U_3(x,\tilde x) \\
  \b_1(x)\pa_{x}U_2(x,\tilde x) & \b_1(x)\pa_{x}U_4(x,\tilde x)+ \b_1(\tilde x)\pa_{\tilde x}U_4(x,\tilde x)
  \end{array}
\right),\ms\\
\dis \Gamma_4^V(x,\tilde x):= -\dbE\big[
  \mathfrak{D}_\mu V(x,\bar X_t)\pa_{(z,\gamma)}^2\cH_V (\bar X_t) \mathfrak{D}_\mu V(\bar X_t,\tilde x)\big];\ms\\
\dis \Gamma_5^V(x,\tilde x):=\sqrt{2\pa_{\gamma}\cH_V (x)}\left(
  \begin{array}{cc}
   U_2(x,\tilde x) & U_4(x,\tilde x) \\
  \pa_{x}U_2(x,\tilde x)& \pa_{x}U_4(x,\tilde x) \\
  \end{array}
\right),\ms\\
\dis\Gamma_6^V(x,\tilde x):= \sqrt{2\pa_{\gamma}\cH_V (\tilde x)} \left(
  \begin{array}{cc}
   U_3(x,\tilde x) & \pa_{\tilde x}U_3(x,\tilde x) \\
  U_4(x,\tilde x)& \pa_{\tilde x}U_4(x,\tilde x) \\
  \end{array}
\right).
\ea\right.
\eea
Here the constant $C$ depends on $c_0$, $L_{F_0}$, $L^*_T$, and $L_{H_0}(R)$ with $R=L^*_T$. 
%the parameters specified in Theorem \ref{thm-Lip} for $L^V_T$.
\end{lem}

The proof is again by straightforward calculation, and we postpone it to Appendix.

\bs
\no{\bf Proof of Theorem \ref{thm-Lip}.} Without loss of generality, it suffices to show that $\mathfrak{D}_\mu V(0,\cd)$ is uniformly bounded by some constant $L_T^V>0$. For this purpose, introduce
\bea
\label{Lip-Lambda}
\L_t &:=& \tr\dbE_{\cF_0\vee \tilde\cF_0}\big[  \mathfrak{D}_\mu V(X_t,\tilde X_t)\pa_{(z,\gamma)}^2\cH_V (\tilde X_t)  \mathfrak{D}_\mu V(X_t,\tilde X_t)^\top\big] \nonumber\\
&\le& -c_R \dbE_{\cF_0\vee \tilde\cF_0}\big[ |\mathfrak{D}_\mu V(X_t,\tilde X_t)|^2\big] \le 0,
\eea
where the first inequality is due to \reff{pa2H0-negative}. We claim that
\bea
\label{Lip-dLambda}
{d\over dt} \L_t \le C  \dbE_{\cF_0\vee \tilde\cF_0}\big[ |\mathfrak{D}_\mu V(X_t,\tilde X_t)|^2+1\big],
\eea
where $C$ depends on the parameters specified for $L^V_T$. Then by \reff{Lip-Lambda} we have
\beaa
{d\over dt} \L_t \le - C \L_t+C.
\eeaa
Since $\L_T \ge -C$, the above implies that $\L_0 \ge -C$. Then, by \reff{Lip-Lambda} again,
\beaa
|\mathfrak{D}_\mu V(0, \xi, \mu, \tilde \xi)|^2 = \dbE_{\cF_0\vee\tilde\cF_0}\big[ |\mathfrak{D}_\mu V(X_0,\tilde X_0)|^2\big]  \le -C\L_0 \le C, \q\mbox{a.s.}
\eeaa
That is, setting $L^V_T := \sqrt{C}$ for the above $C$,
\beaa
|\mathfrak{D}_\mu V(0, x, \mu, \tilde x)| \le L^V_T, \q\mbox{$\mu$-a.e. $x$ and $\mu$-a.e. $\tilde x$}.
\eeaa
Now by the continuity of $\mathfrak{D}_\mu V$, we see that the above holds true for all $x, \mu, \tilde x$.

It remains to prove \reff{Lip-dLambda}. 
Combining \eqref{Lip-claim1} and \eqref{Lip-claim2}, and further taking conditional expectation $\dbE_{\cF_0\vee\tilde\cF_0}[\cd]$, we can derive from the standard It\^o's formula that
\begin{eqnarray}\label{eq:d3U}
&&{d\L_t \over dt} = I_1(t) + I_2(t) + I_3(t) + I_4(t),\q\mbox{where}\\
&&I_1(t) := \tr\dbE_{\cF_0\vee\tilde \cF_0}\Big[2\G^V_2 \pa_{(z,\gamma)}^2\cH_V (\tilde X_t)  \mathfrak{D}_\mu V^\top + \mathfrak{D}_\mu V\G^V_0(\tilde X_t)  \mathfrak{D}_\mu V^\top\Big]; \nonumber\\
&&I_2(t) := 2\tr\dbE_{\cF_0\vee\tilde \cF_0}\Big[\big[\G^V_3 \pa_{(z,\gamma)}^2\cH_V (\tilde X_t) +\Gamma_6^V\Gamma_1^V(\tilde X_t)\big]  \mathfrak{D}_\mu V^\top\Big];\nonumber\\
&&I_3(t) := 2\tr\dbE_{\cF_0\vee\tilde \cF_0}\Big[\G^V_4 \pa_{(z,\gamma)}^2\cH_V (\tilde X_t)  \mathfrak{D}_\mu V^\top\Big];\nonumber\\
&&I_4(t) := \tr\dbE_{\cF_0\vee\tilde \cF_0}\Big[\Gamma_5^V\pa_{(z,\gamma)}^2\cH_V ( \tilde X_t)(\Gamma_5^V)^\top+\Gamma_6^V\pa_{(z,\gamma)}^2\cH_V (\tilde X_t)(\Gamma_6^V)^\top\Big].\nonumber
\end{eqnarray}
Here again we drop the arguments $(X_t, \tilde X_t)$ insider $\mathfrak{D}_\mu V$ and $\G^V_i$. 
We now estimate $I_i(t)$ separately. First, since $V(t,\cd)$ satisfies \eqref{LLmon2} and $\pa_{(z,\gamma)}^2\cH_V $ is symmetric, then
\bea
\label{eq:d3Ugeq1}
&&\dis\!\!\!\!\!\!\! I_3(t) =-2\tr\dbE_{\cF_0\vee\tilde\cF_0}\Big[\mathfrak{D}_\mu V(X_t,\bar X_t)\pa_{(z,\gamma)}^2\cH_V (\bar X_t)\mathfrak{D}_\mu V(\bar X_t,\tilde X_t)\pa_{(z,\gamma)}^2\cH_V (\tilde X_t)\mathfrak{D}_\mu V( X_t,\tilde X_t)^\top\big]\Big]\nonumber\\
&&\!\!\!\!\!\!\!=\! -2\dbE_{\cF_0\vee\tilde \cF_0}\Big[\tr \dbE_{\cF_t \vee \tilde \cF_0 \vee \bar\cF_0}\big[\mathfrak{D}_\mu V(X_t,\bar X_t)\pa_{(z,\gamma)}^2\cH_V (\bar X_t)\mathfrak{D}_\mu V(\bar X_t,\tilde X_t)\pa_{(z,\gamma)}^2\cH_V (\tilde X_t)\mathfrak{D}_\mu V( X_t,\tilde X_t)^\top\big]\Big]\nonumber\\
&&\!\!\!\!\!\!\!\leq 0.
\eea
Next,  by \reff{pa2H0-negative}  we have
\bea
\label{eq:d3Ugeq2}
I_4(t) \le - c_R \dbE_{\cF_0\vee\tilde \cF_0}\Big[|\Gamma_5^V|^2 + |\Gamma_6^V|^2\Big].
%&\le&-\e_0  \dbE_{\tilde \cF_0}\Big[|\mathfrak{D}_\mu V(X_t^x,\tilde X_t)|^2 + \big[|\pa_x U_3|^2 + |\pa_x U_4|^2 + |\pa_{\tilde x} U_2|^2 + |\pa_{\tilde x}U_4|^2\big](X_t^x,\tilde X_t)\Big].
\eea
Moreover, by  \reff{Lip-Gamma} we can easily see that, for any $\e>0$,
\bea
\label{eq:d3Ugeq3}
\left.\ba{lll}
\dis I_1(t) \le  C \dbE_{\cF_0\vee\tilde\cF_0}\big[ |\mathfrak{D}_\mu V|^2+1\big];\ms\\
\dis I_2(t) \le  \dbE_{\cF_0\vee \tilde\cF_0}\big[ C\e^{-1} |\mathfrak{D}_\mu V|^2\big] + \e \big[| \G^V_3|^2 + | \G^V_6|^2\big]\Big] \ms\\
\dis\qq \le  \dbE_{\cF_0\vee \tilde\cF_0}\Big[C\e^{-1} |\mathfrak{D}_\mu V|^2 + C_0\e  \big[| \G^V_5|^2 + | \G^V_6|^2\big]\Big],
\ea\right.
\eea
where we used Assumption \ref{assum-regH0} (iii) for the last inequality. Set $\e := {c_R\over C_0}$ for the above $C_0$, and plug 
 \reff{eq:d3Ugeq1},  \reff{eq:d3Ugeq2},  \reff{eq:d3Ugeq3} into \reff{eq:d3U}, we obtain \reff{Lip-dLambda} immediately.
 \qed

\subsection{Proof of Theorem  \ref{thm-global}.}
\label{sect-proof}
We now establish the global well-posedness of master equation \reff{master0}, provided the local existence in Theorem \ref{thm-local}. The uniqueness is already given by Theorem \ref{thm-stability}, so it suffices to prove the existence. As illustrated in \cite{CD2,CCD,GMMZ, MZ2}, the key to extend a local classical solution to a global one is the a priori $W_1$-Lipschitz continuity estimate obtained in the previous subsection.

%Set $L_0^V:= L^*_T$ for the 
Let $\delta_0>0$ be the constant in Theorem \ref{thm-local}, but with $L_G$
replaced with the $L^V_T$ in Theorem \ref{thm-Lip}. We note that $L_G$ is the bound of $\pa^{(k)}_xG$, $k=1,2,3,4$ and $\mathfrak{D}_\mu G$.
Let $0=T_0<\cdots<T_n=T$ be a partition such that $T_{i+1}-T_i\leq \frac{\delta_0}{2}$, $i=0,\cdots, n-1$. We now construct a solution backwardly. 

 First, since $T_n-T_{n-2}\leq \delta_0$, by Theorem \ref{thm-local}  master equation \reff{master0} on $[T_{n-2}, T_n]$ with terminal condition $G$
 has a unique classical solution $V\in \cC_b^{2,6}([T_{n-2},T_n]\times\dbR\times\cP_2(\dbR))$. For any $(t_0,\mu)\in [T_{n-2},T_n]\times\cP_2(\dbR)$ and $\xi\in\dbL^2(\cF_{t_0};\mu)$, let $X^{t_0,\xi}$ be the  solution to the McKean-Vlasov SDE: 
 \bea
\label{XV0}
\left.\ba{c}
\dis X_t^{t_0,\xi} = \xi + \int_0^t \pa_z\cH_V( s, X^{t_0,\xi}_s, \cL_{X^{t_0,\xi}_s})ds +  \int_0^t \sqrt{2\pa_\g \cH_V( s, X^{t_0,\xi}_s, \cL_{X^{t_0,\xi}_s})}dB_s.
%\dis \mbox{where}\q \wh \f_V(t,x,\mu) := \f(x, \pa_x V(t, x, \mu),  \pa_{xx} V(t,x, \mu), \mu),~ \f = \pa_z H, \pa_\g H. 
\ea\right.
\eea
Then we can verify that $u(t,x):=V(t,x,\cL_{X_{t}^{t_0,\xi}})$ solves HJB equation \eqref{HJB0}, and thus, by Assumption \ref{assum-regu},  $\sup_{1\le k\le 4}\|\pa^{(k)}_x V(t,\cd)\|\leq L^*_T\le L^V_T$, $t\in [T_{n-2},T_n]$. Moreover, it follows from Theorem \ref{thm-Lip} that $\|\mathfrak{D}_\mu V(t,\cd)\|\leq L^V_T$,  for any $t\in [T_{n-2},T_n]$. 

 We next consider master equation \reff{master0} on $[T_{n-3}, T_{n-1}]$ with terminal condition $V(T_{n-1}, \cdot)$.
We emphasize that $V(T_{n-1}, \cdot)$ satisfies all the requirements on $G$ in Assumption \ref{assum-regHG}-(i) with the constant $L^V_T$ in \reff{LG}; then we may apply Theorem \ref{thm-local} with the same $\delta_0$ and obtain a classical solution
$V\in C_b^{2,6}([T_{n-3}, T_{n-1}]\times\dbR\times\cP_2(\dbR))$. Clearly
this extends the classical solution of the master equation to $[T_{n-3}, T_{n}]$. We emphasize again
that, by Assumption \ref{assum-regu} and Theorem \ref{thm-Lip}  we see that $\sup_{1\le k\le 4}\|\pa^{(k)}_x V(t,\cd)\|\leq L^*_T\le L^V_T$ and $\|\mathfrak{D}_\mu V(t,\cd)\|\leq L_T^V$ for all $t\in [T_{n-3}, T_{n}]$. This enables us to consider master equation \reff{master0} on
$[T_{n-4}, T_{n-2}]$ with terminal condition $V(T_{n-2}, \cdot,\cdot)$, and then we obtain a classical solution on
$[T_{n-4}, T_{n}]$ with the desired uniform estimates and additional regularities.\par
Now repeat the arguments backwardly in time, we may construct a classical solution $V$
for master equation \reff{master0} on $[0, T]$ with terminal condition $G$.
\qed

\section{The stability: Proof of Theorem \ref{thm-stability}}
\label{sect-stability}
\setcounter{equation}{0}

First, the uniqueness of classical solution $V$ follows directly from \reff{stability2}. 

\no(i) By the equivalence of the mean field equilibrium and the two MFG systems, it suffices to prove the results concerning \reff{HJB}-\reff{Fokker-Planck}. The uniqueness follows from \reff{stability1}. We thus prove the existence for \reff{HJB}-\reff{Fokker-Planck} as well as the first equality in \reff{decoupling}.

Set $t_0=0$ and fix $\xi \in \dbL^2(\cF_0,\mu)$ in \reff{XV0}.  Under our conditions, clearly \reff{XV0}  has a unique solution. Denote $X:= X^{0,\xi}$, and set
\beaa
m_t := \cL_{X_t},\q u(t, x) := V(t, x, m_t).
\eeaa 
Applying It\^{o}'s formula \reff{Ito} on $V(t, x, \cL_{X_t})$ and  by using the master equation \reff{master}, we obtain
\beaa
{d\over dt} V(t,x, m_t) \!\!\!&=&\!\!\! \pa_t V + \dbE\Big[ \pa_\mu V(t,x, m_t, \tilde X_t) \pa_z\cH_V( t, \tilde X_t, m_t) +\pa_{\tilde x \mu} V(t, x, m_t, \tilde X_t)  \pa_\gamma \cH_V( t, \tilde X_t, m_t)\Big]\\
\!\!\!&=&\!\!\! - H(x, \pa_x V, \pa_{xx} V, m_t).
\eeaa
Note that $\pa_x V = \pa_x u$, $\pa_{xx} V = \pa_{xx} u$. Then 
\beaa
\pa_t u( t,x) = {d\over dt} V(t,x, m_t) =  - H(x, \pa_x u, \pa_{xx} u, m_t).
\eeaa
This, together with the obvious terminal condition, verifies \reff{HJB}.  Moreover, since $u = V$, one can easily see that the $ \cH_V$ in \reff{XV0} identifies with the $\cH_u$ in \reff{Hoptimal}. Then \reff{XV0} implies that $\bm := \cL_X$ is a weak solution of  \reff{Fokker-Planck}. We thus verify \reff{HJB}-\reff{Fokker-Planck} and the first equality in \reff{decoupling}.

\ms
\no (ii) We first prove \reff{stability1}. At below the generic constant $C$ depends only on the parameters specified in the theorem. Recall \reff{Fokker-Planck} and fix $\xi$. Clearly $(\bm', u')$ leads to a solution $\Xi'=(X', Y^{1'},  Y^{2'}, \cL_{X'})$ to the corresponding MFG system \reff{FBSDE}, with
\beaa
\bm' = \cL_{X'}, \q Y^{k'}_t = \pa^{(k)}_xu'( t, X'_t),~k=0,1,2.
\eeaa
In particular,  
\beaa
d X'_t = \pa_z \cH'_{u'}(t, X'_t, m'_t) dt +  \sqrt{2 \pa_\g \cH'_{u'}(t, X'_t, m'_t)}dB_t.
\eeaa
Note that $ \pa_z \cH'_{u'}(t, X'_t, m'_t) = \pa_z H'(\Xi'_t)$ and $ \pa_\g \cH'_{u'}(t, X'_t, m'_t) = \pa_\g H'(\Xi'_t)$. Denote, for $k=0,1,2$, 
\beaa
&\dis Y^{k"}_t := \pa^{(k)}_x V(t, X'_t, m'_t),\q Z^{k"}_t := \pa^{(k+1)}_{x} V(t, X'_t, m'_t)\sqrt{2\pa_\g H'(\Xi'_t)},\\
%&&\dis  Y^{2"}_t := \pa_{xx} V(t, X'_t, m'_t),\q Z^{2"}_t := \pa_{xxx} V(t, X'_t, m'_t)\sqrt{2\pa_\g H'(\Xi'_t)};\\
&\dis \D Y^k := Y^{k"}-Y^{k'},\q \D Z^k := Z^{k"}-Z^{k'}.
\eeaa

Applying It\^{o}'s formula \reff{Ito} on $V(t, X'_t, m'_t)$ we have
\bea
\label{dY0-1}
&&\dis d  Y^{0"}_t = \Big[\pa_{t} V + \pa_{x} V \pa_z \cH'_{u'} + \pa_{xx} V  \pa_\g \cH'_{u'}\Big](t, X'_t, m'_t) dt + Z^{0"}_t dB_t\\
&&\q + \dbE_{\cF_t}\Big[\pa_{\mu} V (t, X'_t, m'_t, \tilde X'_t) \pa_z \cH'_{u'}(t, \tilde X'_t, m'_t) + \pa_{\tilde x \mu} V (t, X'_t, m'_t, \tilde X'_t) \pa_\g \cH'_{u'}(t, \tilde X'_t, m'_t) \Big]dt.\nonumber
\eea
Plug \reff{master} into \reff{dY0-1}, and denote $\D\cH := \cH_V - \cH'_{u'}$. Then, by \reff{FBSDE} we have
\bea
\label{dY0-3}
&&\dis d  \D Y^{0}_t = - \G^0_t  dt + \D Z^{0}_t dB_t,\q\mbox{where}\\
&&\dis \G^0_t := \Big[\D\cH- \pa_z \cH'_{u'} \D Y^1_t-\pa_\g \cH'_{u'}(t, \tilde X'_t, m'_t) \D Y^2_t \Big](t, X'_t, m'_t) \nonumber\\
&&\dis \q + \dbE_{\cF_t}\Big[\pa_{\mu} V (t, X'_t, m'_t, \tilde X'_t) \D\pa_z\cH( t, \tilde X'_t, m'_t) + \pa_{\tilde x \mu} V (t, X'_t, m'_t, \tilde X'_t) \D\pa_\g\cH( t, \tilde X'_t, m'_t) \Big].\nonumber
\eea
Note that, by Assumption \ref{assum-regu}, $\|\pa^{(k)}_x V\|\le L^*_T$, $\|\pa^{(k)}_x u'\|\le L^*_T$, $k=1,2,3,4$. Then, for a constant $C$ depending on $L_H(R)$ with  $R=L^*_T$, we have
\bea
\label{dY0-4}
|\pa_z \cH'_{u'}|, |\pa_\g \cH'_{u'}|\le C,\q \big|\D\cH( t, X'_t, m'_t)\big| \le C\Big[\|\D H\|_R + |\D Y^1_t| + |\D Y^2_t|\Big].
\eea
Similar estimates hold for $\D\pa_z \cH$ and $\D\pa_\g \cH$. Note  that $\tilde X'$ is an independent copy of $X'$, then 
%have the same distribution as $X'$. Then, since $V\in \cC^4_b(\Th)$  and $\D Y^1_T=0$, it follows from standard BSDE arguments that
\bea
\label{dY0-5}
| \G^{0}_t | \le C \Big[\|\D H\|_R  + |\D Y^1_t| + |\D Y^2_t| + \dbE_{\cF_t}\big[|\D \tilde Y^1_t| + |\D \tilde Y^2_t|\big]\Big].
\eea

Next, applying It\^{o}'s formula \reff{Ito} on $\pa_x V(t, X'_t, m'_t)$ we have
\bea
\label{dY1-1}
&&\dis d  Y^{1"}_t = \Big[\pa_{xt} V + \pa_{xx} V \pa_z \cH'_{u'} + \pa_{xxx} V  \pa_\g \cH'_{u'}\Big](t, X'_t, m'_t) dt + Z^{1"}_t dB_t\\
&&\q + \dbE_{\cF_t}\Big[\pa_{x\mu} V (t, X'_t, m'_t, \tilde X'_t) \pa_z \cH'_{u'}(t, \tilde X'_t, m'_t) + \pa_{\tilde x x\mu} V (t, X'_t, m'_t, \tilde X'_t) \pa_\g \cH'_{u'}(t, \tilde X'_t, m'_t) \Big]dt.\nonumber
\eea
On the other hand, differentiate the master equation \reff{master} with respect to $x$, we have
\bea
\label{dY1-2}
&&\dis \partial_{xt} V+ \Big[\pa_x\cH_V + \pa_z\cH_V \pa_{xx} V + \pa_\g\cH_V \pa_{xxx} V\Big](t, x,\mu)\\
&&\dis\q +\int_{\mathbb R^d}\Big[\partial_{x\mu}V(t,x,\mu,\tilde x) \pa_z\cH_V( t, \tilde x,\mu)+\partial_{\tilde xx\mu}V(t,x,\mu,\tilde x) \pa_\g\cH_V( t, \tilde x,\mu)\Big]\mu(d\tilde x)=0.\nonumber
\eea
Plug \reff{dY1-2} into \reff{dY1-1}, by \reff{FBSDE} we have
\bea
\label{dY1-3}
&&\dis d  \D Y^{1}_t = - \G^1_t  dt + \D Z^{1}_t dB_t,\q\mbox{where}\\
&&\dis \G^1_t := \Big[\D\pa_x\cH+ \pa_{xx} V \D\pa_z\cH + \pa_{xxx} V \D\pa_\g\cH \Big](t, X'_t, m'_t) \nonumber\\
&&\dis \q + \dbE_{\cF_t}\Big[\pa_{x\mu} V (t, X'_t, m'_t, \tilde X'_t) \D\pa_z\cH( t, \tilde X'_t, m'_t) + \pa_{\tilde x x\mu} V (t, X'_t, m'_t, \tilde X'_t) \D\pa_\g\cH( t, \tilde X'_t, m'_t) \Big].\nonumber
\eea
Similarly to \reff{dY0-5}, we have
\bea
\label{dY1-5}
| \G^{1}_t | \le C \Big[\|\pa^{(1)}_{(x,z,\g)}\D H\|_R  + |\D Y^1_t| + |\D Y^2_t| + \dbE_{\cF_t}\big[|\D \tilde Y^1_t| + |\D \tilde Y^2_t|\big]\Big].
\eea

Moreover, applying It\^{o}'s formula \reff{Ito} on $\pa_{xx} V(t, X'_t, m'_t)$ we have
\bea
\label{dY2-1}
&&\dis d  Y^{2"}_t = \Big[\pa_{xxt} V + \pa_{xxx} V  \pa_z \cH'_{u'} + \pa_{xxxx} V  \pa_\g \cH'_{u'} \Big](t, X'_t, m'_t) dt + Z^{2"}_t dB_t\\
&&\q + \dbE_{\cF_t}\Big[\pa_{xx\mu} V (t, X'_t, m'_t, \tilde X'_t) \pa_z \cH'_{u'}(t, \tilde X'_t, m'_t) + \pa_{\tilde x xx\mu} V (t, X'_t, m'_t, \tilde X'_t) \pa_\g \cH'_{u'}(t, \tilde X'_t, m'_t) \Big]dt.\nonumber
\eea
On the other hand, differentiate \reff{dY1-2} further with respect to $x$, we have
\bea
\label{dY2-2}
&&\dis \partial_{xxt} V+ \Big[\pa_{xx}\cH_V +\pa_{zz}\cH_V |\pa_{xx}V|^2 + \pa_{\g\g}\cH_V |\pa_{xxx} V|^2+ \pa_z\cH_V \pa_{xxx} V + \pa_\g\cH_V \pa_{xxxx} V \nonumber\\
&&\dis \qq\qq\q + 2\pa_{xz}\cH_V \pa_{xx} V + 2\pa_{x\g}\cH_V \pa_{xxx} V + 2 \pa_{z\g}\cH_V \pa_{xx} V\pa_{xxx} V\Big](t, x,\mu)\\
&&\dis\q +\int_{\mathbb R^d}\Big[\partial_{xx\mu}V(t,x,\mu,\tilde x) \pa_z\cH_V( t, \tilde x,\mu)+\partial_{\tilde x xx\mu}V(t,x,\mu,\tilde x) \pa_\g\cH_V( t, \tilde x,\mu)\Big]\mu(d\tilde x)=0.\nonumber
\eea
Plug \reff{dY2-2} into \reff{dY2-1}, and note that $Z^{2'}_t \slash \sqrt{2\pa_\g H'(\Xi'_t)} = \pa_{xxx} u'(t, X'_t)$,  by \reff{FBSDE} we have
\bea
\label{dY2-3}
&&\dis d  \D Y^{2}_t  =  - \G^2_t dt + \D Z^2_t dB_t,\q\mbox{where}\\
&&\dis \G^2_t := \D\pa_{xx}\cH + \D\pa_z\cH \pa_{xxx} V + \D\pa_\g\cH \pa_{xxxx} V+[\pa_{zz}\cH_V |\pa_{xx}V|^2 - \pa_{zz}\cH'_{u'} |\pa_{xx}u'|^2]\nonumber\\
&&\dis\q  +  [\pa_{\g\g}\cH_V |\pa_{xxx} V|^2- \pa_{\g\g}\cH_{u'} |\pa_{xxx} u'|^2] + 2[\pa_{xz}\cH_V \pa_{xx} V-\pa_{xz}\cH'_{u'}\pa_{xx} u']  \nonumber\\
&&\dis \q  + 2[\pa_{x\g}\cH_V \pa_{xxx} V- \pa_{x\g}\cH_{u'} \pa_{xxx} u']  + 2 [\pa_{z\g}\cH_V \pa_{xx} V\pa_{xxx} V-\pa_{z\g}\cH_{u'} \pa_{xx} u'\pa_{xxx} u'\Big](t, X'_t, m'_t) \nonumber\\
&&\dis\q + \dbE_{\cF_t}\Big[\pa_{xx\mu} V (t, X'_t, m'_t, \tilde X'_t) \D\pa_z\cH( t, \tilde X'_t, m'_t) + \pa_{\tilde x xx\mu} V (t, X'_t, m'_t, \tilde X'_t) \D\pa_\g\cH( t, \tilde X'_t, m'_t) \Big].\nonumber
\eea
Similar to \reff{dY0-4}, we have
\bea
\label{dY2-4}
\big|\D\pa_{xx}\cH (t, X'_t, m'_t)\big| \le C\Big[\|\pa_{xx} \D H\|_R+|\D Y^1_t| + |\D Y^2_t|\Big].
\eea
%and similar estimates hold for $\D\pa_z\cH$ and $\D\pa_\g\cH$. 
Moreover, note that
\beaa
&&\dis\Big|\big(\pa_{xx}V - \pa_{xx}u'\big)(t, X'_t, m'_t)\Big| = |\D Y^2_t|;\\
&&\dis  \Big|\big(\pa_{xxx}V - \pa_{xxx}u'\big)(t, X'_t, m'_t)\Big| =  |\D Z^{2}_t|\slash \sqrt{2\pa_\g H'(\Xi'_t)}\le C|\D Z^2_t|,
 \eeaa
where $C$ depends on $c_0$ in Assumption \ref{assum-regHG}-(iii). Then
\beaa
&&\dis \Big|\pa_{zz}\cH_V |\pa_{xx}V|^2 - \pa_{zz}\cH'_{u'} |\pa_{xx}u'|^2\Big|(t, X'_t, m'_t) \\
&&\dis \q \le \Big[|\D\pa_{zz}\cH| |\pa_{xx}V|^2 + |\pa_{zz}\cH_V | [|\pa_{xx} V|+|\pa_{xx}u'|] |\D Y^2_t|\Big](t, X'_t, m'_t)\\
&&\dis\q  \le C\Big[|\D\pa_{zz}\cH( t, X'_t, m'_t)| + |\D Y^2_t|\Big] \le C\Big[\|\pa_{zz}\D H\|_R+|\D Y^1_t| + |\D Y^2_t|\Big];\\
&&\dis  \Big|\pa_{\g\g}\cH_V |\pa_{xxx} V|^2- \pa_{\g\g}\cH_{u'} |\pa_{xxx} u'|^2\Big|(t, X'_t, m'_t) \le C\Big[\|\pa_{\g\g} \D H\|_R+|\D Y^1_t| + |\D Y^2_t| + |\D Z^2_t|\Big].
\eeaa
We can get the estimates for the other terms in \reff{dY2-3} similarly. Then
\bea
\label{dY2-5}
| \G^{2}_t | \le C \Big[\sum_{k=1}^2 \|\pa^{(k)}_{(x,z,\g)}\D H\|_R  + |\D Y^1_t| + |\D Y^2_t| + |\D Z^2_t|+ \dbE_{\cF_t}\big[|\D \tilde Y^1_t| + |\D \tilde Y^2_t|\big]\Big].
\eea
Note that 
\beaa
|\D Y^k_T| = |\pa^{(k)}_x \D G(X'_T, m'_T)|\le \|\pa^{(k)}_x\D G\|,\q k=0,1,2.
\eeaa
 Consider the BSDEs \reff{dY0-3}, \reff{dY1-3} and \reff{dY2-3}. By  \reff{dY0-5}, \reff{dY1-5}, \reff{dY2-5} and standard BSDE estimates (cf. \cite{Zhang}) we have
\beaa
\sum_{k=0}^2 |\D Y^k_0|  \le C\sum_{k=0}^2\big[ \|\pa^{(k)}_{(x,z,\g)}\D H\|_R +\|\pa^{(k)}_{x}\D G\|\big],\q\mbox{a.s.}
\eeaa
where the $C$ may depend on $T$. Note that $\D Y^k_0 =  \pa^{(k)}_x \D u(0, \xi)$. Then the above implies \reff{stability1}.

Finally, to see \reff{stability2}, we first note that \reff{stability1} implies, for any $\mu\in \cP_2(\dbR)$,
\beaa
%\label{dY0-2}
\sum_{k=0}^2 |\pa^{(k)}_x \D V(0,x,\mu)| \le C \sum_{k=0}^2\big[ \|\pa^{(k)}_{(x,z,\g)}\D H\|_R +\|\pa^{(k)}_{x}\D G\|\big],\mbox{for $\mu$-a.e. $x$}.
 \eeaa
Since $\pa^{(k)}_x \D V$ are continuous, and any $\mu\in \cP_2(\dbR)$ can be approximated by measures with full support. Then the above holds true for all $x\in \dbR$, and thus we prove \reff{stability2} for $t=0$. The estimate for general $t$ follows the same arguments, by considering the MFG on $[t, T]$. 
\qed

\section{The local well-posedness:  Proof of Theorem  \ref{thm-local}.}
\label{sect-local}
\setcounter{equation}{0}
We prove the theorem in several steps, reported in three subsections. The key is a pointwise representation formula for $\pa^{(k)}_x\pa_\mu V$ in Subsection \ref{sect-point}. 
%We shall continue to use the notations in the previous section. 

\subsection{The regularity in $x$}
%\label{sect-localx}

We first establish the well-posedness of the MFG system \reff{FBSDE}, which implies that the MFG has a unique equilibrium measure $\bm$.  For this purpose, we introduce a few notations. First, we extend the solution vector $\Xi^\xi$ to $\Pi^\xi := (\Xi^\xi, Z^{\xi,2})$. That is,
\bea
\label{Pi}
\Pi^\xi_t := (X^\xi_t, Y^{\xi,1}_t, Y^{\xi, 2}_t, \cL_{X^\xi_t}, Z^{\xi,2}_t). %\q \Pi^x_t := (X^x_t, Y^{x,1}_t, Y^{x, 2}_t, m_t, Z^{x,2}_t).
\eea
We remark that, once $\Pi^\xi$ is given, the other components of the solution $Y^{\xi,0}, Z^{\xi, 0}, Z^{\xi, 1}$ are obviously uniquely determined. Thus we may call $ \Pi^\xi$ the solution to the FBSDE \reff{FBSDE}. 

Next, for functions $G$, $H$, and their derivatives, we automatically extend their definition to $\pi = (x, y_1, y_2, \mu, z_2)$. Note that here $(y_1, y_2)$  correspond to $(z, \g)$, we shall keep using notations $\pa_z, \pa_\g$ for the derivatives with respect to them. Moreover, in light of  \reff{FBSDE} again, denote: 
\bea
\label{cHk}
\left.\ba{c}
\cH_0 (\pi) := H_3(x, y_1, y_2, \mu),\q \cH_1(\pi) := \pa_xH(x, y_1, y_2, \mu),\q  \cH_2(\pi) := H_4(\pi). %~ \forall\pi=(x, y_1, y_2, \mu, z_2).
%G_k := \pa^{(k)}_x G,~ k=0,1,2. %\q H_1 := \pa_x H,\q H_2 := \sqrt{2\pa_\g H}. 
\ea\right.
\eea
%Here we abuse the notation $\cH_0$ with \reff{propagation-dX2}.
Then FBSDE  \reff{FBSDE} can be simplified to:
\bea
\label{BSDEk}
\begin{cases}
\dis X^\xi_t = \xi +\int_0^t H_1(\Pi^\xi_s) ds + \int_0^t H_2(\Pi^\xi_s) dB_s;\\ 
\dis Y^{\xi,k}_t = \pa^{(k)}_x G(\Pi^\xi_T) + \int_t^T \cH_k(\Pi^\xi_s) ds - \int_t^T Z^{\xi,k}_s dB_s,\q k=0,1,2.
\end{cases}
\eea
We emphasize that, although we write $ \pa^{(k)}_x G(\Pi^\xi_T)$ for notational simplicity, $\pa^{(k)}_x G$ actually depends only on $(X^\xi_T, \cL_{X^\xi_T})$, in particular not on $Y^{\xi,k}_T$, so there is no confusion for the structure of the BSDE. We shall take this convention for all the FBSDEs below.

Note that in Assumption \ref{assum-regHG}-(ii) we assume only local Lipschitz continuity of $H$ and its derivatives,  the well-posedness of \reff{FBSDE} or \reff{BSDEk} does not follow directly from the standard FBSDE estimates. We thus first establish the following well-posedness result.

\begin{lem}
\label{lem-FBSDElocal}
Under the setting in Theorem \ref{thm-local}, the FBSDE \reff{FBSDE} is well-posed when $T\le \d_0$. 
\end{lem} 
\proof We shall solve \reff{FBSDE} by Picard iteration. Fix an arbitrary $\a^0\in \cA$ and denote $\bm^0 := \cL_{X^{\xi,\a^0}}$ as in \reff{Xxia2}. We define $\bm^l$, $l\ge 1$ recursively as follows. Given $\bm^l$, consider the following standard FBSDE with solution $\Pi^{\xi, l}_t=(X^{\xi,l}, Y^{\xi,l, 1}_t, Y^{\xi,l, 2}_t, m_t, Z^{\xi,l, 2}_t)$ and parameter $\bm = \bm^l$: 
\bea
\label{FBSDEl}
\begin{cases}
\dis X_t^{\xi,l} = \xi +\int_0^t H_1(\Pi^{\xi,l}_s) ds + \int_0^t H_2(\Pi^{\xi,l}_s) dB_s;\\ 
\dis Y^{\xi,l, k}_t = \pa^{(k)}_x G(\Pi^{\xi,l}_T) + \int_t^T \cH_k(\Pi^{\xi,l}_s) ds - \int_t^T Z^{\xi,l, k}_s dB_s,\q k=0,1,2.
\end{cases}
\eea
We note that, under Assumption \ref{assum-regu}, this FBSDE is globally well-posed, with the solution $u$ to the HJB equation \reff{HJB} with parameter $\bm = \bm^l$ as its decoupling field. In particular, by \reff{Yxi1}, \reff{Yxi2}, and \reff{eq:L0u}, we have 
\bea
\label{YZbound}
\|Y^{\xi, l, k}\| \le L^*_1,\q \|Z^{\xi, l, k}\|\le L^*_1,\q k=1,2,~ \forall l\ge 0.
\eea

We next define $\bm^{l+1} := \cL_{X^{\xi,l}}$. By \reff{YZbound} and the desired Lipschitz continuity of the coefficients $H_k$, $G_k$, and $\cH_k$, it follows from standard FBSDE arguments (cf. \cite{Antonelli, Zhang}) that 
\bea
\label{Picard}
\sup_{0\le t\le T}W^2_2(m^{l+1}_t, m^l_t) \le \dbE\Big[\sup_{0\le t\le T}|X^{\xi, l+1}_t-X^{\xi,l}|^2\Big] \le c_T \sup_{0\le t\le T}W_2(m^{l}_t, m^{l-1}_t).
\eea
Moreover, when $T\le \d_0$ is small enough, $c_T \le c_{\d_0}<1$. Then it is clear that $\bm^* := \lim_{l\to \infty}\bm^l$ exists, and the solution to \reff{FBSDEl} with parameter $\bm = \bm^*$ is the unique solution to \reff{FBSDE}. 
\qed

\ms
We now fix $T\le \d_0$ and let $\bm$ be the unique equilibrium measure obtained at above.  For each $x\in \dbR$, we consider the following standard FBSDE with parameter $\bm$ and solution $\Pi^x_t = (X^x_t, Y^{x,1}_t, Y^{x,2}_t, m_t, Z^{x,2}_t)$:
\bea
\label{FBSDEx}
\begin{cases}
\dis X^x_t = x +\int_0^t H_1(\Pi^x_s) ds + \int_0^t H_2(\Pi^x_s) dB_s;\\ 
\dis Y^{x,k}_t = \pa^{(k)}_x G(\Pi^x_T) + \int_t^T \cH_k(\Pi^x_s) ds - \int_t^T Z^{x,k}_s dB_s,\q k=0,1,2.
\end{cases}
\eea
We emphasize that, in $\Pi^x$ the measure variable is still $m_t = \cL_{X^\xi_t}$, rather than $\cL_{X^x_t}$. Moreover, as in \reff{BSDEk}, here $\pa^{(k)}_x G$ actually depends only on $(X^x_T, m_T)$. 

Clearly this FBSDE is well-posed with decoupling field $u(\bm; \cd,\cd)\in C^{1,6}_b([0, T]\times \dbR)$, thanks to Assumption \ref{assum-regu}. Define $V(0, x, \mu):= u(\bm; 0,x)$ as in \reff{MFGV}. The following result is immediate. 

\begin{lem}
\label{lem-Vregx}
Under the setting in Theorem \ref{thm-local}, we have
\bea
\label{Vxrep}
\left.\ba{c}
\dis \pa^{(k)}_xV(0,x,\mu) = Y^{x,k}_0,~k=0,1,2;~ %\q \pa_x V(0,x,\mu) = Y^{x,1}_0,\q \pa_{xx} V(0,x,\mu) = Y^{x,2}_0,\ms\\
 \dis V(0,\cd, \mu) \in C^6_b(\dbR);~ \|\pa^{(k)}_x V(0,\cd,\mu)\| \le L^*_1,~ k=1,2,3,4.
\ea\right.
\eea
%Moreover, by the stability of the involved FBSDEs, $\pa^{(k)}_x V(0,\cd,\cd)$ is continuous in $(x,\mu)$, $0\le k\le 6$.
\end{lem} 

\subsection{The representation formula for $\pa_\mu V$}
 \label{sect-point}

In this subsection we provide pointwise representation formulae for $\pa_x^{(k)}\pa_\mu V$, $k=0,1,2$, which implies that $\pa_x^{(k)}\pa_\mu V$ exist and are continuous. We follow the approach in \cite{MZ2, GMMZ}, however, our proofs here are easier. Fix $\mu\in \cP_2(\dbR)$, $\xi\in \dbL^2(\cF_0,\mu)$, and $\bm:= \cL_{X^\xi}$.

For a function $\Phi$ with variable $\pi$, the following differential operators will be used frequently: 
given $\d \pi= (\d x, \d y_1, \d y_2, \d z_2)$, and $\tilde x_1,\cd, \tilde x_m, \delta\tilde x_1,\cds, \delta\tilde x_m$ for some $m\ge 1$,
 \bea
 \label{operator1}
\left.\ba{lll} %\begin{array}{ll}
 \dis \nabla_x \Phi(\pi;\delta \pi):=\pa_x\Phi (\pi)\delta x+\pa_{y_1}\Phi(\pi) \delta y_1+\pa_{y_2}\Phi(\pi)\delta y_2 + \pa_{z_2} \Phi(\pi)\d z_2;\\
\dis  \nabla_\mu\Phi(\pi; \tilde x_1,\cd, \tilde x_m, \delta\tilde x_1,\cds, \delta\tilde x_m):= \sum_{i=1}^m\pa_{\mu}\Phi(\pi,\tilde x_i)\delta\tilde x_i;\\
\nabla \Phi(\pi; \d\pi; \tilde x_1,\cd, \tilde x_m, \delta\tilde x_1,\cds, \delta\tilde x_m):= \nabla_x\Phi(\pi;\delta \pi) + \nabla_\mu\Phi(\pi; \tilde x_1,\cd, \tilde x_m, \delta\tilde x_1,\cds, \delta\tilde x_m).
\ea\right.
\eea
Note that the above operators are linear in $\d\pi$ and $\d \tilde x_i$, consequently  the FBSDEs  driven by them at below are linear.

We now  consider the following linear FBSDE,  by formally differentiating  \reff{FBSDEx} with respect to $x$: denoting   $\td_x\Pi^x_t := (\td_x X^x_t, \td_x Y^{x,1}_t, \td_x Y^{x,2}_t, \td_x Z^{x,2}_t)$, and for $k=0,1,2$,
\begin{eqnarray}\label{tdFBSDEx}
\begin{cases}
\dis \nabla_x X^{x}_t=  1+\int_0^t \nabla_xH_1(\Pi^x_s; \td_x \Pi^x_s) ds+ \int_0^t \nabla_xH_2(\Pi^x_s; \td_x \Pi^x_s)dB_s,\\
\dis \td_x {Y}^{x,k}_t= \td_x (\pa^{(k)}_x G)(\Pi^x_T; \td_x \Pi^x_T) + \int_t^T \td_x \cH_k ( \Pi^x_s; \td_x \Pi^x_s)ds - \int_t^T \td_x {Z}^{x,k}_sdB_s.
\end{cases}
\end{eqnarray}
The next one is motivated from differentiating \reff{FBSDE} with respect to $\xi$, see \reff{tdXeta} below. However, to obtain the pointwise representation, we need to make some modifications. The solution is denoted as $\td_\mu\Pi^{x}_t$ in the above sense and we denote $\Upsilon^{x}_t := (X^x_t, X^\xi_t, \td_x X^x_t, \td_\mu X^{x}_t)$:
\begin{eqnarray}\label{tdFBSDEmu}
\begin{cases}
\dis \nabla_\mu X^{x}_t=  \int_0^t \dbE_{\cF_s}\big[\nabla H_1(\Pi^\xi_s; \td_\mu \Pi^{x}_s; \widetilde \Upsilon_s^{x})\big]ds +\int_0^t \dbE_{\cF_s}\big[\nabla H_2(\Pi^\xi_s; \td_\mu \Pi^{x}_s; \widetilde \Upsilon_s^{x})\big]dB_s,\\
\dis \td_\mu Y^{x,k}_t=\mathbb{E}_{\cF_T}\big[\nabla (\pa^{(k)}_x G)(\Pi^\xi_T; \td_\mu \Pi^{x}_T; \widetilde \Upsilon_T^{x})\big] +\int_t^T\dbE_{\cF_s}\big[\nabla \cH_k(\Pi^\xi_s; \td_\mu \Pi^{x}_s; \widetilde \Upsilon_s^{x})\big]ds \\
\dis\qq\qq  - \int_t^T \td_\mu Z^{x,k}_sdB_s,\q  k=0,1,2. 
\end{cases}
\end{eqnarray}
Recall \reff{operator1} that the driver operator $\td$ here is the sum $\td_x + \td_\mu$, and the state process is $\Pi^\xi$, instead of $\Pi^x$.   We also note that $\widetilde  \Upsilon^{x}$ involves the law of the unknown $\td_\mu X^{x}$, thus this is a linear McKean-Vlasov FBSDE.

The third equation is a modification of \reff{tdFBSDEmu} with solution denoted as $\td_\mu\Pi^{x, \tilde x}$, by changing the state process from $\Pi^\xi$ to $\Pi^x$ and replacing the $\Upsilon^x$ with $\Upsilon^{\tilde x}$: 
\begin{eqnarray}\label{tdFBSDEmux}
\begin{cases}
\dis \nabla_\mu X^{x, \tilde x}_t=  \int_0^t\dbE_{\cF_s}\big[\nabla H_1(\Pi^x_s;\td_\mu \Pi^{x, \tilde x}_s; \widetilde \Upsilon_s^{\tilde x})\big]ds +\int_0^t \dbE_{\cF_s}\big[\nabla H_2(\Pi^x_s; \td_\mu \Pi^{x, \tilde x}_s; \widetilde \Upsilon_s^{\tilde x})\big]dB_s,\\
\dis \td_\mu Y^{x, \tilde x, k}_t=\mathbb{E}_{\cF_T}\big[\nabla (\pa^{(k)}_x G)(\Pi^x_T; \td_\mu \Pi^{x, \tilde x}_T; \widetilde \Upsilon_T^{\tilde x})\big]+\int_t^T \dbE_{\cF_s}\big[\nabla\cH_k(\Pi^x_s; \td_\mu \Pi^{x, \tilde x}_s; \widetilde \Upsilon_s^{\tilde x})\big]ds\\
\dis \qq\qq  - \int_t^T \td_\mu Z^{x,\tilde x, k}_sdB_s,\q k=0,1,2.
\end{cases}
\end{eqnarray}
We remark that, since $\Upsilon$ does not involve the unknown $\td_\mu X^{x,\tilde x}$, this is a standard FBSDE.  

\begin{thm}
\label{thm-rep}
Under the assumptions of Theorem \ref{thm-local}, the FBSDEs \reff{tdFBSDEx}, \reff{tdFBSDEmu}, and \reff{tdFBSDEmux} are well-posed on $[0,T]$ whenever $T\le \d_0$, and  the following representation formulae hold:
    \bea
    \label{Vmurep}
    \pa^{(k)}_x\pa_\mu  V(0,x,\mu, \tilde x) = \td_\mu Y^{x, \tilde x, k}_0,\q k=0,1,2.
    \eea
\end{thm}
\proof  Since $T\le \d_0$, the well-posedness of the FBSDEs follow from standard FBSDE theory and its extension to McKean-Vlasov equations. We now prove \reff{Vmurep} in three steps.  

\ms
\no\emph{Step 1.} For any $\xi\in \mathbb{L}^2(\mathcal{F}_0;\mu)$ and $\eta\in \mathbb{L}^2(\mathcal{F}_0)$, following standard arguments and by the stability property of the involved FBSDE systems we have
\begin{eqnarray}
\label{Xetaconv}
\left.\ba{c}
\dis\lim_{\varepsilon\rightarrow0}\mathbb{E}\Big[\sup_{0\leq t\leq T}\Big|\frac{1}{\varepsilon}[X^{\xi+\varepsilon\eta}_t-X^{\xi}_t]
-\delta X^{\eta}_t\Big|^2\Big]=0;\ms\\
\dis \lim_{\varepsilon\rightarrow0}\mathbb{E}\Big[\sup_{0\leq t\leq T}\sum_{k=0}^2\Big|\frac{1}{\varepsilon}[{Y}^{\xi+\varepsilon\eta,x,k}_t-{Y}^{\xi,x,k}_t]
-\delta Y^{x,\eta,k}_t\Big|^2\Big]=0;
\ea\right.
\end{eqnarray}
where $\d \Pi^{\eta}$  satisfies the following linear McKean-Vlasov FBSDEs, obtained by formally differentiating \reff{FBSDE} with respect to $\xi$ along the direction $\eta$:
\begin{eqnarray}\label{tdXeta}
\begin{cases}
\dis \delta X^{\eta}_t=\eta+\int^t_{0}\dbE_{\cF_s}\big[\nabla H_1(\Pi^\xi_s, \d \Pi^{\eta}_s; \tilde X^\xi_s, \widetilde{\d X_s^{\eta}})\big]ds +\int^t_{0} \dbE_{\cF_s}\big[\nabla H_2(\Pi^\xi_s, \d \Pi^{\eta}_s; \tilde X^\xi_s, \widetilde{\d X_s^{\eta}})\big]dB_s,\\
\dis \delta {Y}^{\eta,k}_t=\mathbb{E}_{\cF_T}\big[\td (\pa^{(k)}_x G)(\Pi^\xi_T, \d \Pi^{\eta}_T; \tilde X^\xi_T, \widetilde{\d X_T^{\eta}})\big] +\int^{T}_{t}\dbE_{\cF_s}\big[\nabla\cH_k(\Pi^\xi_s, \d \Pi^{\eta}_s; \tilde X^\xi_s, \widetilde{\d X_s^{\eta}})\big]ds\\
\dis\qq\qq -\int^T_s\delta {Z}^{\eta,k}_sdB_s, \q k=0,1,2.
\end{cases}
\end{eqnarray}
and $\d \Pi^{x, \eta}$ satisfies the following linear FBSDE, obtained by formally differentiating \reff{FBSDEx} in $\xi$:
%by replacing $\Pi^\xi$ with $\Pi^x$ and with initial $X^{x,\eta}_0=0$:
\begin{eqnarray}\label{tdXxeta}
\begin{cases}
\dis \delta X^{x,\eta}_t=\int^t_{0}\dbE_{\cF_s}\big[\nabla H_1(\Pi^x_s, \d \Pi^{x,\eta}_s; \tilde X^\xi_s, \widetilde{\d X_s^{\eta}})\big]ds +\int^t_{0} \dbE_{\cF_s}\big[\nabla H_2(\Pi^x_s, \d \Pi^{x,\eta}_s; \tilde X^\xi_s, \widetilde{\d X_s^{\eta}})\big]dB_s,\\
\dis \delta {Y}^{x,\eta,k}_t=\mathbb{E}_{\cF_T}\big[\td (\pa^{(k)}_x G)(\Pi^x_T, \d \Pi^{x,\eta}_T; \tilde X^\xi_T, \widetilde{\d X_T^{\eta}})\big] +\int^{T}_{t}\dbE_{\cF_s}\big[\nabla\cH_k(\Pi^x_s, \d \Pi^{x,\eta}_s; \tilde X^\xi_s, \widetilde{\d X_s^{\eta}})\big]ds\\
\dis\qq\qq -\int^T_s\delta {Z}^{x,\eta,k}_sdB_s, \q k=0,1,2.
\end{cases}
\end{eqnarray}
In particular, \reff{Xetaconv} and \reff{Vxrep} imply
\beaa
%\label{08136}
\lim_{\varepsilon\rightarrow0}\left|\frac{1}{\varepsilon}[\pa^{(k)}_xV(0,x,\mathcal{L}_{\xi+\varepsilon\eta})-\pa^{(k)}_xV(0,x,\mathcal{L}_{\xi})]
-\delta {Y}^{x,\eta,k}_0\right|^2=0,\q k=0,1,2.
\eeaa
Thus, by the definition of $\partial_{\mu}V$ in \reff{pamu},
\begin{eqnarray}\label{Vmurep2}
 \mathbb{E}\big[\partial_{\mu} \pa^{(k)}_xV(0,x,\mu,\xi)\eta\big]=\delta Y^{x,\eta,k}_0,  \q k=0,1,2.
\end{eqnarray}

\emph{Step 2.} In this step we assume $\xi$  is discrete: $p_i=\mathbb{P}(\xi=x_i)$, $i=1,\cdots, n$, and we shall prove
\bea
\label{Vmurepclaim}
\d \Pi^{x_i, \eta_j} = p_j  \td_\mu \Pi^{x_i, x_j}_0, \q\mbox{where}\q \eta_j :=  {\bf 1}_{\{\xi=x_j\}}, \q  i, j=1,\cds, n.
\eea
Then, by \reff{Vmurep2} we have
\bea
\label{Vmurep3}
\delta Y^{x_i,\eta_j,k}_0 = \mathbb{E}\big[\partial_{\mu} \pa^{(k)}_xV(0,x_i,\mu,\xi)\eta_j\big] = p_j \partial_{\mu} \pa^{(k)}_xV(0,x_i,\mu,x_j),  \q k=0,1,2.
\eea
 Combining with \reff{Vmurepclaim}, this implies that, when $\mu=\cL_\xi$ is discrete and for $k=0,1,2$,
\bea
\label{Vmurep4}
\partial_{\mu} \pa^{(k)}_xV(0,x_i,\mu,x_j) = \td_\mu Y^{x_i, x_j, k}_0,~\mbox{or say},~\partial_{\mu} \pa^{(k)}_xV(0,x,\mu,\tilde x) = \td_\mu Y^{x, \tilde x, k}_0,~ \mu\mbox{-a.e.}~x, \tilde x.
\eea

To see \reff{Vmurepclaim}, we first note that \reff{FBSDEx} is a standard FBSDE with parameter $\bm = \cL_{X^\xi}$. By \reff{FBSDE} and \reff{FBSDEx}, one can easily see that
\bea
\label{Pixidecompose1}
\Pi^\xi = \sum_{i=1}^n \eta_i \Pi^{x_i}.
\eea
 This implies that, since $\{\xi=x_i\}_{i=1,\cds, n}$ form a partition of $\O$,  for $k=0,1,2$,
\bea
\label{etapartition}
\left.\ba{lll}
\dis \eta_i \nabla_x\cH_k(\Pi^\xi_s, \td_x \Pi^{x_i}_s) = \sum_{j=1}^n \eta_i\eta_j \nabla_x\cH_k(\Pi^{x_j}_s, \td_x \Pi^{x_i}_s) =\eta_i\nabla_x\cH_k(\Pi^{x_i}_s, \td_x \Pi^{x_i}_s);\\
\dis  \dbE_{\cF_s}\big[\pa_\mu \cH_k(\Pi^\xi_s, \tilde X^\xi_s)\widetilde {\eta_i} \widetilde{\td_x X_s^{x_i}}\big] =  \sum_{j=1}^n \dbE_{\cF_s}\big[\pa_\mu \cH_k(\Pi^\xi_s, \tilde X^{x_j}_s)\widetilde {\eta_i} \widetilde {\eta_j} \widetilde{\td_x X_s^{x_i}}\big]\\
\dis = \dbE_{\cF_s}\big[\pa_\mu \cH_k(\Pi^\xi_s, \tilde X^{x_i}_s)\widetilde {\eta_i}  \widetilde{\td_x X_s^{x_i}}\big] = p_i\dbE_{\cF_s}\big[\pa_\mu \cH_k(\Pi^\xi_s, \tilde X^{x_i}_s) \widetilde{\td_x X_s^{x_i}}\big],
\ea\right.
\eea
where the last equality used the fact that $(\Pi^\xi, \tilde X^{x_i}, \widetilde{\td_x X^{x_i}})$ are independent of $\widetilde\eta_i$, and $\dbE[\tilde \eta_i] =p_i$. 

Next, denote $\d \Pi^{'i}:= \eta_i \td_x \Pi^{x_i} + p_i \td_\mu \Pi^{x_i}$. Then, in light of \reff{tdFBSDEx}, \reff{tdFBSDEmu}, and \reff{tdXeta}, by \reff{etapartition}  and \reff{operator1}  we have
\beaa
&&\dis \eta_i \td_x \cH_k ( \Pi^{x_i}_s; \td_x \Pi^{x_i}_s)  + p_i\dbE_{\cF_s}\big[ \td \cH_k ( \Pi^\xi_s; \td_\mu \Pi^{x_i}_s; \widetilde\Upsilon^{x_i}_s)\big]\\
&=&\dis  \eta_i \td_x \cH_k ( \Pi^{\xi}_s; \td_x \Pi^{x_i}_s)  + p_i\td_x \cH_k ( \Pi^{\xi}_s; \td_\mu \Pi^{x_i}_s) \\
&&+ p_i\dbE_{\cF_s}\big[ \pa_\mu \cH_k ( \Pi^\xi_s, \tilde X^{x_i}_s)\widetilde{\td_x X^{x_i}_s} \big]+  p_i\dbE_{\cF_s}\big[\pa_\mu \cH_k ( \Pi^\xi_s, \tilde X^{\xi}_s)\widetilde{\td_\mu X^{x_i}_s}\big]\\
&=& \dis \td_x \cH_k ( \Pi^{\xi}_s; \d \Pi^{'i}_s)  +  \dbE_{\cF_s}\big[\pa_\mu \cH_k(\Pi^\xi_s, \tilde X^\xi_s)\widetilde {\eta_i} \widetilde{\td_x X_s^{x_i}}\big] + p_i\dbE_{\cF_s}\big[\pa_\mu \cH_k ( \Pi^\xi_s, \tilde X^{\xi}_s)\widetilde{\td_\mu X^{x_i}_s}\big]\\
&=&\dis  \td_x \cH_k ( \Pi^{\xi}_s; \d \Pi^{'i}_s)  +  \dbE_{\cF_s}\big[\pa_\mu \cH_k(\Pi^\xi_s, \tilde X^\xi_s)\widetilde{\d X_s^{'i}}\big] \\
&=& \dis   \dbE_{\cF_s}\big[\td \cH_k(\Pi^\xi_s;  \d \Pi^{'i}_s; \tilde X^\xi_s, \widetilde{\d X_s^{'i}})\big].
\eeaa
This implies  that 
\beaa
\dis d \d Y^{'i, k}_s &=& \eta_i d \td_x Y^{x_i, k}_s + p_i d \td_\mu Y^{x_i,k}_s \\
\dis &=& -\Big[\eta_i \td_x \cH_k ( \Pi^{x_i}_s; \td_x \Pi^x_s)  + p_i\dbE_{\cF_s}\big[ \td \cH_k ( \Pi^\xi_s; \td_\mu \Pi^{x_i}_s; \widetilde\Upsilon^{x_i}_s)\big]\Big]ds +  \d Z^{'i, k}_s dB_s\\
\dis &=&  -\dbE_{\cF_s}\big[\td \cH_k(\Pi^\xi_s;  \d \Pi^{'i}_s; \tilde X^\xi_s, \widetilde{\d X_s^{'i}})\big]ds +  \d Z^{'i, k}_s dB_s
\eeaa
Similarly we can show that
\beaa
&&\dis\delta {Y}^{'i,k}_T=\mathbb{E}_{\cF_T}\big[\td (\pa_x^{(k)}G)(\Pi^\xi_T, \d \Pi^{'i}_T; \tilde X^\xi_T, \widetilde{\d X_T^{'i}})\big];\\
&&\dis  d\delta X^{'i}_s=\dbE_{\cF_s}\big[\nabla(\pa_zH)(\Pi^\xi_s, \d \Pi^{'i}_s; \tilde X^\xi_s, \widetilde{\d X_s^{'i}})\big]ds + \dbE_{\cF_s}\big[\nabla(\sqrt{2\pa_\g H})(\Pi^\xi_s, \d \Pi^{'i}_s; \tilde X^\xi_s, \widetilde{\d X_s^{'i}})\big]dB_s.
\eeaa
Moreover, $\delta X^{'i}_0=\eta_i \td_x X^{x_i}_0 + p_i \td_\mu X^{x_i}_0 =  \eta_i$. So $\d \Pi^{'i}$ satisfies \reff{tdXeta} with $\eta=\eta_i$. Then, by the uniqueness of the FBSDE we have $\d \Pi^{'i} = \d \Pi^{\eta_i}$. That is,
\bea
\label{Pixidecompose2}
\d \Pi^{\eta_i} = \eta_i \td_x \Pi^{x_i} + p_i \td_\mu \Pi^{x_i}.
\eea

Finally, in light of \reff{tdXxeta} and \reff{tdFBSDEmux}, by \reff{operator1}, \reff{etapartition}, and \reff{Pixidecompose2} we have, for $k=0,1,2$,
\beaa
&&\dbE_{\cF_s}\big[\nabla\cH_k(\Pi^{x_i}_s, \d \Pi^{x_i,\eta_j}_s; \tilde X^\xi_s, \widetilde{\d X_s^{\eta_j}})\big] \\
&=& \nabla_x \cH_k(\Pi^{x_i}_s, \d \Pi^{x_i,\eta_j}_s) + \dbE_{\cF_s}\Big[\pa_\mu\cH_k(\Pi^{x_i}_s, \tilde X^\xi_s)\tilde \eta_j \widetilde{\td_x X_s^{x_j}}\big] + \dbE_{\cF_s}\Big[\pa_\mu\cH_k(\Pi^{x_i}_s, \tilde X^\xi_s) p_j \widetilde{\td_\mu X_s^{x_j}}\big]\\
&=& \nabla_x \cH_k(\Pi^{x_i}_s, \d \Pi^{x_i,\eta_j}_s) +p_j  \dbE_{\cF_s}\Big[\pa_\mu\cH_k(\Pi^{x_i}_s, \tilde X^{x_j}_s) \widetilde{\td_x X_s^{x_j}}\big] + p_j \dbE_{\cF_s}\Big[\pa_\mu\cH_k(\Pi^{x_i}_s, \tilde X^\xi_s) \widetilde{\td_\mu X_s^{x_j}}\big].
\eeaa
Then, denoting $\d\Pi^{'ij}:= {1\over p_j}\d \Pi^{x_i,\eta_j}$, we have
\beaa
&&\dbE_{\cF_s}\big[\nabla\cH_k(\Pi^{x_i}_s, \d \Pi^{'ij}_s; \tilde X^\xi_s, \widetilde{\d X_s^{\eta_j}})\big] \\
&=& \nabla_x \cH_k(\Pi^{x_i}_s, \d \Pi^{'ij}_s) +  \dbE_{\cF_s}\Big[\pa_\mu\cH_k(\Pi^{x_i}_s, \tilde X^{x_j}_s) \widetilde{\td_x X_s^{x_j}}\big] +  \dbE_{\cF_s}\Big[\pa_\mu\cH_k(\Pi^{x_i}_s, \tilde X^\xi_s) \widetilde{\td_\mu X_s^{x_j}}\big]\\
&=&\dbE_{\cF_s}\big[\nabla\cH_k(\Pi^{x_i}_s; \d \Pi^{'ij}_s; \widetilde \Upsilon_s^{x_j})\big].
\eeaa
Similarly one can check the other required equalities, so that $\d\Pi^{'ij}$ satisfies \reff{tdFBSDEmux} with $(x, \tilde x) = (x_i, x_j)$. Then by the uniqueness of the FBSDE  \reff{tdFBSDEmux}, we obtain $\d\Pi^{'ij} = \td_\mu \Pi^{x_i, x_j}$. This verifies \reff{Vmurepclaim} and hence \reff{Vmurep4} when $\xi$ is discrete.

\ms
\emph{Step 3.} We now prove \reff{Vmurep} in the general case. Fix $\xi\in \dbL^2(\cF_0,\mu)$ and $\eta\in \dbL^2(\cF_0)$. Let $\{x_n\}_{n\ge 1}\subset \dbR$ be a dense sequence. One can easily construct discrete $\xi_n, \eta_n\in \dbL^2(\cF_0)$ such that 

$\bullet$ $\lim_{n\to\infty}\dbE\big[|\xi_n-\xi|^2 + |\eta_n-\eta|^2\big]=0$;

$\bullet$ $\{x_1,\cds, x_n\}\subset \supp(\xi_n + \th \eta_n)$ for all $n\ge 1$ and all $\th\in [0, 1]$.

\no  Clearly $\xi_n + \th \eta_n$ is also discrete, for any $\th\in [0,1]$. Then by the definition of $\pa_\mu$ in \reff{pamu},  it follows from Step 2 that, for any $n\ge m$ and hence $x_m\in \supp(\xi_n + \th \eta_n)$ for all $\th$,
\beaa
&&\dis \pa^{(k)}_xV(0, x_m,  \cL_{\xi_n + \eta_n}) - \pa^{(k)}_xV(0, x_m,  \cL_{\xi_n }) = \int_0^1 \dbE\big[\pa_\mu \pa^{(k)}_xV(0, x_m,  \cL_{\xi_n +\th \eta_n},  \xi_n + \th \eta_n) \eta_n \big] d\th\\
&&\dis = \int_0^1 \dbE\big[J_k(x_m,\cL_{\xi_n +\th \eta_n},  \xi_n + \th \eta_n) \eta_n \big] d\th,\q\mbox{where}\q J_k(x, \mu, \tilde x) := \td_\mu Y^{x, \tilde x,k}_0.
\eeaa
By the stability of the FBSDEs, it is clear that $J_k$ is continuous in all its variables. Since $ \pa^{(k)}_xV$ is also continuous, then, by sending $n\to\infty$ we obtain
\beaa
 \pa^{(k)}_xV(0, x_m,  \cL_{\xi + \eta}) - \pa^{(k)}_xV(0, x_m,  \cL_{\xi}) =  \int_0^1 \dbE\big[J_k(x_m,\cL_{\xi +\th \eta},  \xi + \th \eta) \eta \big] d\th.
\eeaa
Since $\{x_m\}$ is dense, and since both sides above are continuous, we have
\beaa
 \pa^{(k)}_xV(0, x,  \cL_{\xi + \eta}) - \pa^{(k)}_xV(0, x,  \cL_{\xi}) =  \int_0^1 \dbE\big[J_k(x,\cL_{\xi +\th \eta},  \xi + \th \eta) \eta \big] d\th, \q\mbox{for all $x\in \dbR$}.
\eeaa
Again since $J_k$ is continuous, by \reff{pamu} again  this implies that
\beaa
 \pa_\mu\pa^{(k)}_xV(0, x,  \mu, \tilde x) = J_k(x, \mu, \tilde x),\q\mbox{for $\mu$-a.e. $\tilde x$}. 
  \eeaa
Finally, for $\tilde x\notin \supp(\mu)$, since $ \pa_\mu\pa^{(k)}_xV(0, x,  \mu, \tilde x)$ can be arbitrary in the definition \reff{pamu}, we may naturally use $J_k(x, \mu, \tilde x)$ to define it, then we obtain \reff{Vmurep} for all $(x, \mu, \tilde x)$.
\qed

\subsection{Proof of Theorem \ref{thm-local}}
\label{sect-localproof}
By analyzing Definition \ref{defn-space} carefully, the space $\cC^{2,6}(\Th)$ involves the following derivatives:
\bea
\label{C26}
\left.\ba{c}
\dis\pa^{(k)}_x V,~ 1\le k\le 6;\q \pa^{(k)}_x\pa^{(l)}_{\tilde x} \pa_\mu V,~ 0\le k+l\le 5,~ 0\le l\le 3;\q \pa^{(k)}_x \pa_t V,~ 0\le k\le 4;\\
\dis \pa^{(k)}_x\pa^{(l)}_{\tilde x} \pa^{(\l)}_{\bar x} \pa_{\mu\mu} V,~ 0\le k+l+\l\le 4,~0\le l+\l \le 2,~ 0\le k\le 3, ~0\le l \le 2, ~0\le \l \le 1;\\
\dis \pa^{(k)}_x\pa^{(l)}_{\tilde x} \pa_{t\mu} V,~ 0\le k\le 2,~ 0\le l\le 1.
\ea\right.
\eea
We now verify them in three steps. We shall emphasize that the $\d_0$ in Theorem \ref{thm-local} relies on the the regularity of $H$ and $G$ only up to $\cC^{1,4}$, namely the $L_G$ and $L_H(R)$ in Assumption \ref{assum-regHG}.
\ms

\no{\it Step 1.} By Lemma \ref{lem-Vregx} and Theorem \ref{thm-rep}, we see that 
\bea
\label{xreg1}
\pa^{(k)}_x V(0, x,\mu),~ 1\le k\le 6;\q \pa^{(k)}_x\pa_\mu V(0,x,\mu,\tilde x),~0\le k\le 2.
\eea
 By the stability of the involved FBSDEs, it is clear that these derivatives are continuous in $(x,\mu)$ and in $(x, \mu, \tilde x)$, respectively. Moreover, by differentiating \reff{tdFBSDEx}, \reff{tdFBSDEmu}, and \reff{tdFBSDEmux} formally in $(x,\tilde x)$, and by the assumed regularity in Assumption \ref{assum-regHG}, one can easily obtain representation formulae for the following higher order derivatives at $(0,x,\mu)$ or at $(0,x,\mu, \tilde x)$: 
\bea
\label{xreg}
\pa^{(k)}_x V,~  1\le k\le 6;\qq \pa^{(k)}_x\pa^{(l)}_{\tilde x}\pa_\mu V,~0\le k+l \le 5, ~0\le l\le 3. 
\eea
In particular, we may use the same $\d_0$ relying on $L_G$ and $L_H(R)$ with $R=L^*_1$, but not on the higher order derivatives of $H$ and $G$. Indeed, consider $\pa^{(3)}_x\pa_\mu V$ as an example. Differentiate \reff{tdFBSDEmu} in $x$ and denote the solution as $\td^2_\mu \Pi^x$. Note that \reff{tdFBSDEmu} is linear in $\td_\mu \Pi^x$, then the forward SDE for  $\nabla^2_\mu X^{x}$ becomes:
\begin{eqnarray}\label{tdFBSDEmu2}
\left.\ba{c}
\dis \nabla^2_\mu X^{x}_t=  \int_0^t \Big[\nabla_x H_1(\Pi^\xi_s; \td^2_\mu \Pi^{x}_s) + \dbE_{\cF_s}\big[\pa_\mu H_1(\Pi^\xi; \tilde X^\xi_s) \widetilde{\td^2_\mu X^x_s} + \widetilde\G_s\big] \Big]ds +\int_0^t \big[\cds\big]dB_s\ms\\
\dis\mbox{where}\q \widetilde\G_s:=  \pa_{x\mu} H_1(\Pi^\xi; \tilde X^x_s) |\widetilde{\td_x X^x_s}|^2+ \pa_\mu H_1(\Pi^\xi; \tilde X^x_s) \widetilde{\td^2_x X^x_s},
\ea\right.
\end{eqnarray}
$\td^2_x X^x_s$ is the derivative of \reff{tdFBSDEx} with respect to $x$, which is well-posed due to Assumption \ref{assum-regu}, and the diffusion term has the same structure. The BSDEs for $ \td^2_\mu Y^{x,k}$ also have the same structure. Note that $\d_0$ depends only on the bounds of the coefficients for the solution $\nabla^2_\mu X^{x}$, namely the bounds of $\pa_\pi H_1$ and $\pa_\mu H_1$, but not on that $\tilde \G$. With the similar dependence on the coefficients $H_2, \cH_k$ and $\pa^{(k)}_x G$, we see that for the same $\d_0$, the resulting FBSDEs for $\nabla^2_\mu \Pi^{x}_t$ remain well-posed. Similarly differentiating \reff{tdFBSDEmux} in $x$, we see that the resulting FBSDE with solution $\td^2_\mu X^{x,\tilde x}$ is also well-posed whenever $T\le \d_0$. Thus $\pa^{(3)}_x\pa_\mu V(0,x,\mu,\tilde x) = \td^2_\mu Y^{x,\tilde x, 2}_0$ exists. 
Similarly we can show that all the derivatives in \reff{xreg} exist. Moreover, by the estimates and stability of the involved FBSDEs, these derivatives are bounded and continuous in $(x,\mu, \tilde x)$. 

Similarly, for each $t\in [0, T]$, we may define $V(t, \cd,\cd)$ such that all the derivatives in \reff{xreg} exist and are bounded and continuous. By the stability of the involved FBSDEs,  these derivatives, as well as $V$ itself, are also continuous in $t$.

\ms
\no{\it Step 2.} We next study $\pa_t V$.  Note that, for $u = u(\bm;\cd,\cd)$, 
\beaa
&&V(t, x, \mu) - V(0, x, \mu) = \dbE\Big[[V(t, X^x_t, m_t) - V(0, x, \mu)] - [V(t, X^x_t, m_t)-V(t, x, \mu)]\Big]\\
&&=  \dbE\Big[[u(t, X^x_t) - u(0, x)] - [V(t, X^x_t, m_t)-V(t, x, \mu)]\Big]\\
&& = \dbE\Big[\int_0^t \big[\pa_t u(s, X^x_s) + \pa_x u(s, X^x_s) H_1(\Pi^x_s) + {1\over 2} \pa_{xx} u(s, X^x_s) |H_2(\Pi^x_s)|^2\big]ds\\
&&\dis\qq - \int_0^t \big[\pa_x V(s, X^x_s, m_s) H_1(\Pi^x_s) + {1\over 2} \pa_{xx} V(s, X^x_s, m_s) |H_2(\Pi^x_s)|^2\big] ds\\
&&\dis\qq - \int_0^t \big[\pa_\mu V(s, X^x_s, m_s, \tilde X^\xi_s) H_1(\tilde \Pi^\xi_s) + {1\over 2} \pa_{\tilde x\mu} V(s, X^x_s, m_s, \tilde X^\xi_s) |H_2(\tilde \Pi^\xi_s)|^2\big] ds\Big].
\eeaa
Then clearly $\pa_t V(0, x, \mu)$ exists: denoting $\pi = \Pi^x_0$ and  noting that $m_0 = \mu$, 
\beaa
\pa_t V(0, x, \mu) &=&\big[\pa_t u(0, x) + \pa_x u(0,x) H_1(\pi) + {1\over 2} \pa_{xx} u(0,x) |H_2(\pi)|^2\big] \\
&&-  \big[\pa_x V(0, x, \mu) H_1(\pi) + {1\over 2} \pa_{xx} V(0,x,\mu) |H_2(\pi)|^2\big]\\
&& - \dbE\big[\pa_\mu V(0, x, \mu, \tilde \xi) \pa_z \cH_V(0, \tilde \xi, \mu) +  \pa_{\tilde x\mu} V(0,x,\mu, \tilde\xi) \pa_\g \cH_V(0, \tilde \xi, \mu)\big]\\
&=&\pa_t u(0, x)  - \dbE\big[\pa_\mu V(0, x, \mu, \tilde \xi) \pa_z \cH_V(0, \tilde \xi, \mu) +  \pa_{\tilde x\mu} V(0,x,\mu, \tilde\xi) \pa_\g \cH_V(0, \tilde \xi, \mu)\big],
\eeaa
where the last equality used the fact $\pa^{(k)}_x u(0,x) = \pa^{(k)}_x V(0,x,\mu)$, $k=1,2$. Note further that by \reff{HJB} we have $\pa_t u = -\cH_V$. Then the above equation exactly implies \reff{master} at $(0,x,\mu)$.

 Similarly, for any $t\in [0, T)$, the right time derivative $\pa_{t+} V(t,x,\mu)$ exists and satisfies the master equation \reff{master}. Since $\pa^{(k)}_x  V$ and $\pa_{\tilde x\mu} V$ are continuous, this implies that $\pa_{t+} V$ is continuous. Then the time derivative $\pa_t V$ exists and is continuous, and therefore, $V\in \cC^{1,2}_b(\Th)$ and satisfies the master equation \reff{master}. Moreover, \reff{master} implies
\bea
\label{patV}
\partial_t V= - H(x,\partial_x V,\partial_{xx}V,\mu)- \mathcal{N}V(t,x,\mu).
\eea
Then, by \reff{xreg} it is clear that the following derivatives exist and are bounded and continuous:
\bea
\label{treg}
\pa^{(k)}_x \pa_t V,~ 0\le k\le 4.
\eea

\ms
\no{\it Step 3.} We now investigate $\pa_{\mu\mu} V$ by the bootstrap arguments. The main idea is that, by differentiating \reff{master} in $\mu$ twice and by viewing the derivatives in \reff{xreg} as given, we may obtain a linear master equation for $\pa_{\mu\mu} V$, which leads to a representation formula for $\pa_{\mu\mu} V$ immediately.  Since the notation is rather heavy, in particular it will involve some notational convention we will use for the proof of  \reff{propagation-claim} in Appendix, we postpone this analysis to Appendix as well.

Once we obtain the representation for $\pa_{\mu\mu} V$, following similar arguments as in Step 1, for the same $\d_0$ we may obtain representation formula for the following derivatives which are also bounded and continuous:
\bea
\label{ureg}
\pa_x^{(k)} \pa_{\mu\mu} V,~ \pa_x^{(k)} \pa_{\tilde x\mu\mu} V,~ \pa_x^{(k)} \pa_{\tilde x \bar x\mu\mu} V,\q  k=0,1,2.
\eea
Finally, using \reff{patV} again,  by \reff{xreg} and \reff{ureg} we see that the following derivatives exist and are bounded and continuous:
\bea
\label{treg2}
 \pa^{(k)}_x\pa^{(l)}_{\tilde x} \pa_{t\mu} V,~ 0\le k\le 2,~ 0\le l\le 1.
\eea
Combining \reff{xreg}, \reff{treg}, \reff{ureg}, and \reff{treg2}, we verify \reff{C26} and thus $V\in \cC^{2,6}_b(\Th)$.
\qed

\begin{rem}
\label{rem-bootstrap}
The bootstrap arguments in Step 3 above allows us to increase the regularity of $V$ to any order, provided the corresponding sufficient regularity of the data.
\end{rem}

  \section{Appendix}
\label{sect-appendix}
\setcounter{equation}{0}

\no{\bf Proof of \reff{LL1}-\reff{LL2}-\reff{DP1}-\reff{DP2}.} Clearly \reff{LL1} is a special case of \reff{LL2} with $\phi_1=\phi$, $\phi_2 =0$, and similarly \reff{DP1} is a special case of \reff{DP2}, so it suffices to prove \reff{LL2} and \reff{DP2}. 

We first prove \reff{LL2}. Recall \reff{Xvolatility} and note that  $\mu_t:=\cL_{X_t}$ solves the Fokker-Planck equation \reff{Fokker-Planck} in the weak sense, that is, for any test function $\varphi\in C_b^2(\dbR)$, 
\beaa
%\label{eq:FK}
\int_{\dbR}\varphi(x)\big(\mu_t(dx)-\mu_0(dx)\big)=\int_0^t\int_{\dbR}\big[\phi_1(x) \pa_x\varphi(x)+{1\over 2}|\phi_2(x)|^2\pa_{xx}\varphi(x)\big] \mu_s(dx)ds.
\eeaa
Then, recalling the fact that $\pa_{\tilde x} {\d\over \d\mu}U(x, \mu, \tilde x) = \pa_\mu U(x, \mu, \tilde x)$, we have
\beaa
&&\dis \int_{\dbR^2}  {\d\over \d\mu} U(x, \mu, \tilde x) (\mu_t-\mu)(d\tilde x)(\mu_t-\mu)(dx) \\
&&\dis =\int_\dbR \Big[ \int_0^t\int_{\dbR}\big[\phi_1(\tilde x) \pa_\mu U(x, \mu, \tilde x) +{1\over 2}|\phi_2(\tilde x)|^2\pa_{\tilde x\mu} U(x, \mu, \tilde x)\big]  \mu_s(d\tilde x)ds\Big] (\mu_t-\mu)(dx) \\
&&\dis = \int_0^t\int_0^t\int_{\dbR^2}[\phi_1(x),  {1\over 2}|\phi_2(x)|^2] \ccD_\mu U(x, \mu, \tilde x)[\phi_1(\tilde x),  {1\over 2}|\phi_2(\tilde x)|^2]^\top \mu_{\tilde s}(d\tilde x)\mu_s(dx)  d\tilde s ds.
\eeaa
Since $\mu_t$ is continuous in $t$, we have
\beaa
&&\dis \lim_{t\to 0}{1\over t^2}\int_{\dbR^2}  {\d\over \d\mu} U(x, \mu, \tilde x) (\mu_t-\mu)(d\tilde x)(\mu_t-\mu)(dx) \\
&&\dis =\int_{\dbR^2}[\phi_1(x),  {1\over 2}|\phi_2(x)|^2] \ccD_\mu U(x, \mu, \tilde x)[\phi_1(\tilde x),  {1\over 2}|\phi_2(\tilde x)|^2]^\top \mu(d\tilde x)\mu(dx) .
\eeaa
This implies \reff{LL2} immediately. 

We now prove \reff{DP2}. Recall the $X$ in \reff{Xvolatility2} and $\pa_\mu \cU(\mu, x) = \pa_x U(x, \mu)$. Applying It\^{o}'s formula \reff{Ito} we have
\beaa
{d\over dt} \cU\big(\cL_{X_t}\big)  = \dbE\Big[\pa_x U(X_t, \cL_{X_t}) \eta_1 + {1\over 2}\pa_{xx}  U(X_t, \cL_{X_t}) |\eta_2|^2\Big].
\eeaa
Applying  \reff{Ito} again on the right side above, by straightforward calculation we obtain: 
\beaa
 {d^2\over dt^2} \cU\big(\cL_{X_t}\big)  = \dbE\Big[ (\eta_1, \eta_2') \ccD_{\mu}U(X_t,\cL_{X_t},\tilde X_t)(\tilde\eta_1, \tilde \eta_2')^\top + (\eta_1, \eta_2') \ccD_{x}U(X_t,\cL_{X_t})(\eta_1, \eta_2')^\top\Big],
\eeaa
where $\eta_2' := {1\over 2} |\eta_2|^2$. Setting $t=0$ we obtain \reff{DP2}.
\qed

\bs

\no{\bf Proof of Example \ref{eg-DP}.} Note that
\beaa
&\dis \pa_{x\mu} U(x, \mu, \tilde x) = -\cos (x) \sin(\tilde x),\q \pa_{xx} U(x, \mu) = - \sin(x) \int_\dbR \cos(x') \mu(dx') + \k;\\
&\dis \pa_{\tilde x xx\mu} U(x, \mu, \tilde x) = \sin (x) \cos(\tilde x),\q \pa_{xxxx} U(x, \mu) = \sin(x) \int_\dbR \cos(x') \mu(dx').
\eeaa
For any $\xi, \eta\in \dbL^2(\cF_T)$, it is straightforward to show that
\beaa
 \dis \cM_{DP}^1 U(\xi, \eta) &=& \dbE\Big[ - \cos(\xi) \sin(\tilde \xi) \eta \tilde \eta - \sin(\xi) \int_\dbR\cos(x')\mu(dx')|\eta|^2 +\k |\eta|^2\Big]  \\
 \dis &= & -\dbE[\cos(\xi)\eta] \dbE[\sin(\xi) \eta]  -\dbE[\sin(\xi) |\eta|^2]\dbE[\cos(\xi)]  +\k \dbE[|\eta|^2] \\
\dis &\ge& - \big(\dbE[|\eta|]\big)^2 - \dbE[|\eta|^2] + 2 \dbE[|\eta|^2] \ge 0. 
 \eeaa
 That is, $U$ is displacement monotone. On the other hand, however, by setting $\eta_1 \equiv 0$, we have
 \beaa
 \dis \cM_{DP}^2 U(\xi, 0, \eta_2) &=& \dbE\Big[ \sin(\xi) \cos(\tilde \xi) \eta_2 \tilde \eta_2 +\sin(\xi) \int_\dbR\cos(x')\mu(dx')|\eta_2|^2\Big]  \\
 \dis &= &\dbE[\sin(\xi)\eta_2] \dbE[\cos(\xi) \eta_2]  + \dbE[\sin(\xi) |\eta_2|^2] \dbE[\cos(\xi)].
 \eeaa
Set $\xi \equiv {3\over 4} \pi$ and $\eta_2\equiv 1$, we have 
\beaa
\cM_{DP}^2 U(\xi, 0, \eta_2) = 2 \sin({3\over 4} \pi) \cos({3\over 4} \pi) =-1 <0.
\eeaa
Thus $U$ does not satisfy  \reff{DPmon2}.
\qed

\bs

\no{\bf Proof of \reff{propagation-claim}.} Recall that $u(t,x) = V(t, x, m_t)$ and again we drop the variables $(t, m_t)$. Applying It\^o's formula \reff{Ito} we have
\begin{eqnarray}\label{dV1}
     &&d\partial_{x\mu}V(X_t,\tilde X_t)= \partial_{xx\mu}V\sqrt{2\partial_{\gamma}\cH_V (X_t)}dB_t+\partial_{\tilde{x}x\mu}V\sqrt{2\partial_{\gamma}\cH_V (\tilde{X}_t))}d\tilde{B}_t\nonumber\\
     &&\q +\big[\partial_{tx\mu}V +\partial_{xx\mu}V\partial_{z}\cH_V (X_t)+\partial_{\tilde{x}x\mu}V\partial_{z}\cH_V (\tilde{X}_t)
     +\partial_{xxx\mu}V\partial_{\gamma}\cH_V (X_t)+\partial_{\tilde{x}\tilde{x}x\mu}V\partial_{\gamma}\cH_V (\tilde{X}_t)\big]dt\nonumber\\
     &&\q+\mathbb{E}_{\cF_t\vee \tilde \cF_t}\big[\partial_{x\mu \mu}V(X_t,\tilde{X}_t,\bar{X}_t)\partial_{z}\cH_V (\bar{X}_t)
     +\partial_{\bar{x}x\mu\mu}V(X_t,\tilde{X}_t,\bar{X}_t)\partial_{\gamma}\cH_V (\bar{X}_t)\big]dt.
\end{eqnarray}
Here we dropped the variables $(X_t, \tilde X_t)$ as well inside the derivatives of $V$ when the contexts are clear, and similarly at below we may drop the variables $(x, \tilde x)$ for functions. 

Next,  apply $\pa_{\mu}$ to the master equation \eqref{master0} and recall \reff{propagation-dX2}, we obtain
\bea
\label{dV2}
\left.\ba{lll}
\dis \partial_{\mu t}V(x,\tilde{x})+ \pa_z \cH_V (x) \pa_{x\mu} V+ \pa_\g \cH_V (x) \pa_{xx\mu} V + \pa_\mu F_0 + \pa_{\tilde x\mu} V \pa_z \cH_V (\tilde x) + \pa_\mu V \th_1(\tilde x)  \\
\dis \q + \pa_{\tilde x\tilde x\mu} V \pa_\g \cH_V (\tilde x) + \pa_{\tilde x \mu} V \b_1(\tilde x) + \dbE\Big[\pa_{\mu\mu} V(x,  \bar X_t, \tilde x)\pa_z \cH_V (\bar X_t) + \pa_\mu V(x, \bar X_t) \th_2(\bar X_t, \tilde x) \\
\dis\q + \pa_{\tilde x\mu\mu} V(x,\bar X_t, \tilde x)\pa_\g \cH_V (\bar X_t) + \pa_{\tilde x\mu} V(x, \bar X_t) \b_2(\bar X_t, \tilde x)\Big]=0. 
\ea\right.
\eea
Here we need to explain the notations. First, when taking $\pa_\mu$ on $\mathcal{N}_0V(t,x,\mu)$ in \reff{master0}, we rewrite it  as an integral over $\mu(d\bar x)$ instead of  $\mu(d\tilde x)$. More importantly, for a function $\f(x, \tilde x, \bar x)$ with the variables in natural order $(x, \tilde x, \bar x)$, the index for the derivatives always follow this order, even if we change or switch the variables inside the function. For example, for the terms   $\pa_{\tilde x\mu\mu} V(x,\bar X_t, \tilde x)$ and $\pa_{\tilde x\mu} V(x, \bar X_t)$ at above, the $\pa_{\tilde x}$ refers to the second variable $\bar X_t$, rather than the third variable $\tilde x$, and at below we take the following convention:
\bea
\label{convention}
\pa_{\tilde x} \big[\pa_{\mu\mu} V(x, \bar x, \tilde x)\big] = \pa_{\bar x\mu\mu} V(x, \bar x, \tilde x),\q \pa_{\bar x} \big[\pa_{\mu\mu} V(x, \bar x, \tilde x)\big] = \pa_{\tilde x\mu\mu} V(x, \bar x, \tilde x).
\eea
Moreover, with the above convention, we note that
\bea
\label{commute}
\pa_{\mu\mu} V(x, \tilde x, \bar x) = \pa_{\mu\mu} V(x, \bar x, \tilde x),\q\mbox{and thus}\q \pa_{\bar x\mu\mu} V(x, \tilde x, \bar x) = \pa_{\tilde x\mu\mu} V(x, \bar x, \tilde x).
\eea

Apply $\pa_x$ further to \reff{dV2}, we get
\bea
\label{dV3}
\left.\ba{lll}
\dis \partial_{x\mu t}V(x,\tilde{x}) +\pa_z \cH_V (x) \pa_{xx\mu} V+ \pa_z \cH_V (\tilde x)  \pa_{x\tilde x\mu} V + \pa_\g \cH_V (x) \pa_{xxx\mu} V + \pa_\g \cH_V (\tilde x)\pa_{x\tilde x\tilde x\mu} V \\
\dis \q +  \dbE\Big[\pa_{x\mu\mu} V(x, \bar X_t, \tilde x)\pa_z \cH_V (\bar X_t)  + \pa_{x\tilde x\mu\mu} V(x, \bar X_t, \tilde x)\pa_\g \cH_V (\bar X_t) \Big]\\
\dis \q + [\th_1(x)+ \th_1(\tilde x)] \pa_{x\mu} V +  \b_1(x) \pa_{xx\mu} V +  \b_1(\tilde x) \pa_{x\tilde x\mu} V + \pa_{x\mu} F_0 \\
\dis\q  + \dbE\Big[ \pa_{x\mu} V(x, \bar X_t) \th_2(\bar X_t, \tilde x)+ \pa_{x\tilde x\mu} V(x, \bar X_t) \b_2(\bar X_t, \tilde x)\Big]=0. 
\ea\right.
\eea
Plug this into \reff{dV1} and recall \reff{commute},  we obtain
\begin{eqnarray}\label{dV4}
     &&d\partial_{x\mu}V(X_t,\tilde X_t)= \partial_{xx\mu}V\sqrt{2\partial_{\gamma}\cH_V (X_t)}dB_t+\partial_{\tilde{x}x\mu}V\sqrt{2\partial_{\gamma}\cH_V (\tilde{X}_t))}d\tilde{B}_t\nonumber\\
     &&\q -\Big[[\th_1(X_t)+ \th_1(\tilde X_t)] \pa_{x\mu} V + \b_1(X_t) \pa_{xx\mu} V + \b_1(\tilde X_t) \pa_{x\tilde x\mu} V + \pa_{x\mu} F_0\Big]dt\nonumber\\
     &&\q-\mathbb{E}_{\cF_t\vee \tilde \cF_t}\Big[ \pa_{x\mu} V(X_t, \bar X_t) \th_2(\bar X_t, \tilde X_t)+  \pa_{x\tilde x\mu} V(X_t, \bar X_t) \b_2(\bar X_t, \tilde X_t)\Big]dt.
\end{eqnarray}
We remark that there are lots of cancellations at above, and besides $\pa_{x\mu} F_0$, the remaining terms for $dt$ are exactly those involving $\th_i$ and $\b_i$. 
Now recall \reff{propagation-dX} and rewrite the $\tilde X$ there as $\bar X$, it follows from the standard It\^o's formula that
\beaa
%\label{dV5}
&&\dis d \cM^1_{LL} V(X_t, \d X_t) = \dbE\Big[ \d X_t \widetilde{\d X_t} ~d ( \pa_{x\mu} V(X_t, \tilde X_t) )\Big]\\
&&\dis \q +\dbE\Big[ \pa_{x\mu} V\widetilde{\d X_t} \big[\th_1(X_t) \d X_t + \th_2(X_t, \bar X_t) \overline {\d X_t}~\big] +\pa_{x\mu} V \d X_t\big[\th_1(\tilde X_t) \widetilde{\d X_t} + \th_2(\tilde X_t, \bar X_t) \overline {\d X_t}~\big]\\
&&\dis \q+\partial_{xx\mu}V \big[\b_1(X_t) \d X_t + \b_2(X_t, \bar X_t) \overline {\d X_t}\big]\widetilde{\d X_t} +\partial_{x\tilde x\mu}V \big[\b_1(\tilde X_t)\widetilde{ \d X_t} + \b_2(\tilde X_t, \bar X_t) \overline {\d X_t}\big]\d X_t\Big]dt.
\eeaa
Plug \reff{dV4} into the above, by direct cancellation we have
\beaa
&&\dis {d\over dt} \cM^1_{LL} V( X_t, \d X_t) =  \dbE\Big[ \pa_{x\mu} V(X_t, \tilde X_t) \th_2(X_t, \bar X_t)\widetilde{\d X_t} \overline {\d X_t} +\pa_{x\mu} V(X_t, \tilde X_t) \th_2(\tilde X_t, \bar X_t)  \d X_t\overline {\d X_t}\\
&&\dis \q+\partial_{xx\mu}V (X_t, \tilde X_t)\b_2(X_t, \bar X_t) \overline {\d X_t}\widetilde{\d X_t} +\partial_{x\tilde x\mu}V (X_t, \tilde X_t)\b_2(\tilde X_t, \bar X_t) \overline {\d X_t}\d X_t\\
&&\dis\q  - \big[\pa_{x\mu} F_0(X_t, \tilde X_t)  + \pa_{x\mu} V(X_t, \bar X_t) \th_2(\bar X_t, \tilde X_t) +  \pa_{x\tilde x\mu} V(X_t, \bar X_t) \b_2(\bar X_t, \tilde X_t)\big]\d X_t \widetilde{\d X_t}\Big].
\eeaa
Since $(\tilde X, \widetilde {\d X})$ and $(\bar X, \overline{\d X})$ are independent copies of $(X, \d X)$, we may switch them without changing the expectation. Then we can simplify the above further: 
\beaa
&&\dis {d\over dt} \cM^1_{LL} V(X_t, \d X_t) +  \cM^1_{LL} F_0(X_t, \d X_t)\\
&&\dis = \dbE\Big[ \pa_{x\mu} V(X_t, \tilde X_t) \th_2(X_t, \bar X_t)\widetilde{\d X_t} \overline {\d X_t} +\partial_{xx\mu}V (X_t, \tilde X_t)\b_2(X_t, \bar X_t) \overline {\d X_t}\widetilde{\d X_t} \Big]\\
&&= \dbE\Big\{ \Big[\pa_{x\mu} V(X_t, \tilde X_t)\big[ \pa_{zz}\cH_V (X_t)\pa_{x\mu}V(X_t,\bar X_t)+\pa_{z\gamma}\cH_V (X_t)\pa_{xx\mu}V(X_t,\bar X_t)\big]\\
&&\dis \q+\partial_{xx\mu}V (X_t, \tilde X_t) \big[ \pa_{z\g}\cH_V (X_t)\pa_{x\mu}V(X_t,\bar X_t)+\pa_{\g\gamma}\cH_V (X_t)\pa_{xx\mu}V(X_t,\bar X_t)\big]\Big]\widetilde{\d X_t} \overline {\d X_t} \Big\}.
\eeaa
This implies  \reff{propagation-claim} immediately.
\qed

\bs
\no {\bf Proof of Lemma \ref{lem-Lip}.} We first prove \reff{Lip-claim1}. Let $\Xi=\Xi^\xi$ be the solution to the MFG system \reff{FBSDE}. We see that $\pa_{zz}\cH_V (X_t) = \pa_{zz} H_0(t, \Xi_t)$. Then, applying It\^{o}'s formula \reff{Ito} we have
\beaa
&&\dis d \pa_{zz}\cH_V (X_t) = \Big\{\pa_{tzz} \cH_V   + \pa_{xzz}\cH_V  \pa_z\cH_V  - \pa_{zzz}\cH_V  \pa_x\cH_V -  \pa_{\g zz}\cH_V \Big[\partial_{xx}\cH_V  +2\partial_{xz}\cH_V  Y^{2}_t \\
&&\dis \q +\partial_{zz}\cH_V  |Y^{2}_t|^2+2\big[\partial_{x\gamma}\cH_V  + \partial_{z\gamma}\cH_V Y^{2}_t\big] Z^{2}_t\slash \sqrt{2\partial_{\gamma}\cH_V }+\partial_{\gamma\gamma}\cH_V (Z^{2}_t)^2\slash \big(2\partial_{\gamma}\cH_V \big)\Big]\\
&&\dis \q+ \pa_{xxzz}\cH_V  \pa_\g \cH_V  + {1\over 2} \pa_{zzzz}\cH_V  |Z^1_t|^2 + {1\over 2} \pa_{\g\g zz} \cH_V  |Z^2_t|^2 \\
&&\dis\q +  \pa_{xzzz}\cH_V  \sqrt{2\pa_\g \cH_V } Z^1_t + \pa_{x\g zz}\cH_V  \sqrt{2\pa_\g \cH_V } Z^2_t + \pa_{\g zzz}\cH_V  Z^1_t Z^2_t\Big\}dt \\
&&\dis+ \Big[\pa_{xzz}\cH_V \sqrt{2\pa_\g \cH_V } + \pa_{zzz}\cH_V  Z^1_t + \pa_{\g zz}\cH_V  Z^2_t\Big]dB_t.  
\eeaa
Here we dropped the variable $X_t$ inside $\cH_V $. Note that $Y^1, Y^2, Z^1, Z^2$  and all the involved derivatives of $H_0$ are bounded, then clearly $\pa_{zz} \cH_V $ satisfies \reff{Lip-claim1}. Similarly we can show that $\pa_{z\g}\cH_V $ and $\pa_{\g\g}\cH_V $ also satisfy \reff{Lip-claim1}, thus we prove \reff{Lip-claim1}.

We now prove \reff{Lip-claim2}. Again we shall drop the variables $(X, \tilde X)$ when the contexts are clear, and  denote $ \G^V_2 := \left(\ba{c} \G^V_{2,1} \q \G^V_{2,3}\\ \G^V_{2,2} \q \G^V_{2,4}\ea\right)$. First, by \reff{dV4} and \reff{propagation-dX2} we verify \reff{Lip-claim2} for $U_1$ with
\beaa
-\G^V_{2,1} = [\th_1(X_t)+ \th_1(\tilde X_t)] \pa_{x\mu} V + \b_1(X_t) \pa_{xx\mu} V + \b_1(\tilde X_t) \pa_{x\tilde x\mu} V + \pa_{x\mu} F_0.
\eeaa
Note that $\th_1$ and $\b_1$ are bounded, clearly we have $|\G^V_{2,1}|\le  C\big[|\mathfrak{D}_\mu V|+1\big]$.

Next, similarly to \reff{dV1} we derive from  It\^{o}'s formula \reff{Ito} that
\begin{eqnarray}
\label{dU2-1}
\left.\ba{lll}
 \dis dU_2(X_t,\tilde X_t)= \partial_{xxx\mu}V\sqrt{2\partial_{\gamma}\cH_V (X_t)}dB_t+\partial_{\tilde{x}xx\mu}V\sqrt{2\partial_{\gamma}\cH_V (\tilde{X}_t))}d\tilde{B}_t +\Big[\partial_{txx\mu}V\\
\dis\q  +\partial_{xxx\mu}V\partial_{z}\cH_V (X_t)+\partial_{\tilde{x}xx\mu}V\partial_{z}\cH_V (\tilde{X}_t)
     +\partial_{xxxx\mu}V\partial_{\gamma}\cH_V (X_t)+\partial_{\tilde{x}\tilde{x}xx\mu}V\partial_{\gamma}\cH_V (\tilde{X}_t)\Big]dt\\
\dis \q+\mathbb{E}_{\cF_t\vee \tilde \cF_t}\big[\partial_{xx\mu \mu}V(X_t,\tilde{X}_t,\bar{X}_t)\partial_{z}\cH_V (\bar{X}_t)
     +\partial_{\bar{x}xx\mu\mu}V(X_t,\tilde{X}_t,\bar{X}_t)\partial_{\gamma}\cH_V (\bar{X}_t)\big]dt.
\ea\right.
\end{eqnarray}
Taking $\pa_x$ in \reff{dV3} and recalling that  $\pa_x[\pa_z \cH_V (x)] = \th_1(x)$, $\pa_x[ \pa_\g \cH_V (x)] = \b_1(x)$, we obtain
\bea
\label{dU2-2}
\left.\ba{lll}
\dis \partial_{xx\mu t}V(x,\tilde{x}) +\pa_z \cH_V (x) \pa_{xxx\mu} V+ \pa_z \cH_V (\tilde x)  \pa_{xx\tilde x\mu} V + \pa_\g \cH_V (x) \pa_{xxxx\mu} V + \pa_\g \cH_V (\tilde x)\pa_{xx\tilde x\tilde x\mu} V \\
\dis \q +  \dbE\Big[\pa_{xx\mu\mu} V(x, \bar X_t, \tilde x)\pa_z \cH_V (\bar X_t)  + \pa_{xx\tilde x\mu\mu} V(x, \bar X_t, \tilde x)\pa_\g \cH_V (\bar X_t) \Big]\\
\dis \q + [2\th_1(x)+ \th_1(\tilde x)+  \pa_x\b_1(x)] \pa_{xx\mu} V +  2\b_1(x) \pa_{xxx\mu} V +  \b_1(\tilde x) \pa_{xx\tilde x\mu} V   + \pa_x\th_1(x) \pa_{x\mu} V \\
\dis\q  + \pa_{xx\mu} F_0 + \dbE\Big[ \pa_{xx\mu} V(x, \bar X_t) \th_2(\bar X_t, \tilde x)+ \pa_{xx\tilde x\mu} V(x, \bar X_t) \b_2(\bar X_t, \tilde x)\Big]  =0. 
\ea\right.
\eea
Plug this into \reff{dU2-1} and recall \reff{commute}, we have
\begin{eqnarray}
\label{dU2-3}
\left.\ba{lll}
 \dis dU_2(X_t,\tilde X_t)= \partial_{xxx\mu}V\sqrt{2\partial_{\gamma}\cH_V (X_t)}dB_t+\partial_{\tilde{x}xx\mu}V\sqrt{2\partial_{\gamma}\cH_V (\tilde{X}_t))}d\tilde{B}_t \\
\dis\q -\Big[ [2\th_1(X_t)+ \th_1(\tilde X_t)+\pa_x\b_1(X_t)] \pa_{xx\mu} V + 2 \b_1(X_t) \pa_{xxx\mu} V +  \b_1(\tilde X_t) \pa_{xx\tilde x\mu} V  + \pa_{xx\mu} F_0 \\
\dis\qq  + \pa_x\th_1(X_t) \pa_{x\mu} V+ \dbE_{\cF_t\vee \tilde \cF_t}\big[ \pa_{xx\mu} V(X_t, \bar X_t) \th_2(\bar X_t, \tilde X_t)+ \pa_{xx\tilde x\mu} V(X_t, \bar X_t) \b_2(\bar X_t, \tilde X_t)\big] \Big]dt.
\ea\right.
\end{eqnarray}
Then we see that $U_2$ satisfies \reff{Lip-claim2} with
\beaa
-\G^V_{2,2} =[2\th_1(X_t)+ \th_1(\tilde X_t)+\pa_x\b_1(X_t)] \pa_{xx\mu} V +  \b_1(\tilde X_t) \pa_{xx\tilde x\mu} V  + \pa_x\th_1(X_t) \pa_{x\mu} V  + \pa_{xx\mu} F_0.
\eeaa
Again it is obvious that $|\G^V_{2,2}|\le  C\big[|\mathfrak{D}_\mu V|+1\big]$.

We now verify \reff{Lip-claim2} for $U_3$. Note that 
\begin{eqnarray}
\label{dU3-1}
\left.\ba{lll}
 \dis dU_3(X_t,\tilde X_t)= \partial_{\tilde xxx\mu}V\sqrt{2\partial_{\gamma}\cH_V (X_t)}dB_t+\partial_{\tilde{x}\tilde xx\mu}V\sqrt{2\partial_{\gamma}\cH_V (\tilde{X}_t))}d\tilde{B}_t +\Big[\partial_{t\tilde xx\mu}V\\
\dis\q  +\partial_{\tilde xxx\mu}V\partial_{z}\cH_V (X_t)+\partial_{\tilde{x}\tilde xx\mu}V\partial_{z}\cH_V (\tilde{X}_t)
     +\partial_{\tilde xxxx\mu}V\partial_{\gamma}\cH_V (X_t)+\partial_{\tilde{x}\tilde{x}\tilde xx\mu}V\partial_{\gamma}\cH_V (\tilde{X}_t)\Big]dt\\
\dis \q+\mathbb{E}_{\cF_t\vee \tilde \cF_t}\big[\partial_{\tilde xx\mu \mu}V(X_t,\tilde{X}_t,\bar{X}_t)\partial_{z}\cH_V (\bar{X}_t)
     +\partial_{\bar{x}\tilde xx\mu\mu}V(X_t,\tilde{X}_t,\bar{X}_t)\partial_{\gamma}\cH_V (\bar{X}_t)\big]dt.
\ea\right.
\end{eqnarray}
Take $\pa_{\tilde x}$ in \reff{dV3} and recall the convention \reff{convention}, we have
\bea
\label{dU3-2}
\left.\ba{lll}
\dis \partial_{\tilde xx\mu t}V +\pa_z \cH_V (x) \pa_{\tilde xxx\mu} V+ \pa_z \cH_V (\tilde x)  \pa_{\tilde xx\tilde x\mu} V + \pa_\g \cH_V (x) \pa_{\tilde xxxx\mu} V + \pa_\g \cH_V (\tilde x)\pa_{\tilde x x\tilde x\tilde x\mu} V \\
\dis \q +  \dbE\Big[\pa_{\bar xx\mu\mu} V(x, \bar X_t,, \tilde x)\pa_z \cH_V (\bar X_t)  + \pa_{\bar xx\tilde x\mu\mu} V(x, \bar X_t, \tilde x)\pa_\g \cH_V (\bar X_t) \Big]\\
\dis \q + [\th_1(x)+ 2\th_1(\tilde x)+  \pa_{x} \b_1(\tilde x)] \pa_{\tilde xx\mu} V +  \b_1(x) \pa_{\tilde xxx\mu} V +  2\b_1(\tilde x) \pa_{\tilde xx\tilde x\mu} V + \pa_{\tilde xx\mu} F_0 \\
\dis \q + \pa_{x} \th_1(\tilde x) \pa_{x\mu} V  + \dbE\Big[ \pa_{x\mu} V(x, \bar X_t) \pa_{\tilde x}\th_2(\bar X_t, \tilde x)+ \pa_{x\tilde x\mu} V(x, \bar X_t) \pa_{\tilde x}\b_2(\bar X_t, \tilde x)\Big]=0. 
\ea\right.
\eea
Plug this into \reff{dU3-1} and recall \reff{commute}, we have
\begin{eqnarray}
\label{dU3-3}
\left.\ba{lll}
 \dis dU_3(X_t,\tilde X_t)= \partial_{\tilde xxx\mu}V\sqrt{2\partial_{\gamma}\cH_V (X_t)}dB_t+\partial_{\tilde{x}\tilde xx\mu}V\sqrt{2\partial_{\gamma}\cH_V (\tilde{X}_t))}d\tilde{B}_t \\
\dis\q -\Big[ [\th_1(X_t)+ 2\th_1(\tilde X_t)+\pa_{x}\b_1(\tilde X_t)] \pa_{\tilde xx\mu} V +  \b_1(X_t) \pa_{\tilde xxx\mu} V +  2\b_1(\tilde X_t) \pa_{x\tilde x\tilde x\mu} V  + \pa_{\tilde xx\mu} F_0 \\
\dis\q  + \pa_{x}\th_1(\tilde X_t) \pa_{x\mu} V+ \dbE_{\cF_t\vee \tilde \cF_t}\big[ \pa_{x\mu} V(X_t, \bar X_t) \pa_{\tilde x}\th_2(\bar X_t, \tilde X_t)+ \pa_{x\tilde x\mu} V(X_t, \bar X_t)\pa_{\tilde x} \b_2(\bar X_t, \tilde X_t)\big] \Big]dt.
\ea\right.
\end{eqnarray}
Then we see that $U_3$ satisfies \reff{Lip-claim2} with
\beaa
-\G^V_{2,3} =[\th_1(X_t)+ 2\th_1(\tilde X_t)+\pa_{x}\b_1(\tilde X_t)] \pa_{\tilde xx\mu} V +  \b_1(X_t) \pa_{\tilde xxx\mu} V +  \pa_{x}\th_1(\tilde X_t) \pa_{x\mu} V  + \pa_{\tilde xx\mu} F_0.
\eeaa
Again it is obvious that $|\G^V_{2,3}|\le  C\big[|\mathfrak{D}_\mu V|+1\big]$.

Finally, for $U_4$ we have 
\begin{eqnarray}
\label{dU4-1}
\left.\ba{lll}
 \dis\!\!\!\! dU_4(X_t,\tilde X_t)= \partial_{x\tilde xxx\mu}V\sqrt{2\partial_{\gamma}\cH_V (X_t)}dB_t+\partial_{x\tilde{x}\tilde xx\mu}V\sqrt{2\partial_{\gamma}\cH_V (\tilde{X}_t))}d\tilde{B}_t +\Big[\partial_{xt\tilde xx\mu}V\\
\dis\!\!\!\!  +\partial_{x\tilde xxx\mu}V\partial_{z}\cH_V (X_t)+\partial_{x\tilde{x}\tilde xx\mu}V\partial_{z}\cH_V (\tilde{X}_t)
     +\partial_{x\tilde xxxx\mu}V\partial_{\gamma}\cH_V (X_t)+\partial_{x\tilde{x}\tilde{x}\tilde xx\mu}V\partial_{\gamma}\cH_V (\tilde{X}_t)\Big]dt\\
\dis\!\!\!\! +\mathbb{E}_{\cF_t\vee \tilde \cF_t}\big[\partial_{x\tilde xx\mu \mu}V(X_t,\tilde{X}_t,\bar{X}_t)\partial_{z}\cH_V (\bar{X}_t)
     +\partial_{\bar{x}x\tilde xx\mu\mu}V(X_t,\tilde{X}_t,\bar{X}_t)\partial_{\gamma}\cH_V (\bar{X}_t)\big]dt.
\ea\right.
\end{eqnarray}
Take $\pa_x$ in \reff{dU3-2}, we have
\bea
\label{dU4-2}
\left.\ba{lll}
\dis \partial_{x\tilde xx\mu t}V +\pa_z \cH_V (x) \pa_{x\tilde xxx\mu} V+ \pa_z \cH_V (\tilde x)  \pa_{x\tilde xx\tilde x\mu} V + \pa_\g \cH_V (x) \pa_{x\tilde xxxx\mu} V + \pa_\g \cH_V (\tilde x)\pa_{x\tilde x x\tilde x\tilde x\mu} V \\
\dis \q +  \dbE\Big[\pa_{x\bar xx\mu\mu} V(x, \bar X_t, \tilde x)\pa_z \cH_V (\bar X_t)  + \pa_{x\bar xx\tilde x\mu\mu} V(x, \bar X_t, \tilde x)\pa_\g \cH_V (\bar X_t) \Big]\\
\dis \q + [2\th_1(x)+ 2\th_1(\tilde x) +  \pa_x\b_1(x) +  \pa_{ x} \b_1(\tilde x)] \pa_{x\tilde xx\mu} V +  2\b_1(x) \pa_{x\tilde xxx\mu} V +  2\b_1(\tilde x) \pa_{x\tilde xx\tilde x\mu} V  \\
\dis \q + \pa_x\th_1(x) \pa_{\tilde xx\mu} V + \pa_{x} \th_1(\tilde x) \pa_{xx\mu} V  + \pa_{x\tilde xx\mu} F_0 \\
\dis\q + \dbE\Big[ \pa_{xx\mu} V(x, \bar X_t) \pa_{\tilde x}\th_2(\bar X_t, \tilde x)+ \pa_{xx\tilde x\mu} V(x, \bar X_t) \pa_{\tilde x}\b_2(\bar X_t, \tilde x)\Big]  =0 
\ea\right.
\eea
Plug this into \reff{dU4-1} and recall \reff{commute}, we have
\begin{eqnarray}
\label{dU4-3}
\left.\ba{lll}
 \dis dU_4(X_t,\tilde X_t)= \partial_{x\tilde xxx\mu}V\sqrt{2\partial_{\gamma}\cH_V (X_t)}dB_t+\partial_{x\tilde{x}\tilde xx\mu}V\sqrt{2\partial_{\gamma}\cH_V (\tilde{X}_t))}d\tilde{B}_t  \\
\dis\q -\Big[ [2\th_1(X_t)+ 2\th_1(\tilde X_t) + \pa_x \b_1(X_t)+\pa_{x}\b_1(\tilde X_t)] \pa_{x\tilde xx\mu} V +  2\b_1(X_t) \pa_{x\tilde xxx\mu} V    \\
\dis\q +  2\b_1(\tilde X_t) \pa_{xx\tilde x\tilde x\mu} V+ \pa_{x}\th_1(X_t) \pa_{\tilde xx\mu} V  + \pa_{x}\th_1(\tilde X_t) \pa_{xx\mu} V+ \pa_{x\tilde xx\mu} F_0\\
\dis\q + \dbE_{\cF_t\vee \tilde \cF_t}\big[ \pa_{xx\mu} V(X_t, \bar X_t) \pa_{\tilde x}\th_2(\bar X_t, \tilde X_t)+ \pa_{xx\tilde x\mu} V(X_t, \bar X_t)\pa_{\tilde x} \b_2(\bar X_t, \tilde X_t)\big] \Big]dt.
\ea\right.
\end{eqnarray}
Then we see that $U_4$ satisfies \reff{Lip-claim2} with
\beaa
-\G^V_{2,4} &=& [2\th_1(X_t)+ 2\th_1(\tilde X_t) + \pa_x \b_1(X_t)+\pa_{x}\b_1(\tilde X_t)] \pa_{x\tilde xx\mu} V \\
&&+ \pa_{x}\th_1(X_t) \pa_{\tilde xx\mu} V  + \pa_{x}\th_1(\tilde X_t) \pa_{xx\mu} V+ \pa_{x\tilde xx\mu} F_0.
\eeaa
Once again it is obvious that $|\G^V_{2,4}|\le  C\big[|\mathfrak{D}_\mu V|+1\big]$. The proof is complete now.
\qed

\bs

\no{\bf Section \ref{sect-localproof} Step 3: representation of $\pa_{\mu\mu} V$.}

We now analyze $\pa_{\mu\mu} V$ by using the bootstrap arguments, by assuming all the derivatives in \reff{xreg} are given.  We first note that, similar to the reason that we consider $(Y^{\xi,1}, Y^{\xi,2})$ and that $Y^{\xi,0}$ is decoupled in \reff{FBSDE},  it is more convenient to construct a representation for $\pa_{x\mu\mu} V$. This will lead further to a representation for  $\pa_{\mu\mu} V$, which we skip here. 

Next, to avoid the rather heavy notations, we focus only on the following special case, and we emphasize that the same arguments work for the general case \reff{master}: letting $\hat \xi\in \dbL^2(\cF_t, \mu)$,
\bea
\label{master2}
\pa_t V(t,x,\mu) + H(\pa_{xx} V(t,x,\mu)) + \dbE\big[\pa_{\tilde x\mu} V(t,x,\mu, \hat \xi) \pa_\g H(\pa_{xx}V(t, \hat \xi, \mu))\big]=0.
\eea
 Again dropping the variables $(t,\mu)$, similar to \reff{dV3} we can compute:
\bea
\label{master2mu}
\left.\ba{lll}
\dis \pa_{x\mu t}V(x,\tilde x) + \pa_\g \cH_V(x) \pa_{xxx\mu} V(x,\tilde x) + \pa_x[ \pa_\g \cH_V(x)] \pa_{xx\mu} V(x,\tilde x) \\
\dis\q +  \pa_\g \cH_V(\tilde x) \pa_{\tilde x\tilde x x\mu} V(x,\tilde x)+ \pa_{\tilde x x\mu} V(x, \tilde x) \pa_{\g\g}\cH_V(\tilde x)\pa_{xxx} V(\tilde x) \\
\dis \q + \dbE\big[\pa_{\bar x x\mu\mu} V(x,\tilde x, \hat \xi) \pa_\g \cH_V(\hat \xi)+\pa_{\tilde x x\mu} V(x,\hat \xi) \pa_{\g\g} \cH_V(\hat \xi)\pa_{xx\mu} V(\hat \xi, \tilde x)\big] =0.
\ea\right.
\eea
Here we used the notation convention \reff{convention} and \reff{commute} for the term $\pa_{\bar x x\mu\mu} V(x,\tilde x, \hat \xi)$. 
We remark that, the last term above involves $\pa_{\tilde x x\mu} V\pa_{xx\mu} V$, and thus is not linear in terms of $\pa_{\mu} V$ and its derivatives, and that's why we made the serious efforts to establish the representation formulae in Subsection \ref{sect-point}.

The situation is different for higher order derivatives. Differentiate \reff{master2mu} further in $\mu$ formally and denote $U(x,\tilde x, \bar x) := \pa_{x\mu\mu} V(x, \tilde x, \bar x)$. Then we have:
\bea
\label{master2mumu}
\left.\ba{lll}
\dis \pa_{t}U(x,\tilde x, \bar x) + \pa_\g \cH_V(x) \pa_{xx} U(x,\tilde x, \bar x) +  \pa_\g \cH(\tilde x) \pa_{\tilde x\tilde x} U(x,\tilde x, \bar x) +\pa_\g \cH_V(\bar x) \pa_{\bar x\bar x} U(x,\tilde x, \bar x)\\
\dis +  \pa_{\g\g} \cH_V(\tilde x)\pa_{xxx} V(\tilde x) \pa_{\tilde x} U(x, \tilde x, \bar x) + \pa_{\g\g} \cH_V(\bar x)\pa_{xxx} V(\bar x)\pa_{\bar x} U(x, \tilde x,\bar x) \\
\dis+ \dbE\Big[\pa_\g \cH_V(\hat \xi)\pa_{\hat x\mu} U(x, \tilde x,  \bar x, \hat \xi)  + \pa_{\g\g} \cH_V(\hat \xi)\pa_{xx\mu} V(\hat\xi, \bar x) \pa_{\bar x} U(x,\tilde x, \hat \xi)\\
\dis  + \pa_{\g\g} \cH_V(\hat \xi))\pa_{xx\mu} V(\hat \xi, \tilde x)\pa_{\tilde x} U(x,\hat \xi, \bar x) +\pa_{\tilde x x\mu} V(x,\hat \xi) \pa_{\g\g} \cH_V(\hat \xi)\pa_x U(\hat \xi, \tilde x, \bar x)\Big] + \G(x, \tilde x, \bar x) =0,
\ea\right.
\eea
where $\G$ involves only $\pa^{(k)}_x V$ and $\pa^{(k)}_x \pa^{(l)}_{\tilde x}\pa_\mu V$ in \reff{xreg}, as well as $H$ and its derivatives, and thus can be viewed as a given function.  Then this is a linear master equation for $U$. 

To introduce a representation formula for $U(0, x, \mu, \tilde x, \bar x)$, denote $X := X^\xi$ and note that
\beaa
X_t =\xi+ \int_0^t \sqrt{2\pa_\g H(\pa_{xx} V(s, X_s, \cL_{X_s}))} dB_s.
\eeaa
Denote 
\bea
\label{mumuU}
\left.\ba{c}
\dis Y_t := U(X_t, \tilde X_t, \bar X_t),\q Z_t := \pa_x U(X_t, \tilde X_t, \bar X_t)\sqrt{2\pa_\g \cH_V(X_t)},\\
\dis \tilde Z_t := \pa_{\tilde x} U(X_t, \tilde X_t, \bar X_t)\sqrt{2\pa_\g \cH_V(\tilde X_t)},\q \bar Z_t := \pa_{\bar x} U(X_t, \tilde X_t, \bar X_t)\sqrt{2\pa_\g \cH_V(\bar X_t)}.
\ea\right.
\eea
Then by applying It\^{o}'s formula \reff{Ito} we have, letting $\G_t$ denote a generic bounded process,
\bea
\label{mumuBSDE}
\left.\ba{lll}
\dis d Y_t =  Z_t dB_t + \tilde Z_t d\tilde B_t + \bar Z_t d\bar B_t     - \dbE_{\cF_t\vee \tilde \cF_t \vee \bar \cF_t}\Big[ \G_t \tilde Z_t + \G_t\bar Z_t + \G_t+\G_t \wh{\bar Z}_t   + \G_t \wh{\tilde Z}_t+\G_t\wh Z_t\Big] dt.
\ea\right.
\eea
Here $\wh Z_t := \pa_x U(\hat X_t, \tilde X_t, \bar Z_t) \sqrt{2\pa_\g \cH_V(\hat X_t)}$, which is a conditional independent copy of $Z_t$, conditional on $\tilde \cF_t\vee \bar \cF_t$. The terms $\wh{\tilde Z}$ and $\wh{\bar Z}$ are defined similarly. Note that \reff{mumuBSDE} is a McKean-Valsov BSDE, rather than a coupled FBSDE, with  terminal condition $Y_T = \pa_{x\mu\mu} G(X_T, \cL_{X_T}, \tilde X_T, \bar X_T)$. Thus it is clearly well-posed. We remark that, although \reff{mumuBSDE} was derived formally by assuming $U$ is smooth, once we obtain the representation formula, one can easily verify that
\bea
\label{mumurep}
Y_0 = U(0,\xi, \mu, \tilde \xi, \bar \xi) = \pa_{x\mu\mu} V(0, \xi, \mu, \tilde \xi, \bar \xi),~\mbox{a.s.}
\eea

Next, denote $\bm := \cL_X$, and introduce
\bea
\label{mumuUx}
\left.\ba{c}
\dis X^x_t = x+ \int_0^t \sqrt{2\pa_\g H(\pa_{xx} V(s, X^x_s, m_s))} dB_s;\\
\dis Y^x_t := U(X^x_t, \tilde X_t, \bar X_t),\q Z^x_t := \pa_x U(X^x_t, \tilde X_t, \bar X_t)\sqrt{2\pa_\g \cH_V(X^x_t)},\\
\dis \tilde Z^x_t := \pa_{\tilde x} U(X^x_t, \tilde X_t, \bar X_t)\sqrt{2\pa_\g \cH_V(\tilde X_t)},\q \bar Z^x_t := \pa_{\bar x} U(X^x_t, \tilde X_t, \bar X_t)\sqrt{2\pa_\g \cH_V(\bar X_t)}.
\ea\right.
\eea
Note that $\tilde X$ and $\bar X$ still have initial condition $\tilde \xi, \bar \xi$. Then, applying It\^{o}'s formula again we have
\bea
\label{mumuBSDEx}
\left.\ba{lll}
\dis d Y^x_t =  Z^x_t dB_t + \tilde Z^x_t d\tilde B_t + \bar Z^x_t d\bar B_t   \\
\dis \qq  - \dbE_{\cF_t\vee \tilde \cF_t \vee \bar \cF_t}\Big[ \G_t \tilde Z^x_t + \G_t\bar Z^x_t + \G_t+\G_t \wh{\bar Z^x_t}   + \G_t \wh{\tilde Z^x_t}+\G_t\wh Z_t\Big] dt.
\ea\right.
\eea
We emphasize that the last term above is $\wh Z_t$, rather than $\wh Z^x_t$. Thus we obtain a representation:
\bea
\label{mumurepx}
Y^x_0 = U(0,x, \mu, \tilde \xi, \bar \xi) = \pa_{x\mu\mu} V(0, x, \mu, \tilde \xi, \bar \xi),~\mbox{a.s.}
\eea

Similarly we may obtain  a.s. representation formulae for $U(0, \xi, \mu, \tilde x, \bar \xi)$ and $U(0, \xi, \mu, \tilde \xi, \bar x)$. Moreover, we may continue the arguments and obtain representation formulae for $U(0, x, \mu, \tilde x, \bar \xi)$, $U(0, x, \mu, \tilde \xi, \bar x)$, $U(0, \xi, \mu, \tilde x, \bar x)$, and finally for $U(0, x, \mu, \tilde x, \bar x)$. We skip the details. 
\qed

\end{document}